\documentclass[10pt]{amsart}
\pdfoutput=1
\usepackage{amsmath, amsthm, amssymb,slashed,stmaryrd}
\usepackage{ifpdf}
\usepackage[pdftex]{graphicx}
\usepackage{tikz}
\usetikzlibrary{matrix,arrows,calc}

\usepackage[pdftex,plainpages=false,hypertexnames=false,pdfpagelabels]{hyperref}
\UseRawInputEncoding

 \theoremstyle{plain}
 \newtheorem{thm}{Theorem}[section]
  
    \newtheorem{goal}{Goal}[section]

  \newtheorem{quest}{Question}[section]
 \newtheorem{cor}[thm]{Corollary}
 \newtheorem{lem}[thm]{Lemma}

 \newtheorem{prop}[thm]{Proposition}
 \theoremstyle{definition}
 \newtheorem{defn}[thm]{Definition}
 \newtheorem{notation}[thm]{Notation}

 \newtheorem{ex}[thm]{Example}
 \theoremstyle{remark}
 \newtheorem{rmk}[thm]{Remark}

\def\beq{\begin{eqnarray}}
\def\eeq{\end{eqnarray}}

\DeclareSymbolFont{bbold}{U}{bbold}{m}{n}
\DeclareSymbolFontAlphabet{\mathbbold}{bbold}

 \newcommand{\bp}{\begin{proof}[Proof]}
 \newcommand{\ep}{\end{proof}}

\def\Hess{{\rm Hess}}
\def\HessS{{\rm Hess}(\mathcal{S}_\sigma)}

\def\x{{\bf x}}
\def\X{{\mathcal{X}}}
\def\Y{{\mathcal{Y}}}
\def\L{{\mathcal{L}}}

\def\pt{{\rm pt}}

\def\cs{{\rm c}}
\def\Wit{{{\rm Wit}}}
\def\ch{{\rm ch}}
\def\vol{{\rm vol}}
\def\sdet{{\rm sdet}}

\def\MQ{{\rm MQ}}
\def\Lat{{\sf Lat}}
\def\Rot{{\rm Rot}}

\def\MF{{\rm MF}}

\def\H{{\rm H}}
\def\HP{{\rm HP}}

\def\KO{{\rm KO}}

\def\SL{{\rm SL}}
\def\GL{{\rm GL}}

\def\Map{{\sf Map}}

\def\ev{{\rm ev}}
\def\odd{{\rm odd}}
\def\cl{{\rm cl}}

\newcommand{\sq}{\mathord{/\!\!/}}

\def\SMQ{\mathcal{S}_{\rm MQ}}

\def\R{{\mathbb{R}}}

\def\id{{{\rm id}}}

\def\K{{\rm {K}}}
\def\C{{\mathbb{C}}}

\def\Q{{\mathbb{Q}}}

\def\Z{{\mathbb{Z}}}
\def\E{{\mathbb{E}}} 

\def\pt{{\rm pt}}

\def\End{{\rm End}}

\def\TMF{{\rm TMF}}

\def\Fer{{\sf Fer}}
\def\Bos{{\sf Bos}}
\def\BosT{\Bos(T(M/B))}

\def\path{ {{\mathfrak{p}}_0}}

\def\Iso{ {{\sf Iso}}}

\newcommand{\op}{{\sf{op}}}

\begin{document}

\title[Field theories and the elliptic index theorem]{Supersymmetric field theories and the elliptic index theorem with complex coefficients}

\def\obfuscate#1#2{\rlap{\hphantom{#2}@#1}#2\hphantom{@#1}}
\author{Daniel Berwick-Evans}
\address{Department of Mathematics, University of Illinois at Urbana-Champaign}
\email{\obfuscate{illinois.edu}{danbe}}

\begin{abstract}
We present a cocycle model for elliptic cohomology with complex coefficients in which methods from 2-dimensional quantum field theory can be used to rigorously construct cocycles. For example, quantizing a theory of vector bundle-valued fermions yields a cocycle representative of the elliptic Thom class. This constructs the complexified string orientation of elliptic cohomology, which determines a pushfoward for families of rational string manifolds. A second pushforward is constructed from quantizing a supersymmetric $\sigma$-model. These two pushforwards agree, giving a precise physical interpretation for the elliptic index theorem with complex coefficients. This both refines and supplies further evidence for the long-conjectured relationship between elliptic cohomology and 2-dimensional quantum field theory. Analogous methods in supersymmetric mechanics recover path integral constructions of the Mathai--Quillen Thom form in complexified $\KO$-theory and a cocycle representative of the $\hat{A}$-class for a family of oriented manifolds. 
\end{abstract}

\maketitle

\setcounter{tocdepth}{1}
\tableofcontents

\section{Introduction} \label{sec:intro}

The last 35 years have seen rich cross-fertilization between geometry, topology, and quantum field theory. A key aspect is rooted in the connection between supersymmetric mechanics and the Atiyah--Singer index theorem~\cite{susymorse,Alvarez,GetzlerIndex}. Witten's work from the late `80s~\cite{Witten_Elliptic,WittenDirac} points towards a generalization of the index theorem associated with the analysis of 2-dimensional quantum field theory on the one hand and the algebraic topology of elliptic cohomology on the other---specifically, the string orientation of topological modular forms (TMF)~\cite{AHS,AHR}. Such an index theorem would probe subtle invariants of these field theories, mirroring the known intricacies in TMF~\cite{HopkinsTMF}. In so doing, it would also impart geometric meaning to powerful homotopy invariants. 

The goal of this paper is to provide solid mathematical footing for some of the basic ideas from physics together with a dictionary that translates to the corresponding structures in algebraic topology. 
Specifically, we put the geometry of perturbative 2-dimensional supersymmetric field theories in direct contact with complexified~$\TMF$. The modifiers \emph{perturbative} and \emph{complexified} simplify the physics and algebraic topology, respectively. This allows us to carry out a baby case of the program initiated by Segal~\cite{Segal_Elliptic,SegalCFT} and Stolz--Teichner~\cite{ST04,ST11} which we paraphrase as follows. Techniques in 2-dimensional supersymmetric field theories construct analytic and topological pushforwards in~$\TMF\otimes \C$, and there is an index theorem equating the two. 

To tie in with the standard index-theoretic story and emphasize structural aspects, we present our results in terms of a hierarchy that on the physics side corresponds to dimensions $d=0,1,2$ and on the algebraic topology side corresponds to de~Rham cohomology, complexified $\K$-theory, and complexified $\TMF$.  
A brief (though imprecise) summary of these results is as follows.


\begin{enumerate}
\item[I:] Let $\Phi$ be the space of fields for the $\mathcal{N}=1$ supersymmetric $\sigma$-model in dimension $d\in \{0,1,2\}$ (see Example~\ref{ex:sigma}). The global minima of the classical action functional yields a finite-dimensional and smooth subspace~$\Phi_0\subset \Phi$. Functions on~$\Phi_0$ carry an action by a groupoid of symmetries of the classical field theory, and this structure furnishes a cocycle model for de~Rham cohomology when $d=0$, complexified $\K$-theory when $d=1$, and complexified~$\TMF$ when $d=2$. See~\S\ref{sec:I} for details.
\item[II:] The Hessian of the classical action of the supersymmetric $\sigma$-model (see Example~\ref{ex:sigmahess}) defines a family of operators parametrized by~$\Phi_0$. The regularized determinant of this family yields a function on $\Phi_0$ that is a cocycle representative of the Riemann--Roch factor for the analytic pushforward in the relevant cohomology theory: this factor is~1 when $d=0$, the $\hat{A}$-class when $d=1$, and the Witten class when $d=2$. See~\S\ref{sec:II} for details.
\item[III:] Applying similar methods to the Mathai--Quillen classical field theory (see Example~\ref{ex:MQ}) yields a cocycle refinement of the Thom class in the relevant cohomology theory over~$\C$. This defines a topological pushforward. When applied to trivial bundles, it also determines a cocycle representative of the suspension class and the Thom isomorphism becomes the suspension isomorphism. This reveals a candidate origin for an infinite loop space structure on a space of 2-dimensional supersymmetric field theories. See~\S\ref{sec:III} for details.
\item[IV:] The analytic and topological pushforwards constructed using these techniques agree at the cocycle level. We remark that this is a trivial consequence of formulas for the pushforwards in II and III. See~\S\ref{sec:IV} for details. 
\end{enumerate}

In light of these results, the contribution of this paper can be summarized as having two pieces. First, I-IV gives a robust and precise link in the long-sought connection between elliptic cohomology and 2-dimensional field theories, providing further supporting evidence for an elliptic index theorem. 
 Second, the constructions going into I-IV highlight the objects in physics that need to be clarified mathematically if one is to understand the full depth of the connection between elliptic cohomology and 2-dimensional field theories. The remainder of this introduction examines this second aspect in more detail.


First we observe an obvious shortcoming of the results above: although $\TMF\otimes \C$ encodes the often mysterious modularity properties of certain classes in elliptic cohomology (e.g., the Witten genus), working over~$\C$ obliterates all the rich torsion in~$\TMF$. Ultimately one would hope to refine I-IV to statements over~$\Z$. The paradigm of \emph{extended} field theories (analogous to those appearing in the cobordism hypothesis~\cite{Lurie_cob}) illuminates a path toward such a refinement. Roughly speaking, the results of this paper concern 0-extended 2-dimensional field theories, meaning the constructions take place over a moduli space of closed 2-dimensional (super) manifolds. Some features of a 1-extended refinement were investigated by Pokman Cheung~\cite{PokmanPhD}, who constructed a geometric model for elliptic cohomology at the Tate curve (over~$\Z$) in terms of representations of a certain bordism 1-category. This follows a sketch by Stolz and Teichner~\cite{ST04} who in turn were inspired by Segal~\cite{SegalCFT}. 
A 2-extended generalization is expected to construct~TMF as representations of some bordism 2-category~\cite{ST04,ST11}. However, the very \emph{definition} of the appropriate category of 2-extended field theories remains unclear, and basic features of the physics behind the analytic and topological pushforwards have yet to be understood mathematically (see~\S\ref{sec:context} for a review).
Given these long-standing challenges, an incremental approach to the full elliptic index theorem for~TMF seems prudent.

With the roadmap of extending down in mind, the main goal of this paper is to pin down the 0-extended case. The homotopical counterpart turns out to be the index theorem for~$\TMF\otimes \C$. Pushforwards in~$\TMF\otimes \C$ are built by quantizing specific classical field theories. More sophisticated quantization methods applied to \emph{the same} classical field theories will refine these pushforwards, e.g., one can study higher-categorical refinements involving boundary conditions. These specific examples will help guide the development of a framework to capture the relevant 1- and 2-extended generalizations.

Indeed, by now there are at least~3 distinct languages in which one might try to connect elliptic cohomology with ideas from physics, all of which have offered different insights into the problem. Perhaps the most well-developed approach is the Stolz--Teichner program; the state-of-the art is currently a 1-extended super bordism category~\cite{ST11}, and there are suggestions of how to proceed to the 2-extended case~\cite{ST04}. Conformal nets offer a second possibility~\cite{chris_andre_TMF,CN1,CN2,CN3,CN4}. Factorization algebras suggest a third route~\cite{costello_WG1,costello_WG2}. The methods we use below were inspired by attempting to combine desirable aspects of these languages. 

In a different direction, an \emph{equivariant} refinement in the 0-extended case that incorporates gauge fields \cite{BET0,BET1} leads directly to a cocycle model for Grojnowski's complex analytic equivariant elliptic cohomology~\cite{Grojnowski}. This is a highly nontrivial extension of~$\TMF\otimes \C$, and its natural emergence suggests that the geometry encountered in the 0-extended setup is significantly deeper than just a cocycle model for~$\TMF\otimes \C$. Among other insights, the equivariant refinement clarifies the physical origins of elliptic (formal) group laws in terms of $U(1)$-gauge theories. These group laws lie at the heart of the homotopical machine that constructs elliptic cohomology theories, and will need to feature prominently if one is to identify a candidate geometric model with~TMF. 

Further mathematical development of a geometric model for elliptic cohomology holds great promise. For example, recent progress by physicists suggests that a robust understanding of the connection between equivariant~TMF and supersymmetric field theory will yield generalizations of Donaldson invariants of 4-manifolds~\cite{VafaTMF}. In a different direction, a construction of certain field theories looks to explain features of the image of topological modular forms in integral modular forms, e.g., why $\Delta$ is not in the image but~$24\Delta$ is~\cite{DavideTheo}. This paper is a first step towards a precise dictionary connecting these and other ideas from physics with structures in homotopy theory.

\subsection*{Organization of the paper}

The next section overviews the intuition behind the anticipated elliptic index theorem while also contextualizing the contribution made by this paper. We present our results in~\S\ref{sec:results} in a manner that emphasizes the hierarchy that on the physics side corresponds to dimensions~$d=0$,~$1$, and~$2$ and on the algebraic topology side corresponds to de~Rham cohomology, $\K$-theory, and $\TMF$. The basic ideas going into the constructions are uniform for all~3 cases. The remaining sections address each case in turn, with de~Rham cohomology in~\S\ref{sec:01}, complexified $\K$-theory in~\S\ref{sec:KO} and complexified $\TMF$ in~\S\ref{sec:TMF}.

\subsection*{Conventions}
The main constructions in this paper take place within the category of supermanifolds with structure sheaves defined over $\C$; see~\S\ref{appen:super} for a brief review. This is the category of $cs$-manifolds in the terminology of~\cite{DM}. When unspecified, objects are implicitly ``super": vector space means super vector space, algebra means super algebra, Lie group means super Lie group, the tensor product~$\otimes$ is the $\Z/2$-graded one, etc. Throughout, $S$ will denote a (test) supermanifold and~$M$ an ordinary smooth manifold. In the usual abuse, $M$ will also denote the supermanifold determined by $M$ with its sheaf of complex-valued smooth functions. For a super Lie group $G$ acting on a supermanifold~$N$, we use the notation~$N\sq G$ for the Lie groupoid quotient, and $[N\sq G]$ for the stack underlying this Lie groupoid. For free actions, we use the notation $N/G$ for the quotient in supermanifolds. When working with algebras of functions on supermanifolds or sections of vector bundles, we use the projective tensor product of Fr\'echet spaces. This is a completion of the algebraic tensor product satisfying
\beq
C^\infty(S)\otimes C^\infty(T)\cong C^\infty(S\times T)\cong C^\infty(S;C^\infty(T))\label{eq:HST}
\eeq
for supermanifolds $S$ and $T$; we refer to~\cite[\S6.1]{HST} for an overview. 

\subsection*{Acknowledgements}

It is a pleasure to thank Matt Ando, Kevin Costello, Ryan Grady, Owen Gwilliam, Theo Johnson-Freyd, Eugene Lerman, Charles Rezk, and Arnav Tripathy for helpful discussions and correspondence. The deep influence of Stephan Stolz and Peter Teichner should be clear, and I thank them for their support and encouragement.

\section{Context and motivation}\label{sec:context}

This section provides mathematical and physical context for the core ideas that go into the constructions outlined in the introduction. Broadly speaking, we look to generalize the well-known physics explanation of the Atiyah--Singer index theorem~\cite{susymorse,Alvarez} via supersymmetric quantum mechanics (in dimension~$d=1$) to elliptic index theory and supersymmetric quantum field theory (where $d=2$).

\subsection{Towards an index theorem in elliptic cohomology}

The \emph{$\hat{A}$-genus} is a $\Q$-valued cobordism invariant of oriented manifolds that is $\Z$-valued on spin manifolds. There is a families version of the $\hat{A}$-genus that takes values in real K-theory,~$\KO$. Let $\pi\colon M\to B$ by a family of spin manifolds\footnote{For simplicity, in this introductory discussion we assume~$M$ is compact.} with fiber dimension~$n$. The \emph{analytic pushforward} is the map~$\pi_!^{\rm an}\colon \KO^\bullet(M)\to \KO^{\bullet-n}(B)$ gotten by taking a suitable fiberwise index of the family of Dirac operators associated to the spin family $\pi\colon M\to B$. When $B=\pt$ and $M$ is $4k$-dimensional, $\pi_!^{\rm an}(1)$ agrees with the $\hat{A}$-genus of $M$. 

One can construct a second pushforward in $\KO$-theory using the Thom isomorphism. A choice of embedding $M\hookrightarrow \R^N$ determines a fiberwise embedding $i\colon M\hookrightarrow B\times \R^N$. Let~$\nu$ be the $(N-n)$-dimensional vector bundle on $M$ that is the normal bundle for~$i$. Then the \emph{topological pushforward} $\pi_!^{\rm top}\colon \KO^\bullet(M)\to \KO^{\bullet-n}(B)$ sits in the commutative diagram
\beq
\begin{tikzpicture}[baseline=(basepoint)];
\node (A) at (0,0) {$\KO^\bullet(M)$};
\node (B) at (6,0) {$\KO_{\cs}^{\bullet+N-n}(B\times \R^N)$};
\node (C) at (2.5,-1.5) {$\KO^{\bullet - n}(B)$};
\draw[->] (A) to node [above=1pt] {$\alpha_\nu\smallsmile$} (B);
\draw[->] (A) to node [left=.1in] {$\pi_!^{\rm top}$} (C);
\draw[->] (B) to node [right=.1in] {$\Sigma^{-N}$} (C);
\path (0,-.75) coordinate (basepoint);
\end{tikzpicture}\nonumber
\eeq
where $\Sigma^{-N}$ is the inverse to the suspension isomorphism, $\alpha_\nu$ is the $\KO$-theory Thom class (that depends on the spin structure on $\pi\colon M\to B$), and $\KO_{\cs}$ denotes $\KO$-theory with compact support. The top horizontal arrow is the Thom isomorphism followed by extension by zero for the inclusion~$\nu\hookrightarrow B\times \R^N$ using a tubular neighborhood for the embedding. The index theorem compares these two pushforwards. 

\begin{thm}[\cite{AtiyahSingerfamII}] \label{indexthm}The analytic and topological pushforwards agree, $\pi_!^{\rm top}=\pi_!^{\rm an}$. \end{thm}

\begin{rmk}\label{rmk:Thom}
The $\KO$-Thom class can be constructed completely homotopy theoretically, which is why the associated pushforward is called \emph{topological}. However, geometric constructions can also be quite insightful. For example, a cocycle representative of the Thom class can be described in terms of Clifford modules~\cite{ABS}, and the 8-periodicity of the Clifford algebras sheds light on the geometric origin of the 8-fold periodicity of~$\KO$. 
\end{rmk} 

The index theorem has an expected analog in elliptic cohomology~\cite{WittenDirac,Segal_Elliptic,ST04}. The starting point is the \emph{Witten genus}, a $\Q\llbracket q\rrbracket$-valued cobordism invariant of oriented manifolds. This genus takes values in integral modular forms when the manifold has a \emph{string structure}, meaning the manifold admits a spin structure and the fractional first Pontryagin class vanishes. Just as $\KO$-theory is the natural home for the families $\hat{A}$-genus, Witten suggested that the families Witten genus take values in elliptic cohomology theories via a (not yet constructed) analytic pushforward~\cite{Witten_Elliptic}. To explain the modularity of the Witten genus, Hopkins and collaborators assembled the elliptic cohomology theories attached to each elliptic curve into a global object over the moduli stack of elliptic curves. The result is the universal elliptic cohomology theory of \emph{topological modular forms}~\cite{HopkinsTMF,HopkinsICM94,GoessTMF,Lurie_Elliptic,Lurie_Ell1,Lurie_Ell2}. This cohomology theory has a families Witten genus via a topological pushforward~\cite{AHS,AHR}. Just as for $\KO$, this topological pushforward comes from Thom classes~$\sigma_\nu$ associated with the normal bundle of a fiberwise embedding and the commuting triangle,
\beq
\begin{tikzpicture}[baseline=(basepoint)];
\node (A) at (0,0) {$\TMF^\bullet(M)$};
\node (B) at (6,0) {$\TMF_{\cs}^{\bullet+N-n}(B\times \R^N)$};
\node (C) at (2.5,-1.5) {$\TMF^{\bullet - n}(B).$};
\draw[->] (A) to node [above=1pt] {$\sigma_\nu\smallsmile$} (B);
\draw[->] (A) to node [left=.1in] {$\pi_!^{\rm top}$} (C);
\draw[->] (B) to node [right=.1in] {$\Sigma^{-N}$} (C);
\path (0,-.75) coordinate (basepoint);
\end{tikzpicture}
\label{TMFtoppush}
\eeq
The existence of the Thom class $\sigma_\nu$ requires the family $\pi \colon M\to B$ to have a string structure. 

To date, the only known construction of the Thom class in TMF uses purely homotopy theoretic techniques, and no tandem analytic pushforward has been constructed. The key open questions concerning an index theorem in TMF are as follows. 

\begin{quest} \label{q1}Is there a model for TMF wherein classes can be constructed using differential-geometric and analytic methods? 
\end{quest}
\begin{quest} \label{q2}
Is there an analytic pushforward $\pi_!^{\rm an}\colon \TMF^\bullet(M)\to \TMF^{\bullet-n}(B)$?
\end{quest}
\begin{quest} \label{q3}
Is there a differential-geometric construction of Thom classes in $\TMF$?
\end{quest}
\begin{quest} \label{q4} Do the analytic and topological pushforwards agree?
\end{quest}

Optimism that these questions can be answered in the affirmative comes from physics.

\subsection{Expectations from physics}\label{physmot}

A long-held belief is that some space of 2-dimensional supersymmetric quantum field theories represents~TMF as an infinite loop space or spectrum~\cite{Segal_Elliptic,ST04,HuKriz,ST11,VafaTMF}, thereby answering Question~\ref{q1}. One way of phrasing this expectation is 
\beq
&&\TMF(M)\cong \left\{\begin{array}{c}{\rm families\ of\ 2d\ field\ theories\ with\ } \mathcal{N}= (0,1)\\
{\rm \ supersymmetry\ parametrized\ by\ }M\end{array}\right\} /{\rm deformations},\label{eq:TMFspace}
\eeq
e.g., see~\cite[\S1]{ST11} or~\cite[\S3]{VafaTMF}.
Following Witten's original work on the subject, supporting evidence for~\eqref{eq:TMFspace} comes from several known relationships between $\K$-theory and \emph{1-dimensional} supersymmetric quantum field theory, which is more commonly known as \emph{supersymmetric quantum mechanics}. 

A bare-bones description of a quantum mechanical system consists of a Hilbert space (of \emph{states}) and an operator that generates time-evolution (the \emph{Hamiltonian}). In a \emph{supersymmetric} quantum mechanical system, the Hilbert space of states is $\Z/2$-graded and the Hamiltonian has a specified odd square root~\cite[\S1]{susymorse}; a prototypical example takes as its $\Z/2$-graded Hilbert space the $L^2$-sections of a spinor bundle of a spin manifold, and the Dirac operator as an odd operator on spinors (see the next paragraph). The map from supersymmetric quantum mechanical systems to K-theory views this odd square root of the Hamiltonian as an odd Fredholm operator, and then identifies such Fredholm operators with a representing space for~$\K$-theory~\cite{ST04,PokmanPhD,Kitaev,HST}. The difficulty in applying similar methods to~TMF is that 2-dimensional supersymmetric quantum field theory remains rather nebulous mathematically. In particular, there is no obvious analog to the space of Fredholm operators. As such, the viability of this proposed answer to Question~\ref{q1} is perhaps best tested in relation to the physics pertaining to Questions~\ref{q2} and~\ref{q3}. 

The most significant insight comes from studying the analog of Question~\ref{q2} in $\K$-theory: famously, the analytic pushforward can be interpreted as a quantization map for supersymmetric mechanics~\cite{Wittenconstraints,susymorse,Alvarez,FriedanWindey,WitteninStrings}. We review this briefly. Supersymmetric classical mechanics studies paths $x\colon \R\to M$ in a Riemannian manifold~$M$ together with an odd spinor, $\psi\in \Gamma(\R,\Pi \mathbb{S}\otimes x^*TM)$ equipped with the action functional\footnote{The classical action is best thought of as a Lagrangian density when integrating over non-compact manifolds, as these integrals need not converge. We commit this standard abuse throughout the introduction. We are also a bit cavalier about measures of integration; these are essentially determined by the chosen spin structure (on $\R$ in the above case), but we prefer not to belabor this point in the present discussion.}
\beq
&&\mathcal{S}_\sigma(x,\psi)=\frac{1}{2}\int (\langle \dot x,\dot x\rangle+\langle \psi,\nabla_t\psi\rangle ) \label{eq:11sigmaact}
\eeq
where $\nabla_t$ is the covariant derivative along the standard vector field $\partial_t$ on $\R$ using the pullback of the Levi-Civita connection on $M$.
As alluded to above, quantization of this supersymmetric classical theory constructs a $\Z/2$-graded Hilbert space of states together with an operator that generates time-evolution with a specified odd square root. In this case, the correct Hilbert space of states is sections of the spinor bundle on~$M$ with its $\Z/2$-grading, and so in particular $M$ needs to be a spin manifold for this quantization to make sense (e.g., see~\cite[pg.~43]{5lectures}). The Hamiltonian is the Dirac Laplacian $\slashed{D}^2$ that has the Dirac operator $\slashed{D}$ as an odd square root. The index of the Dirac operator also has an interpretation in this quantum mechanical theory. The time-evolution semigroup is $\exp(-t\slashed{D}^2)$ for $t\in \R_{>0}$ the time parameter, and the \emph{partition function} is the (super) trace of $\exp(-t\slashed{D}^2)$. The McKean--Singer formula shows this is independent of~$t$ and equal to the index of the Dirac operator~\cite{McKeanSinger}. 

Witten computed the partition function of an analogous supersymmetric quantum field theory one dimension higher~\cite{WittenDirac}, thereby generalizing the index of the Dirac operator. The starting classical theory has as fields maps~$x\colon \R^2/\Z\to M$ from cylinders (i.e., certain 2-manifolds rather than 1-manifolds) into a Riemannian manifold~$M$ together with an odd spinor $\psi\in \Gamma(\R^2/\Z,\Pi \overline{\mathbb{S}}\otimes x^*TM)$. Equip these fields with the action functional
\beq
&&\mathcal{S}_\sigma(x,\psi)=\frac{1}{2}\int (\langle \partial x,\bar\partial x\rangle+\langle \psi,\partial_\nabla\psi\rangle ),\label{eq:21sigmaact}
\eeq
where $\mathcal{S}_\sigma$ depends on an identification $\R^2/\Z\cong \C/\tau\Z$ for $\tau\in \mathbb{H}\subset\C$ in the upper-half plane to define the Cauchy--Riemann operators $\partial$ and $\bar\partial$.
Cylinders mapping to~$M$ should be thought of as propagating strings: using the diffeomorphism $\R^2/\Z\cong \R\times S^1$, one can view~$x\colon \R^2/\Z\to M$ as a map from~$\R$ into the free loop space of~$M$. This sets up a correspondence between the 2-dimensional theory associated to~\eqref{eq:21sigmaact} for fields valued in~$M$ and the 1-dimensional theory considered previously, but with fields valued in the free loop space of~$M$. For this reason, Witten suggestively calls the partition function of the associated 2-dimensional quantum field theory \emph{the index of the Dirac operator in loop space}. Using physical arguments (namely, supersymmetric localization) Witten was able to compute this partition function as a function of~$\tau$ in the upper-half plane. This is the original construction of the Witten genus. Witten predicted that this function would be a modular form when~$M$ has a string structure, which was verified computationally by Zagier~\cite{Zagiermodular}. One might hope that a deeper mathematical understanding of this quantum field theory would give a construction of the analytic pushforward in Question~\ref{q2}. 

A second key insight from physics is that Thom classes in $\KO$ also afford a description in terms of a quantum mechanical system; this is present (at least philosophically) in the work of Mathai and Quillen~\cite[3$^{\rm rd}$ paragraph of~\S0]{MathaiQuillen}. For a vector space~$V$ and element $v\in V$, we can consider the space $C^\infty(\R,\Pi \mathbb{S}\otimes V)$ of $V$-valued (odd) spinors equipped with the classical action functional
\beq
&&\mathcal{S}_\MQ(\psi)=\frac{1}{2}\int (\langle \psi,D \psi\rangle+\langle \psi,v\rangle)\qquad \psi\in C^\infty(\R,\Pi\mathbb{S}\otimes V),\label{eq:11MQ}
\eeq
where $D$ is the Dirac operator on~$\R$ acting on $V$-valued spinors. This is the classical action for the free fermion with an added source term, $\langle v,\psi\rangle$. The quantization of this classical theory constructs a class in~$\KO_\cs^{{\rm dim}(V)}(V)$ via a $V$-family of Clifford-linear operators. Indeed, the translational action of $ \Pi \mathbb{S}\otimes V$ on $C^\infty(\R,\Pi \mathbb{S}\otimes V)$ preserves the classical action, which implies that the quantum Hilbert space must carry an action by the Weyl algebra of~$\Pi\mathbb{S}\otimes V$. This Weyl algebra can be identified with~${\rm Cl}(V)$, the Clifford algebra of~$V$. The time-evolution operator commutes with this Clifford action, and so letting $v\in V$ vary in the action functional~\eqref{eq:11MQ} we obtain a $V$-family of Clifford-linear operators. A suitable families-version of this, starting with a vector \emph{bundle} rather than a vector \emph{space}, constructs the cocycle representative of the~$\KO$ Thom class from Remark~\ref{rmk:Thom}. 

One can consider the analog of the above theory one dimension higher. For a vector space~$V$, we have the space $C^\infty(\R^2/\Z,\Pi\mathbb{S}\otimes V)$ of $V$-valued odd spinors on a cylinder equipped with the classical action 
\beq
&&\mathcal{S}_\MQ(\psi)=\frac{1}{2}\int (\langle \psi,\bar\partial \psi\rangle+\langle \psi,v\rangle)\qquad \psi\in C^\infty(\R^2/\Z,\Pi\mathbb{S}\otimes V)\label{eq:21MQ}
\eeq
where now $\bar\partial$ is the chiral Dirac operator that again depends on an identification $\R^2/\Z\cong \C/\tau\Z$. 
When $v=0$, the action functional~\eqref{eq:21MQ} is the well-studied \emph{free fermion} in 2-dimensions, whose quantization is expected to construct the TMF Euler class~\cite{ST04,chris_andre_TMF}. Letting $v\in V$ vary, the expectation is that the quantization will construct the TMF Thom class. 

Unfortunately, promoting the above ideas from physics to rigorous constructions in algebraic topology presents significant technical challenges. At the outset, one must give a precise mathematical definition of 2-dimensional supersymmetric quantum field theory that is robust enough to encode topological modular forms. In light of the complicated torsion in~$\pi_*\TMF$, such a definition is likely to be subtle. To analogize, the appropriate topology on the space of Fredholm operators is delicate in the analytic model for~$\KO$. It seems that one would have to be either very clever or very lucky to happen upon the ``correct" definition for 2-dimensional supersymmetric quantum field theories from the point of view of~TMF. 

Rather than face this problem head-on, our approach below is to remain agnostic about the ultimate definition, and instead start small. We study basic aspects of the relevant physical theories above to see what geometry and analysis emerges naturally. Specifically we develop a mathematically rigorous framework for \emph{perturbative} quantization of the classical supersymmetric field theories, outlined in~\S\ref{sec:pertquant} below. The homotopy-theoretic counterpart also simplifies: we recover a model for \emph{complexified}~TMF together with analytic pushforwards and geometrically-constructed Thom classes, giving answers to Questions~\ref{q1}, \ref{q2}, \ref{q3} and~\ref{q4} in this simplified setting. 

\subsection{The complexified index theorem}\label{sec:cindex}

We start by describing the context for this simplified index theorem within algebraic topology. One can characterize complexifications of the analytic and topological pushforwards in Theorem~\ref{indexthm} using versions of the Riemann--Roch theorem.
The Chern character $\ch\colon \KO^\bullet(M)\to \KO^\bullet(M)\otimes \C\cong \HP^\bullet(M)$ in this case lands in 4-periodic ordinary cohomology with $\C$-coefficients. This affords descriptions of the complexified Thom class and the complexified analytic pushforward in terms of the commuting squares
\beq
\begin{tikzpicture}[baseline=(basepoint)];
\node (A) at (0,0) {$\KO^\bullet(M)$};
\node (B) at (5,0) {$\HP^\bullet(M)$};
\node (C) at (0,-1.5) {$\KO^{\bullet +m}_\cs(V)$};
\node (D) at (5,-1.5) {$\HP^{\bullet +m}_\cs(V),$}; 
\draw[->] (A) to node [above=1pt] {$\ch$} (B);
\draw[->] (A) to node [left=1pt] {$ \alpha_V\smallsmile$} (C);
\draw[->] (C) to node [above=1pt] {$\ch$} (D);
\draw[->] (B) to node [right=1pt] {$\smallsmile [u_V]\smallsmile [\hat{A}(V)]^{-1}$}(D);
\path (0,-.75) coordinate (basepoint);
\end{tikzpicture}\label{diag:RRThom}
\eeq
and 
\beq
&&
\begin{tikzpicture}[baseline=(basepoint)];
\node (A) at (0,0) {$\KO^\bullet(M)$};
\node (B) at (5,0) {$\HP^\bullet(M)$};
\node (C) at (0,-1.5) {$\KO^{\bullet - n}(B)$};
\node (D) at (5,-1.5) {$\HP^{\bullet -n}(B),$}; 
\draw[->] (A) to node [above=1pt] {$\ch$} (B);
\draw[->] (A) to node [left=1pt] {$\pi_!^{\rm an}$} (C);
\draw[->] (C) to node [above=1pt] {$\ch$} (D);
\draw[->] (B) to node [right=1pt] {$\int_{M/B} -\smallsmile [\hat{A}(T(M/B))]$}(D);
\path (0,-.75) coordinate (basepoint);
\end{tikzpicture}\label{diag:RR}
\eeq
respectively, where $V\to M$ is an $m$-dimensional vector bundle with spin structure,~$\HP_\cs$ is compactly supported cohomology, $[u_V]$ is the Thom class in ordinary cohomology, $[\hat{A}(V)]^{-1}$ is the inverse $\hat{A}$-class of~$V$, and is $[\hat{A}(T(M/B))]\in \HP^0(M)$ is the $\hat{A}$-class of the spin family. The crucial point is that the characteristic classes $[\hat{A}(T(M/B))]$ and $[u_V]\smallsmile[\hat{A}(V)]^{-1}$ control the complexified analytic and topological pushforwards, respectively. We observe that these characteristic classes can be defined for any oriented manifold or oriented vector bundle; the spin condition is unnecessary. 

\begin{rmk}\label{rmk:complexifiedindex11}
The equality of complexified pushforwards $\pi_!^{\rm top}\otimes \C=\pi_!^{\rm an}\otimes \C$ can easily be proven independent from the full index theorem. Indeed, it is a simple consequence of the Thom isomorphism in ordinary cohomology together with the equality of characteristic classes,~$[\hat{A}(\nu)]^{-1}=[\hat{A}(T(M/B))]$, which in turn follows from the fact that $\nu\oplus T(M/B)$ is a trivial vector bundle on~$M$. 
\end{rmk}

In the case of elliptic cohomology, the Chern--Dold character ${\rm ch}\colon \TMF(M)\to \TMF(M)\otimes \C\cong \H(M;\MF)$ lands in ordinary cohomology with values in weakly holomorphic modular forms over~$\C$ (see Remark~\ref{rmk:tens}). The Thom class in TMF and any would-be analytic pushforward sit in the diagrams
\beq
\begin{tikzpicture}[baseline=(basepoint)];
\node (A) at (0,0) {$\TMF^\bullet(M)$};
\node (B) at (5,0) {$\H^\bullet(M;\MF)$};
\node (C) at (0,-1.5) {$\TMF_\cs^{\bullet +m}(V)$};
\node (D) at (5,-1.5) {$\H^{\bullet +m}_\cs(V;\MF).$}; 
\draw[->] (A) to node [above=1pt] {$\ch$} (B);
\draw[->] (A) to node [left=1pt] {$\sigma_V\smallsmile $} (C);
\draw[->] (C) to node [above=1pt] {$\ch$} (D);
\draw[->] (B) to node [right=1pt] {$\smallsmile [u_V]\smallsmile [\Wit(V)]^{-1}$}(D);
\path (0,-.75) coordinate (basepoint);
\end{tikzpicture}\label{diag:WitRR1}
\eeq
and
\beq
\begin{tikzpicture}[baseline=(basepoint)];
\node (A) at (0,0) {$\TMF^\bullet(M)$};
\node (B) at (5,0) {$\H^\bullet(M;\MF)$};
\node (C) at (0,-1.5) {$\TMF^{\bullet - n}(B)$};
\node (D) at (5,-1.5) {$\H^{\bullet -n}(B;\MF).$}; 
\draw[->] (A) to node [above=1pt] {$\ch$} (B);
\draw[->,dashed] (A) to node [left=1pt] {$\pi_!^{\rm an}$} (C);
\draw[->] (C) to node [above=1pt] {$\ch$} (D);
\draw[->] (B) to node [right=1pt] {$\int_{M/B} -\smallsmile [\Wit(T(M/B))]$}(D);
\path (0,-.75) coordinate (basepoint);
\end{tikzpicture}\label{diag:WitRR2}
\eeq
respectively, where $V\to M$ is now a vector bundle with string structure, $\pi\colon M\to B$ is a family of string manifolds, $[\Wit(V)]^{-1}$ is the inverse Witten class and $[\Wit(T(M/B))]$ is the Witten class of the family. These characteristic classes now live in $\TMF^0(M)\otimes \C\cong \H^0(M;\MF)$. 
The key point is that $[\Wit(T(M/B))]$ and $[u_V]\smallsmile[\Wit(V)]^{-1}$ control the complexified analytic and topological pushforwards, respectively. We remark that these characteristic classes can be defined for any oriented manifold or oriented vector bundle with $[p_1(TM)]=0$ or $[p_1(V)]=0$ rationally; in this case we say that the manifold or vector bundle admits a \emph{rational string structure}. This requirement of a rational string structure is crucial for modularity. If the condition fails, $[\Wit(V)]^{-1}$ and $[\Wit(T(M/B))]$ only define classes with coefficients in quasimodular forms. In particular, the downward arrows on the right in~\eqref{diag:WitRR1} and~\eqref{diag:WitRR2} are not defined, at least not with the stated target.

\begin{rmk}\label{rmk:complexifiedindex21}
The equality of pushforwards in complexified $\TMF$ in this case is again a simple consequence of the Thom isomorphism in ordinary cohomology together with the equality of characteristic classes,~$[\Wit(\nu)]^{-1}=[\Wit(T(M/B))]$. 
\end{rmk}

The complexified index theorem leads us to the following goals, which are more modest versions of Questions~\ref{q1}-\ref{q4}.

\begin{goal} \label{g1}Is there a cocycle model for $\TMF\otimes \C$ that has an interpretation in terms of the supersymmetric $\sigma$-model?
\end{goal}
\begin{goal} \label{g2}
Can a cocycle representative of the Witten class in~$\TMF\otimes \C$ be constructed by applying quantization methods to the classical supersymmetric $\sigma$-model? 
\end{goal}
\begin{goal} \label{g3}
Can a cocycle representative of the Thom class in $\TMF\otimes \C$ be constructed by applying quantization methods to the Mathai--Quillen classical field theory?
\end{goal}
\begin{goal} \label{g4} Do the associated analytic and topological pushforwards for $\TMF\otimes \C$ agree at the level of cocycles?
\end{goal}

This paper provides affirmative (and mathematically precise) answers to the above questions, thereby giving a bridge between homotopy theory and the ideas from physics. The constructions of cocycles we have in mind involve path integral quantization of the supersymmetric theories from~\S\ref{physmot}. A formal version of the steepest descent approximation allows one to make the standard path integral constructions in physics rigorous. We begin by describing our mathematical approach to the relevant \emph{classical} supersymmetric field theories, inspired by the formalism developed in~\cite{Supersolutions,5lectures}.

\subsection{Classical (supersymmetric) field theory}\label{sec:susyft}
A vague---and perhaps overly simplistic---description of a classical field theory consists of a pair of data: (1) a \emph{space of fields} $\Phi$, and (2) a \emph{classical action functional} $\mathcal{S}\colon \Phi\to \C$. A \emph{symmetry} of the classical field theory is an automorphism of~$\Phi$ that preserves~$\mathcal{S}$. \emph{Classical observables} are functions on~$\Phi$, and these carry an action by the symmetries.

Typically, classical field theories have infinite-dimensional spaces of fields. Examples include spaces of maps between manifolds, or spaces of sections of vector bundles. A smooth structure on these types of spaces is often important in analyzing the field theory; for example, \emph{classical solutions} are the subspace of~$\Phi$ on which~$d\mathcal{S}=0$ and so, in particular, $\mathcal{S}$ needs to be a differentiable function on~$\Phi$ for this to make sense. One approach to infinite-dimensional manifolds is as sheaves on the site of manifolds. From this point of view, one can refine the data of a classical field theory to: (1) a sheaf $\Phi$ on the site of manifolds, and (2) a morphism of sheaves $\Phi\to C^\infty(-)$, where the target is the sheaf that assigns to a manifold its algebra of smooth functions.  Symmetries are then automorphisms in the category of sheaves preserving~$\mathcal{S}$, and observables are morphism of sheaves $\Phi\to C^\infty(-)$. Here is an example to keep in mind.

\begin{ex} \label{ex:easyfields}
Fix Riemannian manifolds~$M$ and~$\Sigma$ and let $\Phi=\Map(\Sigma,M)$ denote the sheaf on the site of smooth manifolds that assigns to~$S$ the set of maps $S\times \Sigma\to M$. Let $\mathcal{S}$ be the energy functional: 
$$
\mathcal{S}(\phi)=\frac{1}{2}\int_\Sigma \|d\phi\|^2,\qquad \phi\colon S\times \Sigma\to M,
$$
where $\|d\phi\|^2$ makes use of the Riemannian structures on $\Sigma$ and $M$. The integral is fiberwise over the projection $S\times \Sigma \to S$, so that the above assigns a function on $S$ to a map $\phi\colon S\times \Sigma\to M$, and therefore a morphism of sheaves $\mathcal{S}\colon \Map(\Sigma,M)\to C^\infty(-)$. 
Let $\Iso(\Sigma)$ be the isometry group of $\Sigma$, meaning the subgroup of the diffeomorphism group preserving the Riemannian structure on~$\Sigma$. Then $\Iso(\Sigma)$ acts by precomposition on $\Map(\Sigma,M)$ and preserves $\mathcal{S}$. The classical observables therefore furnish a representation of this isometry group. It can be useful to think of this representation as a sheaf on the stacky quotient~$[\Map(\Sigma,M)\sq \Iso(\Sigma)]$. Stacks like this one will feature prominently below. 
\end{ex}

One established approach to \emph{supersymmetric} field theory uses the formalism of \emph{supermanifolds} e.g., see~\cite{Supersolutions,5lectures}. Roughly speaking, supermanifolds are manifolds with even and odd coordinates; we review the basics in~\S\ref{appen:super}. The simple-minded definition of a classical supersymmetric field theory is just as above, replacing manifolds as test objects with supermanifolds, i.e., passing from the site of manifolds to the site of supermanifolds. Hence, $\Phi$ is a sheaf on the site of supermanifolds and~$\mathcal{S}$ is a morphism of sheaves~$\Phi\to C^\infty(-)$. Observables and symmetries are as before. The upshot is that when~$\Phi$ is a sheaf on supermanifolds, it can have ``odd" symmetries---a hallmark of supersymmetric field theory.

As outlined in~\S\ref{physmot} there are two types of classical supersymmetric field theories we wish to consider: (1) the supersymmetric $\sigma$-model that depends on a smooth manifold~$M$, and (2) the Mathai--Quillen field theory that depends on a real oriented vector bundle $V\to M$. It will be instructive to consider these field theories in dimensions $d=0,1,2$, with the cases $d=1,2$ related to K-theory and TMF, respectively, as outlined in~\S\ref{physmot}. The~$d=0$ case is the (somewhat degenerate, yet instructive) super-geometric interpretation of de~Rham cohomology. To treat these classical field theories rigorously, we work with the \emph{superspace} formulation of the theory. Essentially, this repackages the space of fields to be in the form of Example~\ref{ex:easyfields}, except that~$\Sigma$ is now a supermanifold; see~\cite{Supersolutions,5lectures} for introductory accounts to this formalism.

\begin{ex}[Supersymmetric $\sigma$-models]\label{ex:sigma}
Fix a Riemannian manifold $M$ and let $d\in \{0,1,2\}$. The \emph{space of fields} for the supersymmetric $\sigma$-model with target~$M$ is the sheaf whose value on a supermanifold~$S$ consists of pairs $(\ell,\phi)$ where $\ell\colon S\times \Z^d\hookrightarrow S\times \R^{d}$ is an $S$-family of based lattices in~$\R^d$ and $\phi\colon (S\times \R^{d|1})/\Z^d\to M$ is a map, where the $\Z^d$-quotient is determined by~$\ell$. Define an odd vector field~$D$ on $\R^{d|1}$ as
\beq
D:=\left\{\begin{array}{ll} \partial_\theta & d=0\\ \partial_\theta+i\theta\partial_t & d=1\\ \partial_\theta-\theta\partial_{\bar z} & d=2,\end{array}\right.\label{eq:Ddef}
\eeq
where $\theta$ is the odd coordinate on $\R^{d|1}$, $t$ is the standard coordinate on $\R$ when $d=1$, and $(z,\bar z)$ are complex coordinates on $\R^2\cong \C$ when $d=2$. These vector fields descend to~$(S\times \R^{d|1})/\Z^d$, and we use the same symbol $D$ to denote the vector field on this quotient. 
Equip the space of fields with the \emph{classical action} 
\beq
\mathcal{S}_\sigma(\ell,\phi)=\left\{\begin{array}{cl} 0 & d=0\\ \displaystyle\int\langle \partial_t\phi,D\phi\rangle & d=1\\ \displaystyle\int\langle \partial_z\phi,D\phi\rangle & d=2\end{array}\right.\label{eq:superLag}
\eeq
where the integrals are Berezinian integrals (see~\S\ref{appen:super}) over the (compact) fibers $(S\times \R^{d|1})/\Z^d\to S$, and $\langle-,-\rangle$ denotes the Riemannian pairing on~$TM$. We will consider the symmetries of the supersymmetric $\sigma$-model in detail in~\S\ref{sec:deRham}, \S\ref{sec:KC} and~\S\ref{sec:TMFC}. The most important symmetries are of the flavor of Example~\ref{ex:easyfields}, where the analog of~$\Iso(\Sigma)$ is a super Lie group of symmetries acting by precomposition. 
\end{ex}

\begin{rmk}\label{rmk:fields}
We indicate how the example above connects with the classical field theories from~\S\ref{physmot} with action functionals ~\eqref{eq:11sigmaact} and~\eqref{eq:21sigmaact}. The fiberwise inclusion $i_0$ of reduced manifolds
$$
x\colon (S\times \R^d)/\Z^d\stackrel{i_0}{\hookrightarrow} (S\times \R^{d|1})/\Z^d\stackrel{\phi}{\to} M
$$
determines a map $x\colon (S\times\R^d)/\Z^d\to M$ from an $S$-family of (ordinary) circles or tori to~$M$. Using the odd vector field $D$ on $(S\times \R^{d|1})/\Z^d$ we obtain a section
$$
\psi=i_0^*(D\phi)\in \Gamma((S\times \R^d)/\Z^d;x^*\Pi TM)
$$ 
that can be identified with a section of the twisted spinor bundle for the non-bounding spin structure (for which the spinor bundle on $\R^d/\Z^d$ is trivial). We often return to this description of fields---called \emph{component fields}---in computations. We refer to~\cite[pp.~649-656 and~663-665]{strings1} for derivations of the component versions of the classical actions~\eqref{eq:11sigmaact} and~\eqref{eq:21sigmaact} from the superspace versions~\eqref{eq:superLag}. 
\end{rmk}

\begin{ex}[The Mathai--Quillen theory]Let $d\in \{0,1,2\}$ and fix a real oriented vector bundle $p\colon V\to M$ with metric $\langle-,-\rangle$ and compatible connection~$\nabla$. We can pull back~$V$ over itself
\beq
\begin{tikzpicture}[baseline=(basepoint)];
\node (A) at (0,0) {$p^* V$};
\node (B) at (3,0) {$V$};
\node (C) at (0,-1.5) {$V$};
\node (D) at (3,-1.5) {$M$}; 
\draw[->] (A) to (B);
\draw[->, bend left] (C) to node [left=1pt] {$\x$} (A);
\draw[->] (A) to (C);
\draw[->] (C) to node [above=1pt] {$p$} (D);
\draw[->] (B) to (D);
\path (0,-.75) coordinate (basepoint);
\end{tikzpicture}\nonumber
\eeq
and the pullback has a tautological section~$\x$. The \emph{space of fields} for the Mathai--Quillen classical field theory consists of triples $(\ell,\phi,\psi)$ where $\ell\colon S\times \Z^d\hookrightarrow S\times \R^d$ is an $S$-family of based lattices, $\phi\colon (S\times \R^{d|1})/\Z^d\to V$ is a map, and $\psi\in \Gamma((S\times \R^{d|1})/\Z^d,\phi^*\Pi p^*V)$ is a section, where $\Pi$ denotes the parity reversal functor. The \emph{classical action} functional on sections is
\beq
&&\SMQ(\ell,\phi,\psi)=\int \left(\frac{1}{2}\langle \psi,\nabla_D\psi\rangle+\frac{i}{\vol^{1/2}}\left\langle \psi,\phi^*\x \right\rangle\right),\label{eq:functional}
\eeq
where $\vol$ is the fiberwise volume of the family~$(S\times \R^d)/\Z^d\to S$ (set to $1$ when $d=0$), $D$ is as in~\eqref{eq:Ddef}, $\langle-,-\rangle$ is the pullback along~$\phi$ of the metric pairing on sections of~$V$, and the integral is again a Berezinian integral over the fibers of $(S\times \R^{d|1})/\Z^d\to S$. 
\end{ex}

\begin{ex}[A combined theory] \label{ex:het} For~$d\in \{0,1,2\}$, $M$ a Riemannian manifold, and $V\to M$ a real oriented vector bundle with metric and connection, consider the Mathai--Quillen space of fields with the action functional
$$
\mathcal{S}(\ell,\phi,\psi)=\SMQ(\ell,\phi,\psi)+\mathcal{S}_\sigma(\ell,\phi). 
$$
For $V=\{0\}$, this classical field theory recovers the supersymmetric $\sigma$-model from Example~\ref{ex:sigma}. If instead we consider pairs $(\ell,\phi)$ with $\mathcal{S}_\sigma(\ell,\phi)=0$, we recover part of the Mathai--Quillen classical field theory. 
\end{ex}

\begin{rmk} When $d=2$, the classical field theory in Example~\ref{ex:het} is an important part of heterotic string theory, e.g., see~\cite[\S2.8]{WitteninStrings}. \end{rmk}

\subsection{Perturbative quantization via steepest descent}\label{sec:pertquant}

Given a classical field theory with fields $\Phi$ and classical action $\mathcal{S}$ one can try to construct various pieces of an associated \emph{quantum} field theory. A basic object is the assignment of a quantum expectation value to a classical observable. Naively, this is the evaluation of the integral
\beq
&&\langle f\rangle:=\int f e^{-\mathcal{S}}[d\phi] \qquad f\in C^\infty(\Phi) \qquad\qquad {\rm (naive)}\label{pathint1}
\eeq
for some measure $[d\phi]$. Of course, these sorts of integrals in quantum field theory are notoriously difficult to make rigorous. In the examples of interest below, we can use a formal steepest descent approximation to define quantum expectation values. In finite-dimensions this approximation takes the form
\beq
&&\int_\Phi f e^{-\mathcal{S}}[d\phi] \approx \int_{\Phi_0} \frac{f e^{-\mathcal{S}}}{ {\det}'(\Hess(\mathcal{S}))^{1/2}}[d\phi]_0\label{eq:steepest}
\eeq
where $\Phi_0={\rm min}(\mathcal{S})\subset \Phi$ denotes the minima of $\mathcal{S}$, $\det'(\Hess(\mathcal{S}))$ denotes the determinant of the Hessian of $\mathcal{S}$ on the orthogonal complement of its kernel, and $[d\phi]_0$ is a measure on~$\Phi_0$. We observe that using this as a formal approximation in infinite dimensions results in a well-defined integral provided that~$\Phi_0$ is finite-dimensional and smooth (so that a volume form $[d\phi]_0$ can be chosen), and that some appropriate (regularized) version of $\det'(\Hess(\mathcal{S}))^{-1/2}$ can be defined. 

\begin{rmk} Within the calculus of Feynman diagrams and perturbative quantum field theory, the approximation~\eqref{eq:steepest} picks out the 1-loop contribution to the path integral, e.g., see~\cite[Theorem~3.5]{EtingofQFT}. In this sense, formally applying the approximation~\eqref{eq:steepest} implements a version of perturbative quantization.
\end{rmk}

\begin{rmk} By physical reasoning---namely supersymmetric localization---the approximation~\eqref{eq:steepest} should be exact for the supersymmetric field theories considered in this paper. In effect, this is an infinite-dimensional incarnation of the Duistermaat--Heckman formula~\cite{DuistermaatHeckman}. When~$d=1$ this was explained by Atiyah~\cite{Atiyahast}, giving a link between supersymmetric quantum mechanics and the (usual) index theorem. The $d=2$ case is the subject of~\cite{BElocalization}.
\end{rmk} 

We formally apply the Fubini theorem and steepest descent approximation described above to the field theory in Example~\ref{ex:het}, 
\beq
\langle f\rangle&=&\int f\exp(-(\mathcal{S}_\sigma+\SMQ))[d\psi][ d\phi]\nonumber\\
&\stackrel{{\rm Fubini}}{=}&\left(\int f\exp(-\mathcal{S}_\sigma)\left(\int \exp(-\SMQ)[d\psi]\right)[d\phi]\right)\nonumber\\
&\stackrel{{\rm steepest}}{\approx}&\left(\int \frac{f\exp(-\mathcal{S}_\sigma)}{{\det}_{\rm reg}'(\HessS)^{1/2}} \left(\int \frac{\exp(-\SMQ)}{{\det}_{\rm reg}(\Hess(\SMQ))^{1/2}} [d\psi]_0\right)[d\phi]_0\right).\nonumber
\eeq
The final line above can be given rigorous meaning, which we will take as a definition for the quantum expectation value~$\langle f\rangle$. The constructions of analytic and topological pushforwards concern the case that~$f$ is independent of~$\psi$, which we have assumed in the formal approximation above; however, the methods below can be applied to general observables~$f$. We start by outlining the case of just the~$\sigma$-model, i.e., when~$V=\{0\}$. 


\begin{ex}[Steepest descent for the $\sigma$-model]\label{ex:sigmahess}
In terms of the action functionals~\eqref{eq:11sigmaact} and~\eqref{eq:21sigmaact}, the critical locus on which the classical action is minimized consists of those maps $x\colon (S\times \R^d)/\Z^d\to M$ that are a family of constant maps (i.e., factor through the projection to $S$) and those spinors $\psi\in \Gamma((S\times\R^d)/\Z^d;\Pi \mathbb{S}\otimes x^*TM)$ that are covariantly constant (when $d=2$ we use that an anti-holomorphic function on a compact complex manifold is constant). We observe that $\mathcal{S}_\sigma$ vanishes on these maps. In terms of the superspace formulation, these maps are characterized as the ones that factor as
$$
(S\times \R^{d|1})/\Z^d\to S\times \R^{0|1}\to M.
$$
Hence, we have
$$
\Phi_0=\{{\rm based\ lattices\ in\ }\R^d\}\times \Map(\R^{0|1},M)\subset \Phi
$$
which is a finite-dimensional and smooth subspace of $\Phi$. The Hessian of the classical action on $\Phi_0$ is the quadratic pairing
\beq
\HessS(X,Y)=\left\{\begin{array}{ll} 0 & d=0\\ \langle X,\nabla_{i\partial_t}\nabla_DY\rangle & d=1\\ \langle X,\nabla_{\partial_z}\nabla_DY\rangle & d=2\end{array}\right.\label{eq:sigmahess}
\eeq
with the vector field $D$ as in~\eqref{eq:Ddef}. For a physics derivation of these claims (stated in terms of component fields), see~\cite[\S4]{Alvarez} or~\cite[\S1.2.3]{WitteninStrings} for~$d=1$ and \cite[\S3.1]{Alvarezetal} for~$d=2$. We will construct a map from classical observables to quantum expectation values via the formula
\beq
\langle f\rangle=\int \frac{f}{{\sdet}_{\rm reg}'(\HessS)^{1/2}}[d\phi]_0\label{eq:sigmapush}
\eeq
where the integral is along the finite-dimensional fibers of the projection $\Phi_0\to \{{\rm based\ lattices}\}$, and the (super) determinant will be interpreted as a regularized product of finite-dimensional determinants. 
\end{ex}

\begin{rmk}
The above characterization of ${\rm min}(\mathcal{S}_\sigma)=\Phi_0$ is a finesse: algebras of functions on supermanifolds have nilpotent elements, so minima don't quite make sense in the superspace formulation of the $\sigma$-model from Example~\ref{ex:sigma}. The convention for $\Phi_0$ in the physics literature follows the above argument in component fields, which is equivalent to considering those maps on which $\mathcal{S}_\sigma$ vanishes~\cite{susymorse,Alvarez,Alvarezetal,WitteninStrings}. This agrees with our stated characterization above. 
\end{rmk}

\begin{rmk} Costello's construction of the Witten genus of a complex manifold~\cite{costello_WG1,costello_WG2} uses similar ideas from perturbative quantum field theory, and indeed, the geometry of our construction is very much inspired by his work. However our formalism is different, both in that it connects directly with the classical field theory employed by Witten's original construction, and it applies to all smooth manifolds without the requirement of a complex structure. The output is also different: we obtain pushforwards in complexified~TMF, whereas Costello builds a more algebraically sophisticated object, namely a \emph{factorization algebra}. 
\end{rmk}

\begin{ex}[Steepest descent for the Mathai--Quillen theory] \label{ex:MQ}The Mathai--Quillen classical action~\eqref{eq:functional} defines a family of field theories parameterized by maps $\phi\colon (S\times \R^{d|1})/\Z^d\to V$. Following the steepest descent prescription above, we first restrict attention to those maps $\phi$ that factor through $S\times \R^{0|1}$. On this restriction, we consider a finite-dimensional integral over the constant sections of~$\phi^*V$, which are called the \emph{zero modes}. The integrand we take is the Mathai--Quillen classical action modified by the determinant of the Hessian of the action, 
\beq
\int \frac{\exp(-\SMQ)}{\sdet_{\rm reg}'(\nabla_D)^{1/2}}[d\psi]_0={\det}_{\rm reg}'(\nabla_D)^{-1/2}\int \exp(-\SMQ)[d\psi]_0,\label{MQsplit}
\eeq
where the equality follows from the fact that $\sdet_{\rm reg}'(\nabla_D)^{1/2}$ turns out to be independent of the zero modes. We will find that the integral of~$\exp(-\SMQ)$ over the zero modes exactly contributes a cocycle representative $u_V$ of the Thom class in de~Rham cohomology (with $\C$ coefficients).
The regularized determinant of $\nabla_D$ supplies the correct Riemann--Roch factor to give a cocycle representative of the Thom class in complexified $\KO$ or complexified $\TMF$: $\hat{A}(V)^{-1}$ when $d=1$ and $\Wit(V)^{-1}$ when~$d=2$. 
\end{ex}

\section{Detailed statements of results}\label{sec:results}

We now elaborate on the brief summary of results from the beginning of~\S\ref{sec:intro}. As stated there, we divide the constructions into areas I-IV corresponding to \S\ref{sec:I}-\S\ref{sec:IV} below. When $d=2$, these results realize Goals~\ref{g1}-\ref{g4}, and so give affirmative answers to Questions~\ref{q1}-\ref{q4} for complexified TMF. 

\subsection{Results I: the cocycle models}\label{sec:I}

The main definition that leads to our cocycle model is an elaboration on the space of fields for the supersymmetric $\sigma$-model (Example~\ref{ex:sigma}) together with the action by certain symmetries. 
A toy version of this (without supermanifolds) is the stack $[\Map(\Sigma,M)\sq \Iso(\Sigma)]$ from Example~\ref{ex:easyfields}.
The symmetries relevant in the cases below come from \emph{rigid conformal maps} $(S\times \R^{d|1})/\Z^d\to (S\times \R^{d|1})/\Z^d$ between $S$-families, which are a certain class of fiberwise isomorphism. Following the usual definition of symmetry in classical field theory (see~\S\ref{sec:susyft}), these maps preserve the action functional~\eqref{eq:superLag} in the following sense. Given a map $\phi\colon (S\times \R^{d|1})/\Z^d\to M$, precomposition with a rigid conformal map leaves the value of the classical action~\eqref{eq:superLag} unchanged. We refer to~\S\ref{sec:deRham}, \S\ref{sec:KC} and~\S\ref{sec:TMFC} for the precise definition of fiberwise rigid conformal maps. 

\begin{defn} \label{defn:fields}
For $d=0,1,2$, let~$\L^{d|1}(M)$ be the stackifcation of the prestack whose objects over~$S$ are pairs $(\ell,\phi)$ where $\ell\colon S\times \Z^d\to S\times \R^d$ is an $S$-family of based lattices and $\phi\colon (S\times \R^{d|1})/\Z^d\to M$ is a map of supermanifolds for $\Z^d$-action determined by~$\ell$. Morphisms over a base change $S\to S'$ are commuting triangles
\beq
\begin{tikzpicture}[baseline=(basepoint)];
\node (A) at (-1,0) {$(S\times \R^{d|1})/\Z^d$};
\node (B) at (4,0) {$(S'\times \R^{d|1})/\Z^d$};
\node (C) at (1.5,-1.5) {$M$};
\draw[->] (A) to (B);
\draw[->] (A) to node [left=1pt]{$\phi$} (C);
\draw[->] (B) to node [right=1pt]{$\phi'$} (C);
\path (0,-.75) coordinate (basepoint);
\end{tikzpicture}\nonumber
\eeq
where the horizontal arrow is a fiberwise rigid conformal map between families of supermanifolds. 
\end{defn}

\begin{rmk}
The stacks $\L^{d|1}(M)$ are super geometric versions of the $d$-fold loop stack of~$M$, whence the notation. 
\end{rmk}

The steepest descent approximation~\eqref{eq:steepest} and its application in Example~\ref{ex:sigmahess} prompts the following definition. 

\begin{defn}\label{defn:factor}
Define $\L^{d|1}_0(M)$ as the full substack of $\L^{d|1}(M)$ with objects $(\ell,\phi)$ where $\phi\colon (S\times \R^{d|1})/\Z^d\to M$ factors through the projection to $S\times \R^{0|1}$, 
$$
(S\times \R^{d|1})/\Z^d\twoheadrightarrow S\times \R^{0|1}\to M,
$$
and the first arrow is induced by the standard projection $\R^{d|1}\to \R^{0|1}$.
\end{defn}

\begin{rmk} The assignment $M\mapsto \L^{d|1}_0(M)$ is natural in~$M$: a map $M\to M'$ determines a functor $\L^{d|1}_0(M)\to \L^{d|1}_0(M')$. \end{rmk}

\begin{rmk} When $d=0$ the factorization condition in Definition~\ref{defn:factor} is vacuous: $\L^{0|1}_0(M)=\L^{0|1}(M)$. This case and its relation to de~Rham cohomology is due to~\cite{HKST}. 
\end{rmk}

For each~$d$ there are line bundles $\omega^{1/2}$ over $\L^{d|1}(\pt)=\L^{d|1}_0(\pt)$. These are analogs of (an odd square root of) the Hodge bundle on the moduli stack of elliptic curves. These lines pull back to $\L^{d|1}_0(M)$ along the functor $\L^{d|1}_0(M)\to \L^{d|1}_0(\pt)$ induced by the canonical map $M\to \pt$, and we denote tensor powers of this pullback by $\omega^{\otimes \bullet/2}$. The assignment $M\mapsto \Gamma(\L_0^{d|1}(M);\omega^{\otimes\bullet/2})$ defines a sheaf on the site of smooth manifolds. We will compare these with sheaves of cocycles for various cohomology theories defined over~$\C$. 

\begin{rmk}\label{rmk:tens}
Define the cohomology theories $\KO\otimes \C$ and $\TMF\otimes \C$ by smashing the spectra representing $\KO$ and $\TMF$ with the spectrum representing $\H(-;\C)$, cohomology with $\C$ coefficients. The resulting cohomology theories (being defined over $\C$) are determined by their value on a point. We recall that $\KO(\pt)\otimes \C\cong \C[\beta,\beta^{-1}]$ where $|\beta|=-4$ and $\TMF(\pt)\otimes \C\cong \MF$, the graded ring of weakly holomorphic modular forms. Hence, $\KO\otimes \C$ is ordinary cohomology with coefficients~$\C[\beta,\beta^{-1}]$ and $\TMF\otimes\C$ is ordinary cohomology with coefficients~$\MF$. We use the notation $\C[\beta,\beta^{-1}]^j$ and $\MF^j$ for the $j$th graded piece of these coefficient rings (in the cohomological grading). When $M$ is compact, we have isomorphisms of rings $\KO(M)\otimes_\Z \C\cong \H(M;\C[\beta,\beta^{-1}])$ and $\TMF(M)\otimes_\Z \C\cong \H(M;\MF)$. When $M$ is not compact, in a mild abuse of notation we still write $\KO(M)\otimes \C$ and $\TMF(M)\otimes \C$ to denote the cohomology theories $\KO\otimes \C$ and $\TMF\otimes \C$ applied to~$M$, respectively. 
\end{rmk}

\begin{thm}\label{thm1}
There are isomorphisms of sheaves of graded algebras 
\beq
\Gamma(\L^{0|1}(-),\omega^{\otimes\bullet/2})&\cong & \Omega^\bullet_{\cl}(-)\label{deRhamcocycle}\\
\Gamma(\L^{1|1}_0(-),\omega^{\otimes\bullet/2})&\cong&\bigoplus_{i+j=\bullet} \Omega^i_{\cl}(-;\C[\beta,\beta^{-1}]^j) \label{K-theorycocycle} \\
\mathcal{O}(\L^{2|1}_0(-),\omega^{\otimes\bullet/2})&\cong&\bigoplus_{i+j=\bullet} \Omega^i_{\rm cl}(-; \MF^j).\label{TMFcocycle}
\eeq
This gives differential cocycle models in the sense of Hopkins--Singer~\cite{HopSing} for ordinary cohomology with $\C$-coefficients, $\KO\otimes \C$ and $\TMF\otimes \C$, respectively. 
\end{thm}

\begin{rmk}
The part of the statement regarding Hopkins--Singer differential cocycles follows easily from the isomorphism of sheaves. Indeed, by Remark~\ref{rmk:tens}, the relevant cohomology theories are determined by their coefficient rings. An application of the de~Rham theorem along with the isomorphisms~\eqref{deRhamcocycle}-\eqref{TMFcocycle} presents functions on $\L^{d|1}_0(M)$ as a differential cocycle model for the given cohomology theory applied to~$M$. 
\end{rmk}

\begin{rmk} The definition of holomorphic sections $\mathcal{O}(\L^{2|1}_0(M),\omega^{\otimes\bullet/2})\subset \Gamma((\L^{2|1}_0(M),\omega^{\otimes\bullet/2}))$ in~\eqref{TMFcocycle} records holomorphic dependence on the moduli of (super) tori. Such holomorphic sections are the ones relevant in (chiral) supersymmetric field theory; see Remark~\ref{rmk:holo}. 
\end{rmk}

To prove the isomorphisms in Theorem~\ref{thm1}, we show that the stacks have covers by ordinary supermanifolds 
\beq
u\colon \Map(\R^{0|1},M)&\twoheadrightarrow&\L^{0|1}(M)\label{eq:01atlasin}\\
u\colon \R^\times\times \Map(\R^{0|1},M)&\twoheadrightarrow &\L^{1|1}_0(M)\label{eq:11atlasin}\\
u\colon \Lat\times \Map(\R^{0|1},M)&\twoheadrightarrow & \L^{2|1}_0(M)\label{eq:21atlasin}
\eeq
where $\Lat$ is the space of based, oriented lattices in $\C\cong \R^2$ (see~\S\ref{sec:TMFC} and~\S\ref{appen:modform}), and $\R^\times=\R\setminus \{0\}$ is the space of based lattices in~$\R$. This allows us to construct objects over the corresponding stacks by first constructing them on these covers and then studying descent along~\eqref{eq:01atlasin}-\eqref{eq:21atlasin}. For now we observe that the isomorphism
\beq
C^\infty(\Map(\R^{0|1},M))\cong \Omega^\bullet(M),\label{eq:piTM}
\eeq
yields inclusions of algebras
\beq
&&C^\infty(\L^{0|1}_0(M))\hookrightarrow C^\infty(\Map(\R^{0|1},M))\cong \Omega^\bullet(M),\label{incl1}\\ 
&&C^\infty(\L^{1|1}_0(M))\hookrightarrow C^\infty(\R^\times\times \Map(\R^{0|1},M))\cong \Omega^\bullet(M;C^\infty(\R^\times))\\
&&C^\infty(\L^{2|1}_0(M))\hookrightarrow C^\infty(\Lat\times \Map(\R^{0|1},M))\cong \Omega^\bullet(M;C^\infty(\Lat))\label{incl3}
\eeq
where the isomorphisms are as Fr\'echet vector spaces. There are similar statements for sections of $\omega^{\otimes\bullet/2}$. These inclusions are what connect the stacks $\L^{d|1}_0(M)$ with cocycle models for various cohomology theories over~$\C$, and indeed, the proof of Theorem~\ref{thm1} boils down to computing the images of~\eqref{incl1}-\eqref{incl3} and similar statements for sections of $\omega^{\otimes \bullet /2}$.

\subsection{Results II: analytic pushforwards}\label{sec:II}

The analytic and topological pushforwards come from variations on a general construction of families of operators over $\L^{d|1}_0(M)$. The input to this construction is a real vector bundle~$V\to M$. From this we can construct two types of vector bundles over $\L^{d|1}_0(M)$. The first is a bundle of $V$-valued bosons, $\Bos(V)\to \L^{d|1}_0(M)$ that has $\Gamma((S\times \R^{d|1})/\Z^d,\phi^*V)$ as its sections at~$(\ell,\phi)$. The second is a bundle of $V$-valued fermions, $\Fer(V)\to \L^{d|1}_0(M)$ that has $\Gamma((S\times \R^{d|1})/\Z^d,\phi^*\Pi V)$ as its sections at~$(\ell,\phi)$, where $\Pi$ denotes the parity reversal functor. 
If we endow~$V$ with a connection $\nabla$ we get families of operators acting on the sections of $\Bos(V)$ and $\Fer(V)$. The family of operators acting on sections of $\Bos(V)$ is 
\beq
\Delta_V=\left\{\begin{array}{ll} 0 & d=0\\ \nabla_{\partial_t}\nabla_D & d=1\\ \nabla_{\partial_z}\nabla_D & d=2\end{array}\right.\label{eq:sigmaops}
\eeq
for the vector field $D$ in~\eqref{eq:Ddef}. The family of operators acting on sections of $\Fer(V)$ is 
\beq
D_V=\nabla_D.
\eeq
Regularized determinants of these families of operators for particular choices of vector bundles will construct analytic and topological pushforwards. 

The input for the analytic pushforward is $\pi\colon M\to B$, a family of oriented Riemannian manifolds with fiber dimension~$n$. The vector bundle~$\BosT\to \L^{d|1}_0(M)$ has $\Gamma((S\times \R^{d|1})/\Z^d,\phi^*T(M/B))$ as its sections at~$(\ell,\phi)$ where $T(M/B)$ is the vertical tangent bundle of $\pi\colon M\to B$. More geometrically, the bundle~$\BosT$ is the vertical tangent bundle of the projection~$\L^{d|1}(M)\to \L^{d|1}(B)$ restricted to $\L^{d|1}_0(M)\subset \L^{d|1}(M)$. Motivated by~\eqref{eq:sigmahess}, define the operators from~\eqref{eq:sigmaops}
$$
\HessS=\Delta_{T(M/B)}
$$
acting on sections of this bundle. We construct a pushforward relative to the covers~\eqref{eq:01atlasin}-\eqref{eq:21atlasin} using the formula~\eqref{eq:sigmapush} and then study descent to the stack. So for $d=0,1,2$, let $U_M\to \L^{d|1}_0(M)$ denote the relevant cover from~\eqref{eq:01atlasin}-\eqref{eq:21atlasin}; this cover turns out to be natural in~$M$, giving a 2-commutative square
\begin{equation}
\begin{array}{c}
\begin{tikzpicture}[node distance=3.5cm,auto]
  \node (A) {$U_M$};
  \node (B) [node distance= 4.5cm, right of=A] {$\L^{d|1}_0(M)$};
  \node (C) [node distance = 1.5cm, below of=A] {$U_B$};
  \node (D) [node distance = 1.5cm, below of=B] {$\L^{d|1}_0(B).$};
  \draw[->>] (A) to node [above] {$u$} (B);
  \draw[->] (A) to node [swap] {$\widetilde{\pi}$} (C);
  \draw[->>] (C) to (D);
  \draw[->] (B) to node {$\pi$} (D);
\end{tikzpicture}\end{array}
\label{eq:atlas}
\end{equation}
Using~\eqref{eq:piTM}, integration of differential forms gives a pushforward along~$\widetilde{\pi}$. This descends to the stacks as a linear map
$$
\int\colon \Gamma(\L^{d|1}_0(M);\omega^{\otimes\bullet/2})\to \Gamma(\L^{d|1}_0(B);\omega^{\otimes(\bullet-n)/2}),
$$
where the shift in line bundles is a reflection of the fact that integration changes the degree of differential forms. Sections of the pullback $u^*\BosT$ to the atlas have a $d$-dimensional torus action, which decomposes sections into finite-dimensional weight spaces. We will define $\sdet_{\rm reg}'(\HessS)^{-1/2}$ as the (infinite) product of these (finite-dimensional) super determinants restricted to each weight space, suitably normalized. Super determinants are also sometimes called \emph{Berezinian determinants.} 

To define the normalization on each weight space, consider the operator $D_n:=D_{\underline{\R}^n}$ acting on sections of $\Fer(n):=\Fer(\underline{\R}^n)$ for $\underline{\R}^n\to M$ the trivial bundle with the standard metric and trivial connection. These sections again decompose into an infinite sum of finite-rank vector bundles according to the weight spaces for the torus action. When $d=1$ these weight spaces are indexed by $k\in \Z$, and when $d=2$ they are indexed by $(k,l)\in \Z\times \Z$. To fix the sign of the square root of the determinant, we follow a standard procedure wherein the regularized determinant is the product over ``half" of the weight spaces. To this end, define
$$
\Z^2_+=\{(k,l)\in \Z^2\mid k>0 \ {\rm or} \ k=0, l>0\},\quad \Z_+=\{k\in \Z \mid k>0\}. 
$$

\begin{defn}\label{defn:21sigmareg}
Define the \emph{regularized super determinant} as
\beq
(d=1)\qquad \sdet_{\rm reg}'(\HessS)^{1/2}&:=&\prod_{k\in \Z_+} \sdet\left(\HessS|_{\begin{smallmatrix} k{\rm th \ weight} \\ {\rm space}\end{smallmatrix}}\right)\cdot \sdet\left((D_n)|_{\begin{smallmatrix}k{\rm th\ weight}\\ {\rm space}\end{smallmatrix}}\right)\nonumber \\
(d=2)\qquad \sdet_{\rm reg}'(\HessS)^{1/2}&:=&\prod_{(k,l)\in \Z^2_+} \sdet\left(\HessS|_{\begin{smallmatrix} (k,l){\rm th\ weight}\\ {\rm space}\end{smallmatrix}}\right)\cdot \sdet\left((D_n)|_{\begin{smallmatrix}(k,l){\rm th\ weight}\\ {\rm space}\end{smallmatrix}}\right)\nonumber 
\eeq
where 
$\sdet$ is the usual finite-dimensional super determinant of the restrictions of $\HessS$ and $D_n$ to each weight space. 
\end{defn}

\begin{rmk}
The sum of the weight spaces indexed by $k\in \Z_+$ and $(k,l)\in \Z^2_+$ can be viewed as defining a Lagrangian polarization of fields, which is a necessary choice to be made when quantizing. The choice of~$\Z_+$ when $d=1$ is standard (and nearly canonical). The choice of $\Z^2_+$ when~$d=2$ comes from applying the $d=1$ polarization for each fixed~$l$ (with a special case when $k=0$). This guarantees a compatibility between the regularized determinants, e.g., via dimensional reduction. Roughly stated, the constant term in the $q$-expansion of $\sdet_{\rm reg}'(\HessS)^{1/2}$ when $d=2$ coincides with the regularized determinant when $d=1$. This agrees with the constant term in the $q$-expansion of the Witten class being the $\hat{A}$-class. In part because of the compatibility with standard index theory, this choice of polarization is also standard in the physics literature, e.g., see~\cite[\S2.3]{WitteninStrings}. 
\end{rmk}

\begin{rmk}\label{rmk:ferrmk} The regularized determinant above is the partition function for a theory of bosons in $T(M/B)$ and free fermions in $\R^n$, omitting the zero modes. Hence, the regularization procedure can be understood physically as adding $n$ free fermions (valued in $\R^n$) to the original supersymmetric $\sigma$-model. This is closely related to Stolz and Teichner's definition of a spin structure on a vector bundle when $d=1$~\cite[Definition~2.3.1]{ST04}.\end{rmk}

\begin{rmk} One can also define the regularized determinants above using $\zeta$-regularization techniques. However, when~$d=2$, one must confront regularized products of \emph{complex} numbers, rather than the more standard theory that applies to positive real numbers. Such a $\zeta$-regularization can still be carried out in the complex case (e.g., using results from~\cite{zetaprod}), but the techniques are not widely known. For the sake of clarity and accessibility, we have adopted the above more concrete approach to regularization. 
\end{rmk}

We now apply Definition~\ref{defn:21sigmareg} and formula~\eqref{eq:sigmapush} to produce (analytic) pushforwards 
\beq
&&\widehat{\pi}_!^{\rm an}(f):=\int f\cdot \sdet_{\rm reg}'(\HessS)^{-1/2}\quad {\rm for} \ d=0,1, \ {\rm modified\ when} \ d=2\label{defn:analyticpushgen}
\eeq
on atlases along $\widetilde{\pi}$ in~\eqref{eq:atlas}, where the definition of $\widehat{\pi}_!^{\rm an}$ when $d=2$ requires a modification of the integrand depending on a rational string structure; see Theorem~\ref{thm:1}. We use the notation $\widehat{\pi}_!^{\rm an}$ to denote a cocycle (or differential) refinement of the analytic pushforward~$\pi_!^{\rm an}$. We recall that $n={\rm dim}(M)-{\rm dim}(B)$. 

\begin{prop}\label{01analytic}
When~$d=0$, $\sdet_{\rm reg}'(\HessS)^{-1/2}=1$ and $\widehat{\pi}_!^{\rm an}$ descends to the stack as integration of differential forms under the isomorphism~\eqref{deRhamcocycle},
\beq
&&\widehat{\pi}_!^{\rm an}(f):=\int f,\qquad\qquad \widehat{\pi}_!^{\rm an}\colon \Gamma(\L^{0|1}(M);\omega^{\otimes k/2})\to \Gamma(\L^{0|1}(B);\omega^{\otimes (k-n)/2}).\label{eq:defn:analyticpush}
\eeq
\end{prop}

\begin{thm} \label{prop:Aclass} When $d=1$, $\sdet_{\rm reg}'(\HessS)^{-1/2}$ converges absolutely and determines  a cocycle representative of the $\hat{A}$-class for the family $\pi \colon M\to B$. Consequently, $\widehat{\pi}_!^{\rm an}$ descends to the stack as a cocycle refinement of the complexification of the analytic pushforward in $\KO\otimes \C$
under the isomorphism~\eqref{K-theorycocycle},
$$
\widehat{\pi}_!^{\rm an}(f):=\int f\cdot \sdet_{\rm reg}'(\HessS)^{-1/2}=\int f\cdot \hat{A}(T(M/B)),
$$$$
\qquad \widehat{\pi}_!^{\rm an}\colon \Gamma(\L^{1|1}_0(M);\omega^{\otimes k/2})\to \Gamma(\L^{1|1}_0(B);\omega^{\otimes (k-n)/2}).
$$
\end{thm}

\begin{thm}\label{thm:1} When $d=2$, $\sdet_{\rm reg}'(\HessS)^{-1/2}$ converges conditionally. 
\begin{enumerate}
\item[(i)]
If the first Pontryagin form $p_1(T(M/B))$ is nonzero, then there is no choice of ordering of the conditionally convergent product for which $\sdet_{\rm reg}'(\HessS)^{-1/2}$ descends to a holomorphic function on the stack. 

\item[(ii)] For any particular choice of ordering, $\sdet_{\rm reg}'(\HessS)^{-1/2}$ descends to a section 
$$
\sdet_{\rm reg}'(\HessS)^{-1/2} \in \Gamma(\L^{2|1}_0(M);\mathcal{A}(p_1)),
$$ 
where the line bundle $\mathcal{A}(p_1)$ depends on the differential form~$p_1=p_1(T(M/B))$ and the choice of ordering for the conditionally convergent product. 

\item[(iii)] A choice of rational string structure on $\pi\colon M\to B$ determines a trivializing concordance of the line bundle $\mathcal{A}(p_1)$, meaning a line bundle with section over $\L^{2|1}_0(M\times \R)$ whose restriction to $\L^{2|1}_0(M\times \{1\})$ is~$\mathcal{A}(p_1)$ with section $\sdet_{\rm reg}'(\HessS)^{-1/2}$, and whose restriction to $\L^{2|1}_0(M\times \{0\})$ is the trivial line with section $\Wit(M/B)\in \mathcal{O}(\L^{2|1}_0(M))$, a cocycle representative of the Witten class of the family $\pi\colon M\to B$ in $\TMF^0(M)\otimes \C$. The differential analytic pushforward is then defined by
$$
\widehat\pi_!^{\rm an}(f) :=\int f\cdot \Wit(M/B),\qquad\qquad \widehat\pi_!^{\rm an}\colon \mathcal{O}(\L^{2|1}_0(M);\omega^{\otimes k/2})\to \mathcal{O}(\L^{2|1}_0(B);\omega^{\otimes(k-n)/2}),
$$
which is~\eqref{defn:analyticpushgen} modified by the concordance. This pushforward does not depend on the choice of ordering of the conditionally convergent regularized product nor does it depend on the choice of rational string structure.
\end{enumerate}
\end{thm}

\begin{rmk} A choice of \emph{rational string structure} on a family of oriented Riemannian manifolds $\pi \colon M\to B$ is a 3-form $H\in \Omega^3(M)$ such that $dH=p_1(T(M/B))$. \end{rmk}

\begin{rmk}
In physical language, the $p_1(T(M/B))$ obstruction to constructing the analytic pushfoward is an example of an \emph{anomaly}. The literature on anomalies is vast. A good mathematical introduction is~\cite{Freed_Det}. One reference relevant to path integral computations in the supersymmetric $\sigma$-model considered below is~\cite{Pilch}. 
\end{rmk}

\subsection{Results III: Mathai--Quillen Thom forms}\label{sec:III}

As before, let $d\in \{0,1,2\}$ and $p\colon V\to M$ be a real $m$-dimensional oriented vector bundle with metric and compatible connection, $\Fer(V)\to \L^{d|1}_0(M)$ the vector bundle of $V$-valued fermions described at the beginning of~\S\ref{sec:II} and $u\colon U_M\to \L^{d|1}_0(M)$ the previously chosen atlas~\eqref{eq:01atlasin}-\eqref{eq:21atlasin}. 
For~$p^*V\to V$ the pullback of~$V$ along itself, define the \emph{Mathai--Quillen classical action} functional on sections of~$u^*\Fer(p^*V)$ by the formula~\eqref{eq:functional}. Applying the steepest descent approximation in the form of~\eqref{MQsplit} constructs a function on the atlas, and we can study descent to the stack. 

Since Thom classes are compactly supported, we make the following definition.

\begin{defn}
Let $M$ be a smooth manifold, not necessarily compact. A section of a vector bundle over $\L^{d|1}_0(M)$ has \emph{compact support} if when pulled back to the cover~\eqref{eq:01atlasin}-\eqref{eq:21atlasin} it has compact support along the fibers of the projection to the cover with $M=\pt$, i.e., along the fibers of maps
$$
\Map(\R^{0|1},M)\to \pt,\ \ \ \R^\times\times \Map(\R^{0|1},M)\to \R^\times, \ \ \ \Lat\times \Map(\R^{0|1},M)\to \Lat
$$
in the respective cases, where we used $\Map(\R^{0|1},\pt)=\pt$. 
\end{defn}

\begin{rmk}\label{rmk:compactsupp}
The cocycle representatives of the Thom classes that we construct are only rapidly decreasing along the fibers of a vector bundle rather than compactly supported. We use the standard trick (e.g., see~\cite[pg.~56]{BGV}) of choosing a diffeomorphism from a vector bundle to its unit disk bundle to convert rapidly decaying sections into compactly supported ones. 
We commit the usual abuse of notation, letting $\Gamma_{\rm c}$ or $\mathcal{O}_c$ refer to either compactly supported or rapidly decreasing sections.
\end{rmk}

 As in the case of the $\sigma$-model, sections of the pullback of $\Fer(p^*V)$ to the atlases~\eqref{eq:01atlasin}-\eqref{eq:21atlasin} have an action by a $d$-dimensional torus where each weight space is finite-dimensional. The \emph{zero modes} are the sections fixed by this torus action. 

 \begin{prop}\label{prop:01MQ}\label{lem:d1zeromodes}
Let $\epsilon({\rm dim}(V))=1$ for ${\rm dim}(V)$ even and $i=\sqrt{-1}$ for ${\rm dim}(V)$ odd. For $d\in \{0,1,2\}$, the Berezinian integral $\int\exp(-\SMQ)[d\psi]_0$ of the restriction of $\exp(-\SMQ)$ to the zero modes of~$u^*\Fer(p^*V)$ equals the image of the Mathai--Quillen Thom form in ordinary de~Rham cohomology,
$$
\frac{\epsilon({\rm dim}(V))}{(2\pi)^{{\rm dim}(V)}}\int\exp(-\SMQ)[d\psi]_0=u_V\in \Gamma_c(\L^{d|1}_0(V);\omega^{{\rm dim}(V)/2})
$$ 
using the isomorphisms ~\eqref{deRhamcocycle},~\eqref{K-theorycocycle} and~\eqref{TMFcocycle}.
\end{prop}
\begin{rmk}
When $d=0$, the calculation behind Proposition~\ref{prop:01MQ} above is well-trodden terrain with the main ideas essentially due to Mathai and Quillen~\cite{MathaiQuillen}. These ideas are developed further in~\cite[\S1.6]{BGV}. The super-geometric formalism used below is in the same spirit as~\cite{Wu} and references therein.
\end{rmk}

The nonzero weight spaces define a family of operators over the stack whose regularized super determinant contributes a Riemann--Roch factor when $d=1,2$. We define the regularization using operators $\Delta_n=\Delta_{\underline{\R}^n}$ acting on sections of $\Bos(n)=\Bos(\underline{\R}^n)$ for $\underline{\R}^n\to V$ the trivial rank $n$ vector bundle with the standard metric and trivial connection. As in the discussion preceding Definition~\ref{defn:21sigmareg}, the weight spaces are indexed by $\Z$ when $d=1$ and~$\Z^2$ when $d=2$.

 \begin{defn}\label{defn:21MQreg}
Define the \emph{regularized super determinant} as
\beq
(d=1)\qquad \sdet_{\rm reg}'(\nabla_D)^{1/2}&:=&\prod_{k>0} \sdet\left(\nabla_D|_{\stackrel{k{\rm th \ weight}}{{\rm space}}}\right)\cdot \sdet\left((\Delta_n)|_{\stackrel{k{\rm th\ weight}}{{\rm space}}}\right)\nonumber \\
(d=2)\qquad \sdet_{\rm reg}'(\nabla_D)^{1/2}&:=&\prod_{(k,l)\in \Z^2_+} \sdet\left(\nabla_D|_{\stackrel{(k,l){\rm th\ weight}}{{\rm space}}}\right)\cdot \sdet\left((\Delta_n)|_{\stackrel{(k,l){\rm th\ weight}}{{\rm space}}}\right)\nonumber 
\eeq
where $\sdet$ is the usual finite-dimensional super determinant of the restrictions of $\nabla_D$ and $\Delta_n$ to each weight space. 
\end{defn}

\begin{rmk} This parallels Remark~\ref{rmk:ferrmk}. The regularized determinant above is the partition function for a theory of fermions in $p^*V$ and (a supersymmetric version of) free bosons in~$\R^n$. \end{rmk}

\begin{prop}\label{thm:thomkthy}
When $d=1$, $\sdet_{\rm reg}'(\nabla_D)^{-1/2}$ converges absolutely and the section 

$$
\sdet_{\rm reg}'(\nabla_D)^{-1/2} \int\exp(-\SMQ)[d\psi]_0=\hat{A}(V)^{-1}\cdot u_V \in \Gamma_{\rm c} (\L_0^{1|1}(V);\omega^{\otimes{\rm dim}(V)/2})\label{eq:11Thom}
$$
coincides with the Mathai--Quillen cocycle representative of the Thom form in complexified $\KO$-theory.
\end{prop}

\begin{thm}\label{thm:2}
When $d=2$, $\sdet_{\rm reg}'(\nabla_D)^{-1/2}$ converges conditionally. 
\begin{enumerate}
\item[(i)] If the first Pontryagin form $p_1(V)$ is nonzero, there is no choice of ordering the product in the regularized super determinant so that 
the function $\sdet_{\rm reg}'(\nabla_D)^{-1/2}$ descends to a holomorphic function on the stack. 
\item[(ii)] Any particular choice of ordering descends to a section of a line bundle
\beq
&&\sdet_{\rm reg}'(\nabla_D)^{-1/2}\int\exp(-\SMQ)[d\psi]_0\in \Gamma_{\rm c} (\L^{2|1}_0(V);\omega^{\otimes{\rm dim}(V)/2}\otimes\mathcal{A}(p_1)^\vee) \label{eq:21Thom}
\eeq
where $\mathcal{A}(p_1)$ depends on the differential form~$p_1=p_1(V)$ and the choice of ordering, and $\mathcal{A}(p_1)^\vee$ is the dual of $\mathcal{A}(p_1)$ from Theorem~\ref{thm:1}. 

\item[(iii)] A choice of rational string structure on~$V$ determines a trivializing concordance, meaning a line bundle with section over $\L^{2|1}_0(V\times \R)$ whose restriction to $\L^{2|1}_0(V\times \{1\})$ is~\eqref{eq:21Thom}, and whose restriction to $\L^{2|1}_0(V\times \{0\})$ is 
$$\Wit(V)^{-1}\cdot u_V\in \mathcal{O}_c(\L^{2|1}_0(V);\omega^{\otimes{\rm dim}(V)/2}),
$$ 
a cocycle representative of the Thom class in~$\TMF_{\rm c}^{{\rm dim}(V)} (V)\otimes \C$. The cocycle $\Wit(V)^{-1}\cdot u_V$ is independent of both the choice of rational string structure and the choice of ordering the conditionally convergent product. 
\end{enumerate} 
\end{thm}

\begin{defn}
Call the cocycle $\Wit(V)^{-1}\cdot u_V\in \mathcal{O}_c(\L^{2|1}_0(V);\omega^{\otimes {\rm dim}(V)/2})$ the \emph{elliptic Mathai--Quillen form}. 
\end{defn}

\begin{rmk}
The above construction also defines a cocycle for each real inner product space~$V$, viewed as a vector bundle over $M=\pt$. The result is a cocycle representative of the suspension class for the sphere obtained from the one point compactification of~$V$. This gives a physical interpretation for the suspension isomorphism in ordinary cohomology when $d=0$, complexified $\KO$-theory when $d=1$, and $\TMF\otimes \C$ when $d=2$. One might hope to refine this suspension isomorphism by considering more sophisticated quantizations of the Mathai--Quillen field theory, e.g., by extending down as discussed in~\S\ref{sec:intro}. Ultimately this should explain how the space of 2-dimensional field theories with $\mathcal{N}=(0,1)$ supersymmetry has the structure of an orthogonal ring spectrum. This gives a candidate answer to the question of where such a structure could originate, e.g., as asked in~\cite[\S3.1]{VafaTMF}.
\end{rmk} 
\subsection{Results IV: Index theorems}\label{sec:IV}

 Let $\pi \colon M\to B$ be a family of oriented Riemannian manifolds and $i\colon M\hookrightarrow B\times \R^N$ a fiberwise isometric embedding with normal bundle~$\nu$. Define the \emph{topological pushforward}, denoted $\widehat{\pi}_!^{\rm top}$ as the composition
$$
\Gamma(\L^{d|1}_0(M);\omega^{\otimes k/2}){\longrightarrow} \Gamma_c(\L^{d|1}_0(B\times \R^N);\omega^{\otimes (k+N-n)/2})\stackrel{\int}{\to} \Gamma(\L^{d|1}_0(B);\omega^{\otimes (k-n)/2}),\quad d=0,1
$$
$$
\mathcal{O}(\L^{2|1}_0(M);\omega^{\otimes k/2}){\longrightarrow} \mathcal{O}_c(\L^{2|1}_0(B\times \R^N);\omega^{\otimes (k+N-n)/2})\stackrel{\int}{\to} \mathcal{O}(\L^{2|1}_0(B);\omega^{\otimes (k-n)/2}),\quad d=2
$$
where the first arrow multiplies by the Mathai--Quillen representative of the Thom class constructed in the previous subsection (which requires a rational string structure when~$d=2$), and extends by zero to define a function on $\L^{d|1}_0(B\times \R^N)$ using a tubular neighborhood. The second arrow is the integration along $\L^{d|1}_0(B\times \R^N)\to \L^{d|1}_0(B)$ associated with integration of differential forms. As in the case of the analytic pushforward, we decorate $\widehat{\pi}_!^{\rm top}$ to indicate that it is a map on cocycles, i.e., a differential pushforward.

\begin{thm} \label{thm:3}
Let $\pi \colon M\to B$ be a family of oriented Riemannian manifolds and $i\colon M\hookrightarrow B\times \R^N$ a fiberwise isometric embedding. The analytic and topological pushforwards are equal 
$$
\widehat\pi_!^{\rm an}=\widehat\pi_!^{\rm top}\colon \Gamma(\L^{0|1}(M);\omega^{\otimes k/2})\to \Gamma(\L^{0|1}(B);\omega^{\otimes(k-n)/2})\qquad d=0,
$$
$$
\widehat\pi_!^{\rm an}=\widehat\pi_!^{\rm top}\colon \Gamma(\L^{1|1}_0(M);\omega^{\otimes k/2})\to \Gamma(\L^{1|1}_0(B);\omega^{\otimes(k-n)/2})\qquad d=1,
$$
$$
\widehat\pi_!^{\rm an}=\widehat\pi_!^{\rm top}\colon \mathcal{O}(\L^{2|1}_0(M);\omega^{\otimes k/2})\to \mathcal{O}(\L^{2|1}_0(B);\omega^{\otimes(k-n)/2})\qquad d=2,
$$
where $n={\rm dim}(M)-{\rm dim}(B),$ and we require a rational string structure on the family $\pi\colon M\to B$ when $d=2$. 
\end{thm}

\begin{rmk} A physical interpretation of the index theorem above is an identification between quantum expectation values for the linear supersymmetric $\sigma$-model (where the target is $\R^N$) and the nonlinear supersymmetric $\sigma$-model (where the target is $M$). 
\end{rmk}

\section{Pushforwards in de~Rham cohomology from $0|1$-dimensional field theories}\label{sec:01}

In this section we review well-known connections between de~Rham cohomology and $0|1$-dimensional field theories. At the heart of this story is the isomorphism~\eqref{eq:piTM} between functions on the superspace $\Map(\R^{0|1},M)$ and differential forms on $M$. The supergeometry of $\Map(\R^{0|1},M)$ then repackages the de~Rham complex of~$M$. This idea seems to have first appeared in Kontsevich's deformation quantization~\cite[\S7.2]{kont}. We also drew inspiration from how these ideas appear within the Stolz--Teichner program~\cite{HKST}. 

The development of the Mathai--Quillen formalism from this supergeometric point of view is a bit more diffuse. Mathai and Quillen themselves highlighted the ``super" aspect of their construction~\cite{MathaiQuillen} (see also~\cite[\S1.6]{BGV}), but at the time it seems that supergeometry was largely confined to the physics literature. Atiyah and Jeffrey's approach to topological quantum field theory~\cite{AtiyahJeffrey} sparked further developments that made this more explicit, both in the physics and mathematics literature; for example~\cite{BlauThompson2,BlauThompson,ZhengDong}, and see also~\cite{Wu} for a concise overview. 

\subsection{The superspace $\Map(\R^{0|1},M)$ and de~Rham cohomology}\label{sec:deRham}

Define the super Lie group~$\E^{0|1}\rtimes \R^\times$ as follows. Let $\E^{0|1}$ be $\R^{0|1}$ as a supermanifold with group structure gotten from the usual addition (i.e., vector space structure) on $\R^{0|1}$. Define the semidirect product $\E^{0|1}\rtimes \R^\times$ via the action $\theta\mapsto \mu \cdot \theta$, for $\theta\in\E^{0|1}(S)$ and $\mu \in\R^\times(S)$. We note that $\E^{0|1}\rtimes \R^\times$ has a Lie algebra generated by an odd vector field $D\in {\rm Lie}(\E^{0|1})$ and an even vector field $N \in {\rm Lie}(\R^\times)$, satisfying the relations
$$
[D,D]=2D^2=0, \quad [N,D]=D.
$$
We emphasize that since $D$ is odd, the first relation is nontrivial (whereas $[N,N]=0$ is automatic). There is an evident left action of~$\E^{0|1}\rtimes \R^\times$ on~$\R^{0|1}$. 
We will need a suitable families-version of this. 

\begin{defn}\label{defn:conf01} A \emph{fiberwise conformal map} $S\times \R^{0|1}\to S'\times \R^{0|1}$ over a base change $S\to S'$ is determined by a map
$$
S\times \R^{0|1}\to \E^{0|1}\rtimes \R^\times \times \R^{0|1}\to \R^{0|1}
$$
where the first arrow comes from an $S$-point of $\E^{0|1}\rtimes \R^\times$ and the second arrow is the action of this super Lie group on $\R^{0|1}$.
\end{defn}

\begin{rmk} The above is a simplified version of fiberwise isometries in Stolz and Teichner's rigid geometries~\cite[\S4.2]{ST11}; see also~\cite[\S6.3]{HST}.\end{rmk}

We introduce a little more terminology that will be useful below. An $S$-point of $\E^{0|1}\rtimes \R^\times$ is determined by an $S$-point of~$\E^{0|1}$ and of~$\R^\times$. We refer to the $S$-point of $\E^{0|1}$ as the \emph{super translation} part of the fiberwise conformal map and the $S$-point of $\R^\times$ as the \emph{dilation} part. 

The following is the specialization of Definition~\ref{defn:fields} to the $d=0$ case.

\begin{defn} \label{defn:01fields}The \emph{fields for the $0|1$-dimensional $\sigma$-model with target $M$}, denoted $\L^{0|1}(M)$, is the stack associated to the prestack whose objects over $S$ are maps $\phi\colon S\times \R^{0|1}\to M$ and whose morphisms over a base change $S\to S'$ are commuting triangles, 
\beq
\begin{tikzpicture}[baseline=(basepoint)];
\node (A) at (0,0) {$S\times \R^{0|1}$};
\node (B) at (3,0) {$S'\times \R^{0|1}$};
\node (C) at (1.5,-1.5) {$M$};
\draw[->] (A) to node [above=1pt] {$\cong$} (B);
\draw[->] (A) to node [left=1pt]{$\phi$} (C);
\draw[->] (B) to node [right=1pt]{$\phi'$} (C);
\path (0,-.75) coordinate (basepoint);
\end{tikzpicture}\label{01triangle}
\eeq
where the horizontal arrow is a fiberwise conformal map. For a map $M\to M'$, postcomposition $S\times \R^{0|1}\to M\to M'$ defines a morphism of stacks $\L^{0|1}(M)\to \L^{0|1}(M')$. 
\end{defn}


Define a morphism of stacks $\L^{0|1}(\pt)\to [\pt\sq \R^\times]$ whose value on objects is constant to $\pt$ and to a fiberwise conformal map~$S\times \R^{0|1}\to S'\times \R^{0|1}$ associates the dilation part,~$\mu\colon S\to \R^\times$ in the terminology explained after Definition~\ref{defn:conf01}. 

\begin{rmk} \label{rmk:stackad} Above we implicitly used that stackification is the left adjoint to the forgetful functor from stacks to prestacks~\eqref{eq:leftad}. Hence it suffices to specify the value of a morphism $\L^{0|1}(\pt)\to [\pt\sq \R^\times]$ on the prestack in Definition~\ref{defn:01fields}. \end{rmk} 

As reviewed in~\S\ref{appen:super}, there is a canonical odd complex line bundle over $[\pt\sq \R^\times]$. 

\begin{defn}\label{defn:01Hodge} Let $\omega^{\otimes k/2}$ denote the line bundle over $\L^{0|1}(M)$ that is the pullback of the $k$th tensor power of the dual of the canonical odd line bundle over $[\pt\sq \R^\times]$ along the composition,
$$
\L^{0|1}(M)\to \L^{0|1}(\pt)\to [\pt\sq\R^\times]
$$
where the first arrow is the functor induced by $M\to \pt$. 
\end{defn}

Equation~\eqref{deRhamcocycle} in Theorem~\ref{thm1} characterizes the sections of $\omega^{\otimes k/2}$. To compute these sections, we find a groupoid presentation of $\L^{0|1}(M)$. To this end, let $\Map(\R^{0|1},M)$ denote the (representable) presheaf on supermanifolds whose $S$-points are the set of maps~$S\times \R^{0|1}\to M$. There is an evident functor~$u\colon \Map(\R^{0|1},M)\to \L^{0|1}(M)$ associated with the evaluation map 
$$
\phi=\ev\colon \Map(\R^{0|1},M)\times\R^{0|1}\to M. 
$$
Directly from the definitions, we see that~$u$ is essentially surjective on $S$-points. Together with the fact that an $S$-family of fiberwise conformal maps is an $S$-point of $\E^{0|1}\rtimes \R^\times$, we get the following. 

\begin{lem}\label{lem:01present} There is an action groupoid presentation, $\L^{0|1}(M)\simeq [\Map(\R^{0|1},M)\sq \E^{0|1}\rtimes \R^\times]$ for the left action of $\E^{0|1}\rtimes \R^\times$ on $\Map(\R^{0|1},M)$ by the precomposition action on~$\R^{0|1}$. 
\end{lem}

\begin{rmk}\label{rmk:leftright} Usually precomposition actions are from the right. Turning the precomposition action of $\E^{0|1}\rtimes \R^\times$ on $\Map(\R^{0|1},M)$ into a left action uses inversion on the group. One can verify that this is the correct action by inspecting the source and target data in the commuting triangle~\eqref{01triangle}: if an isomorphism in $\L^{0|1}(M)$ over $S$ associated to a fiberwise conformal map $f\colon S\times \R^{0|1}\to S\times \R^{0|1}$ has source $\phi\colon S\times \R^{0|1}\to M$ and target $\phi'\colon S\times \R^{0|1}\to M$, then $\phi'=\phi\circ f^{-1}$, so the target does indeed come from precomposition with the inverse fiberwise conformal map.
\end{rmk}

\begin{proof}[Proof of Theorem~\ref{thm1}, $d=0$ case] 
The projection homomorphism $\E^{0|1}\rtimes \R^\times \to \R^\times$ defines a functor $\Map(\R^{0|1},M)\sq \E^{0|1}\rtimes \R^\times\to \pt\sq \R^\times$ between quotient Lie groupoids whose image in stacks is the morphism~$\L^{0|1}(M)\to [\pt\sq\R^\times]$ defining~$\omega^{1/2}$. Global sections of $\omega^{\otimes k/2}$ can therefore be computed as functions on $\Map(\R^{0|1},M)$ invariant under the $\E^{0|1}$-action and equivariant for the $\R^\times$-action (e.g., see~\cite[Corollary~7.19]{HKST}). Under the identification~\eqref{eq:piTM} between functions on $\Map(\R^{0|1},M)$ and differential forms on $M$, the $\E^{0|1}$-invariant functions are closed differential forms, and $\mu\in \R^\times$ acts on a $k$-form~$\alpha$ by $\alpha\mapsto \mu^{-k}\alpha$ (the sign comes from Remark~\ref{rmk:leftright}). Hence sections of~$\omega^{\otimes k/2}$ are closed differential forms of degree~$k$, 
\beq
\Gamma(\L^{0|1}(M);\omega^{\otimes k/2})\cong \Omega^k_{\rm cl}(M),\label{eq:01deRham}
\eeq
where to get forms in positive degree we emphasize that we use the \emph{dual} of the canonical bundle in Definition~\ref{defn:01Hodge}. Under this isomorphism the tensor product of sections coincides with multiplication of forms, proving the proposition. 
\ep

\subsection{The analytic pushforward}\label{sec:01sigmamodel}

Let $\pi\colon M\to B$ be a family of oriented manifolds, and let $\pi\colon \L^{0|1}(M)\to \L^{0|1}(B)$ denote the induced functor between stacks. The analytic pushforward is defined by the formula~\eqref{defn:analyticpushgen}, which a priori constructs a pushforward along~$\widetilde{\pi}$ in the diagram
\begin{equation}
\begin{array}{c}
\begin{tikzpicture}[node distance=3.5cm,auto]
  \node (A) {$\Map(\R^{0|1},M)$};
  \node (B) [node distance= 4.5cm, right of=A] {$\L^{0|1}(M)$};
  \node (C) [node distance = 1.5cm, below of=A] {$\Map(\R^{0|1},B)$};
  \node (D) [node distance = 1.5cm, below of=B] {$\L^{0|1}(B).$};
  \draw[->>] (A) to (B);
  \draw[->] (A) to node [swap] {$\widetilde{\pi}$} (C);
  \draw[->>] (C) to (D);
  \draw[->] (B) to node {$\pi$} (D);
\end{tikzpicture}\end{array}
\label{eq:01int}
\end{equation}
To define integration along $\widetilde{\pi}$, we recall that integration on supermanifolds depends on a choice of \emph{Berezinian measure}; see~\S\ref{appen:super} for a quick review. An orientation of~$M$ (in the usual sense) determines a Berezinian measure on~$\Map(\R^{0|1},M)$: integration of functions with respect to this Berezinian measure is simply integration of differential forms on~$M$. More generally, an orientation on the family $M\to B$ defines a Berezinian measure on the fibers~$\Map(\R^{0|1},M)\to \Map(\R^{0|1},B)$ for which integration is the usual fiberwise integration of differential forms.

\begin{proof}[Proof of Proposition~\ref{01analytic}]
With $d=0$, the analytic pushforward associated with the $\sigma$-model is particularly degenerate: the operator $\HessS$ from~\eqref{eq:sigmahess} is the zero operator. So the orthogonal complement of the zero eigenspace is the zero vector space, and hence $\sdet_{\rm reg}'(\HessS)^{-1/2}=1$. So all together, the analytic pushforward~\eqref{defn:analyticpushgen} at the atlas level is simply fiberwise integration of differential forms
\beq
\widehat{\pi}_!^{\rm an}(f)=\int_{M/B}f,\qquad \widehat{\pi}_!^{\rm an}\colon \Omega^\bullet(M) \to \Omega^\bullet(B).\label{eq:01analyticdef}
\eeq
using the isomorphisms $C^\infty(\Map(\R^{0|1},M))\cong \Omega^\bullet(M)$ and $C^\infty(\Map(\R^{0|1},B))\cong\Omega^\bullet(B)$. 
It remains to study descent of this analytic pushforward to $\pi \colon \L^{0|1}(M)\to \L^{0|1}(B)$ in~\eqref{eq:01int}.

Descent of the map~\eqref{eq:01analyticdef} amounts to unpacking standard facts about integration of differential forms. First we observe that $\E^{0|1}$-invariant functions on $\Map(\R^{0|1},M)$ are sent to $\E^{0|1}$-invariant functions on $\Map(\R^{0|1},B)$: the fiberwise integral of a closed form is closed. However, the transformation properties of a function under the $\R^\times$-action typically change under this map, corresponding to the shift in degree of a differential form by the dimension of the fiber of $\pi\colon M\to B$. 
Instead, we obtain a map between sections of line bundles
\beq
\begin{tikzpicture}[baseline=(basepoint)];
\node (A) at (0,0) {$\Gamma(\L^{0|1}(M),\omega^{\otimes \bullet/2})$};
\node (B) at (3,0) {$\Omega_{\cl}^\bullet(M)$};
\node (C) at (6,0) {$\Omega_{\cl}^{\bullet-n}(B)$};
\node (D) at (10,0) {$\Gamma(\L^{0|1}(B),\omega^{\otimes (\bullet-n)/2})$};
\draw[->] (A) to node [below=1pt] {$\cong$} (B);
\draw[->] (B) to node [below=1pt]{$\int_{M/B}$} (C);
\draw[->] (C) to node [below=1pt]{$\cong$} (D);
\draw[->, bend left=10] (A) to node [above=1pt] {$\widehat{\pi}_!^{\rm an}$} (D);
\path (0,-.75) coordinate (basepoint);
\end{tikzpicture}\nonumber
\eeq
as claimed.
\ep

\subsection{Mathai--Quillen Thom forms} \label{sec:01MQ}\label{sec:01MQmeasure} Let $V\to M$ be a real oriented vector bundle over a smooth manifold $M$ with metric and compatible connection. This determines a vector bundle over the supermanifold associated to~$M$, which we also denote by~$V\to M$.

\begin{rmk} Although it is standard to use the same notation for the vector bundle $V\to M$ in either manifolds or supermanifolds, it hides some subtleties: the structure sheaf of~$M$ as a supermanifold consists of \emph{complex-valued} smooth functions, and the sections of~$V$ (being defined as a locally free module over these functions) are the \emph{complexification} of the usual sheaf of sections of~$V$. This complexification will be implicit in the notation below. \end{rmk}

Partitions of unity imply that the sheaf of sections of~$V$ is determined by the $C^\infty(M)$-module of its global sections, which we denote by~$\mathcal{V}$. The metric on~$V$ defines a $C^\infty(M)$-bilinear pairing on~$\mathcal{V}$ with values in $C^\infty(M)$, denoted $\langle-,-\rangle$. Similarly, the input connection defines a covariant derivative operator $\nabla\colon \mathcal{V}\to \Omega^1(M)\otimes \mathcal{V}$. 
Let $\Pi \mathcal{V}$ denote the parity reversal of~$\mathcal{V}$. 
The pairing $\langle-,-\rangle$ also defines a pairing~$\langle-,-\rangle\colon \Pi \mathcal{V}\otimes \Pi \mathcal{V}\cong \mathcal{V}\otimes \mathcal{V}  \to C^\infty(M)$ and the covariant derivative determines an operator~$\nabla\colon \Pi\mathcal{V}\to \Omega^1_M\otimes \Pi\mathcal{V}$.
%

A map $\phi\colon S\times \R^{0|1}\to M$, determines a vector bundle~$\phi^*\Pi V$ over~$S\times \R^{0|1}$ by pullback. Global sections of this pullback are the $C^\infty(S\times \R^{0|1})$-module 
\beq
\Gamma(S\times \R^{0|1};\phi^*\Pi V)=\phi^{-1}\Pi \mathcal{V}=C^\infty(S\times \R^{0|1})\otimes_{C^\infty(M)} \Pi \mathcal{V},\label{eq:vb01}
\eeq
where $\phi^{-1}\Pi \mathcal{V}$ denotes the inverse image sheaf, and the $C^\infty(M)$-action on $C^\infty(S\times \R^{0|1})$ is through the algebra homomorphism $\phi^*\colon C^\infty(M)\to C^\infty(S\times \R^{0|1})$. 
 We observe that the commutative triangle~\eqref{01triangle} gives a morphism of vector bundles $\phi^*\Pi V\to \phi'^*\Pi V$ over $S\times \R^{0|1}\to S'\times \R^{0|1}$, by naturality of pullbacks. We use this naturality to define a vector bundle over the stack~$\L^{0|1}(M)$.

\begin{defn} For $V\to M$ a real oriented vector bundle, define a vector bundle $\Fer(V)$ over $\L^{0|1}(M)$ whose space of sections at an $S$-point $\phi\colon S\times \R^{0|1}\to M$ is the $C^\infty(S)$-module
$$
\Gamma(S\times \R^{0|1},\phi^*\Pi V)=C^\infty(S\times \R^{0|1})\otimes_{C^\infty(M)}\Pi\mathcal{V}.
$$
In terms of a sheaf of sections, this is the direct image sheaf along the projection $S\times \R^{0|1}\to S$ of the sheaf of sections of the pullback~$\phi^*\Pi V$. For a fiberwise conformal map~\eqref{01triangle} over a base change $S\to S'$, define a morphism of vector bundles as the one induced on the direct image sheaf by the composition of the pullback by the fiberwise conformal map $f\colon S\times \R^{0|1}\to S'\times \R^{0|1}$ followed by the action of scalar multiplication by the dilation factor~$\mu\in \R^\times(S)\subset C^\infty(S)^\ev$,
$$
\Gamma(S'\times \R^{0|1},\phi'^*\Pi V)\stackrel{f^*}{\to} \Gamma(S\times \R^{0|1},\phi^*\Pi V)\stackrel{\mu\cdot}{\to} \Gamma(S\times \R^{0|1},\phi^*\Pi V).
$$
\end{defn}

In light of Lemma~\ref{lem:01present}, the vector bundle $\Fer(V)$ can also be characterized as a vector bundle over the quotient Lie groupoid $\Map(\R^{0|1},M)\sq \E^{0|1}\rtimes \R^\times$. This has as data a vector bundle $u^*\Fer(V)$ over the objects $\Map(\R^{0|1},M)$ and an isomorphism ${\sf s}^*u^*\Fer(V)\cong {\sf t}^*u^*\Fer(V)$ over the morphisms~$\E^{0|1}\rtimes \R^\times\times \Map(\R^{0|1},M)$ where ${\sf s}$ and ${\sf t}$ denote the source and target maps in the groupoid. The explicit $C^\infty(\Map(\R^{0|1},M))$-module defined by~$u^*\Fer(V)$ is 
$$
C^\infty(\Map(\R^{0|1},M)\times\R^{0|1},\ev^*\Pi V)\cong C^\infty(\Map(\R^{0|1},M)\times \R^{0|1})\otimes_{C^\infty(M)}\Pi\mathcal{V},
$$ 
where the $C^\infty(M)$-module structure on $C^\infty(\Map(\R^{0|1},M)\times \R^{0|1})$ is inherited from pullback along the evaluation map $\ev\colon \Map(\R^{0|1},M)\times \R^{0|1}\to M$. 

The Mathai--Quillen classical action is a functional on sections of $u^*\Fer(p^*V)$ for $p^* V\to V$ the pullback of $V$ along itself and $u\colon \Map(\R^{0|1},V)\to \L^{0|1}(V)$ the atlas. To emphasize, here we are viewing the base~$V$ as a supermanifold (a locally ringed space) and $p^* V$ as a vector bundle over~$V$ (a locally free module over the structure sheaf) when forming $\L^{0|1}(V)$ and $\Fer(p^*V)$. We will define the classical action on sections of $u^*\Fer(p^*V)$ over objects of the groupoid presentation of~$\L^{0|1}(V)$, and then show that this functional descends to the stack. 
In the following, let $\x\in \Gamma(V,p^*\Pi V)$ denote the tautological odd section, i.e., the tautological section viewed as an element of the parity reversal of the module~$\Gamma(V,p^*V)$ in which the standard tautological section lives.  

\begin{defn} For $\phi=\ev\colon \Map(\R^{0|1},V)\times \R^{0|1}\to V$, and $\psi$ a section of $u^*\Fer(p^*V)$, define the value of the \emph{Mathai--Quillen classical action} as the function on $\Map(\R^{0|1},V)$
\beq
&&\SMQ(\phi,\psi):=\int \left(\frac{1}{2}\langle \psi,\nabla_D\psi\rangle+i\langle \psi,\phi^*\x\rangle\right)[d\theta]\label{eq:01MQaction}
\eeq
where $D=\partial_\theta$ and the integral is over the fibers of the projection 
\beq
\Map(\R^{0|1},V)\times \R^{0|1}\to \Map(\R^{0|1},V)\label{eq:intfibers01}
\eeq
with respect to the fiberwise Berezinian measure~$[d\theta]$ (characterized by $\int \theta [d\theta]=1$). In a slight abuse of notation, in the above $\nabla_D=\iota_D\phi^*\nabla$ denotes the pullback of the covariant derivative operator evaluated along the vector field~$D$.
\end{defn}

To connect the Mathai--Quillen classical action with the Mathai--Quillen Thom form, we make use of \emph{component fields} for sections of~$u^*\Fer(V)$ over $\Map(\R^{0|1},M)$, defined as 
\beq
&&\begin{array}{lll} \psi_1&:=&i_0^*\psi\in \Gamma(\Map(\R^{0|1},M),i_0^*\ev^*\Pi V), \\ 
\psi_0&:=&i_0^* (\nabla_D\psi)\in \Gamma(\Map(\R^{0|1},M),i_0^*\ev^*V).\end{array} \label{eq:component1}
\eeq
where $i_0\colon \Map(\R^{0|1},M)\hookrightarrow \Map(\R^{0|1},M)\times \R^{0|1}$ is the inclusion at~$0\in \R^{0|1}$, or equivalently the map on functions sending the coordinate $\theta\in C^\infty(\R^{0|1})$ to zero. This gives an isomorphism of $C^\infty(\Map(\R^{0|1},M))$-modules 
$$
\Gamma(\Map(\R^{0|1},M)\times \R^{0|1},\ev^*\Pi V)\stackrel{\sim}{\to} \Gamma(\Map(\R^{0|1},M),i_0^*\ev^*\Pi V)\oplus \Gamma(\Map(\R^{0|1},M),i_0^*\ev^*V)
$$
sending $\psi\mapsto (\psi_1,\psi_0)$. The inverse isomorphism is 
$$
(\psi_1,\psi_0)\mapsto {\rm pr}^*\psi_1+\theta {\rm pr}^*\psi_0\qquad {\rm pr}\colon \Map(\R^{0|1},M)\times \R^{0|1}\to \Map(\R^{0|1},M)$$ 
where ${\rm pr}$ is the projection. We caution that this projection is not invariant under the action of super translations (as these act nontrivially on $\Map(\R^{0|1},M)$), and these Taylor components are similarly not invariant. 
To simplify the notation below, we will omit the pullback along the projection, writing $\psi=\psi_1+\theta\psi_0$. We view this as a Taylor expansion in the odd variable~$\theta$ of the section~$\psi$, generalizing the standard formulas in supergeometry for Taylor expansions of functions in odd variables. 

Applying this to $p^*V\to V$, we similarly obtain components for $\ev^*\x$ as
\beq
\begin{array}{lll}
\x_0&:=&i_0^*\ev^*\x \in \Gamma(\Map(\R^{0|1},V),i_0^*\ev^*p^* V),\\
 \x_1&:=&i_0^*(\nabla_D \ev^* \x)\in \Gamma(\Map(\R^{0|1},V),i_0^*\ev^*p^* \Pi V).
 \end{array}\label{eq:component2}
\eeq

To connect component fields with vector-valued differential forms, we observe that the composition
$$
\Map(\R^{0|1},M)\stackrel{i_0}{\hookrightarrow} \Map(\R^{0|1},M)\times \R^{0|1}\stackrel{\ev}{\to}M
$$
is the projection $\Map(\R^{0|1},M)\to M$ induced by the inclusion $\pt\stackrel{0}{\hookrightarrow} \R^{0|1}$. This gives an isomorphism of $C^\infty(\Map(\R^{0|1},M))$-modules
$$
\Gamma(\Map(\R^{0|1},M),i_0^*\ev^*V)= C^\infty(\Map(\R^{0|1},M))\otimes_{C^\infty(M)} \mathcal{V}\cong \Omega^\bullet(M)\otimes_{\Omega^0(M)}\mathcal{V}\cong \Omega^\bullet(M;V)
$$
with $V$-valued differential forms, using the isomorphism $C^\infty(\Map(\R^{0|1},M))\cong \Omega^\bullet(M)$. Applying this to $p^*V\to V$, we identify component fields as elements
$$
\psi_0\in \Omega^\bullet(V,p^*V),\qquad \psi_1\in \Omega^\bullet(V,p^*\Pi V). 
$$
The value of the Mathai--Quillen classical action on these component fields is an element of~$\Omega^\bullet(V)$, which we can now calculate in terms of more familiar quantities.

\begin{lem}\label{lem:calccomp01}
The value of~$\SMQ$ on component fields $(\psi_1,\psi_0)$ is
$$
\SMQ(\psi_1,\psi_0)=\left(\langle \psi_0,\frac{1}{2}\psi_0+i\x\rangle+\langle \psi_1,\frac{1}{2}p^*F \psi_1+i\nabla \x\rangle\right)\in \Omega^\bullet(V)
$$
where~$F\in \Omega^2(M;\End(V))$ is the curvature 2-form of~$V$, $\x\in \Omega^0(V,p^*V)$ is the canonical section, and $\nabla\x\in \Omega^1(V,p^*V)$ is its covariant derivative. 
\end{lem}

\bp We calculate component fields of $\nabla_D\psi$ on $\Map(\R^{0|1},V)$ as 
\beq
(\nabla_D\psi)_0&=&i_0^*(\nabla_D\psi)=\psi_0,\nonumber\\
(\nabla_D\psi)_1&=&i_0^*(\nabla_D (\nabla_D\psi))=i_0^*\big(\frac{1}{2}(\nabla_D\nabla_D+\nabla_D\nabla_D)\psi\big) =i_0^*\big(\frac{1}{2}(\ev^*p^*F)(D,D)\psi\big)\nonumber\\
&=&\frac{1}{2} i_0^*((\ev^*p^*F)(D,D)) \psi_1\nonumber
\eeq
where $F\in \Omega^2(M;\End(V))\cong \Omega^2(M;\End(\Pi V))$ is the curvature 2-form, $i_0\colon \Map(\R^{0|1},V)\hookrightarrow \R^{0|1}\times \Map(\R^{0|1},V)$ is induced by the map on functions that sends $\theta\in C^\infty(\R^{0|1})$ to zero, and we used $[D,D]=0$. Then performing the Berezinian integral over the fibers~\eqref{eq:intfibers01} (using~$\int \theta [d\theta]=1$) we get
\beq
\SMQ(\psi_1,\psi_0)&=&\int \left(\frac{1}{2}\langle \psi_1+\theta\psi_0,\psi_0+\theta \frac{1}{2}i_0^*((\ev^*p^*F)(D,D))\psi_1\rangle+i\langle \psi_1+\theta\psi_0,\x_0+\theta \x_1 \rangle \right)[d\theta]\nonumber\\
&=&\left(\langle \psi_0,\frac{1}{2}\psi_0+i\x_0\rangle+\langle \psi_1,-\frac{1}{4} i_0^*((\ev^*p^*F)(D,D)) \psi_1-i\x_1\rangle\right). \nonumber
\eeq
The remaining identification is achieved by the next lemma. 
\ep

\begin{lem} \label{lem:curvpull}
Under the isomorphism
$$
\Omega^\bullet(M;\End(V))\cong \Gamma(\Map(\R^{0|1},M);i_0^*\ev^*\End(V))
$$
the section $-\frac{1}{2} i_0^*\ev^*F(D,D)\in \Gamma(\Map(\R^{0|1},M);\End(p^*V))$ is identified with the curvature 2-form $
F\in \Omega^2(M;\End(V)).$
Similarly, under the isomorphisms
$$
\Omega^\bullet(V;p^*\Pi V)\cong \Gamma(\Map(\R^{0|1},V);i_0^*\ev^*p^*\Pi V),\quad \Omega^\bullet(V;p^*V)\cong \Gamma(\Map(\R^{0|1},V);i_0^*\ev^*p^*V),
$$
the sections 
$$
\x_1=i_0^*(\nabla_D\ev^*\x)\in \Gamma(\Map(\R^{0|1},V);i_0^*\ev^*p^*\Pi V),\quad \x_0=i_0^*\x \in\Gamma(\Map(\R^{0|1},V);i_0^*\ev^*p^*V)
$$ are identified with $-\nabla\x\in \Omega^1(V;p^*\Pi V)$ and $\x\in \Omega^0(V;p^*V)$, respectively. 
\end{lem}
\bp 
First we establish a little notation. For $f\in C^\infty(M)=\Omega^0(M)\subset \Omega^\bullet(M)$ and $df\in \Omega^1(M)\subset \Omega^\bullet(M)$, let $f,\delta f\in C^\infty(\Map(\R^{0|1},M))$ denote the associated functions on $\Map(\R^{0|1},M)$. 
Then the evaluation map $\ev\colon \Map(\R^{0|1},M)\times \R^{0|1}\to M$ on functions is
$$
f\mapsto f+\delta f \theta \in C^\infty(\Map(\R^{0|1},M)\times \R^{0|1}),\quad f\in C^\infty(M).
$$
The map $i_0$ in this case is the inclusion $\Map(\R^{0|1},M)\hookrightarrow \Map(\R^{0|1},M)\times\R^{0|1}$ that on functions sets the coordinate~$\theta\in C^\infty(\R^{0|1})$ to zero. 

For a 2-form $F=fdgdh\in \Omega^2(M)$, pull back along evaluation is
$$
\ev^*F=\ev^*(fdgdh)=(f+\delta f\theta)d(g+ \delta g\theta)d(h+ \delta h\theta)\in \Omega^2(\Map(\R^{0|1},M)\times\R^{0|1}),
$$
and so contracting with~$D$ twice and then setting $\theta=0$ we get
\beq
i_0^*((\ev^*F)(D,D))&=&i_0^* (\iota_D\iota_D\ev^*F)=-2f\delta g\delta h\nonumber\\
&=&-2F\in \Omega^0(\Map(\R^{0|1},M))= C^\infty(\Map(\R^{0|1},M)),\nonumber
\eeq
where in the last line (in an abuse of notation) $F$ denotes the corresponding function on $\Map(\R^{0|1},M)$. 
From this (and, e.g., working in a local trivialization of $V$), we yield the claimed formula. The arguments for $\nabla\x$ and $\x$ are similar. 
\ep


Following the usual quantization prescription from physics, we would like to integrate $\exp(-\SMQ(\psi))$ over sections,
$$
\int e^{-\SMQ(\psi)}[d\psi]_0=\int e^{-\SMQ(\psi_1,\psi_0)} d\psi_1d\psi_0=\int e^{-\langle \psi_0,\frac{1}{2}\psi_0+i\x\rangle} d\psi_0 \int e^{- \langle \psi_1, \frac{1}{2}p^*F\psi_1+i\nabla \x \rangle} d\psi_1,
$$
which requires a measure $[d\psi]_0$. Using the identification $\psi_0\in \Omega^\bullet(V,p^*V)$ and $\psi_1\in \Omega^\bullet(V,p^*\Pi V)$, we define such a measure in terms of a pair of measures, $d\psi_0$ and $d\psi_1$.

\begin{defn} We define the integral $\int e^{-\SMQ(\psi)}[d\psi]_0$ as follows. 
Let $\int e^{-\langle \psi_0,\frac{1}{2}\psi_0+i\x\rangle} d\psi_0$ be determined by integration with respect to the real volume form on the fibers of~$p^*V\to V$ inherited from the metric and orientation on~$V$.
Define $\int e^{- \langle \psi_1, \frac{1}{2}p^*F\psi_1+i\nabla \x \rangle} d\psi_1$ in terms of the Berezinian integral in~\cite[\S1.6]{BGV}. Explicitly, the functional is a section of $\Lambda^\bullet p^*V^\vee$, and define the Berezinian integral of such a functional as the projection
$$
\Omega^\bullet(V,\Lambda^\bullet p^*V^\vee)\to \Omega^\bullet(V,\Lambda^{\rm top}V^\vee)\cong \Omega^\bullet(V),
$$
using the orientation on~$V$. 
\end{defn}

\begin{proof}[Proof of Proposition~\ref{prop:01MQ} for $d=0$] \label{d0case}
All modes are zero modes in this case, so we need to perform the integral over the whole of $u^*\Fer(p^*V)$. We evaluate the usual Gaussian integral to get
\beq
\int e^{-\SMQ(\psi)}[d\psi]_0&=&\int e^{-\SMQ(\psi_1,\psi_0)} d\psi_1d\psi_0=\int e^{-\langle \psi_0,\frac{1}{2}\psi_0+i\x\rangle} d\psi_0 \int e^{- \langle \psi_1, \frac{1}{2}p^*F\psi_1+i\nabla \x \rangle} d\psi_1\nonumber \\
&=&(2\pi)^{{\rm dim}(V)/2}e^{-\|\x\|^2/2} \int e^{-\langle \psi_1, \frac{1}{2}p^*F\psi_1+i\nabla \x\rangle} d\psi_1\nonumber
\eeq
where we have used
$$
\int e^{-\langle \psi_0,\frac{1}{2}\psi_0+i\x\rangle}d\psi_0=(2\pi)^{{\rm dim}(V)/2}e^{-\|\x\|^2/2}\in C^\infty(V)=\Omega^0(V)\subset \Omega^\bullet(V). 
$$
The remaining Berezinian integral over the odd parameter space is precisely the Mathai--Quillen Thom form~\eqref{eq:MQThom} as defined in~\cite[\S1.6]{BGV} up to the claimed normalizing factor
$$
u_V=\frac{\epsilon({\rm dim}(V))}{(2\pi)^{{\rm dim}(V)}}\int e^{-\SMQ(\psi)}[d\psi]_0\in \Gamma_{\rm c}(\L^{0|1}_0(V);\omega^{\otimes{\rm dim}(V)/2})\cong \Omega^{{\rm dim}(V)}_{\rm cl,c}(V)
$$
where $\epsilon({\rm dim}(V))=1$ for ${\rm dim}(V)$ even and $i=\sqrt{-1}$ for ${\rm dim}(V)$ odd.
\ep

\subsection{An equality of pushforwards}\label{sec:deRhampush} 

For a family of oriented manifolds, $\pi\colon M\to B$, the topological pushforward, denoted $\widehat{\pi}_!^{\rm top}$, was defined in~\S\ref{sec:IV}. 

\begin{proof}[Proof of Theorem~\ref{thm:3} for $d=0$]

This follows directly from the Thom isomorphism in de~Rham cohomology: the integral of the Thom form of~$\nu$ over the fibers is~1, and the remaining integral for the topological pushforward is over the fibers~$\pi\colon M\to B$. This agrees with the analytic pushforward. \ep

\section{Pushforwards in complexified $\KO$ from $1|1$-dimensional field theories}\label{sec:KO}

In this section we prove the $\KO$ versions of the main theorems. The relevant theory from physics is supersymmetric mechanics, i.e., a theory of maps from $1|1$-dimensional supermanifolds to a Riemannian manifold~$M$. 

\subsection{Super circles, super loop spaces, and complexified K-theory}\label{sec:KC}

Let~$\E^{1|1}$ denote the super Lie group with underlying supermanifold $\R^{1|1}$ and multiplication given by the functor of points formula
\beq
(t,\theta)\cdot (t',\theta')=(t+t'-i\theta\theta',\theta+\theta'), \quad (t,\theta),(t',\theta')\in \R^{1|1}(S).\label{11grplaw}
\eeq
In contrast to the $0|1$-dimensional case, this super Lie group is not abelian. The Lie algebra of left-invariant vector fields on~$\E^{1|1}$ is freely generated (as a Lie algebra) by a single odd generator, $D:=\partial_{\theta}+i\theta\partial_t$. As a vector space, this Lie algebra is spanned by~$D$ and~$\frac{1}{2}[D,D]=i\partial_t$. Next, let $\widetilde{\R}^\times\subset \C^\times$ denote the subgroup generated by the 4th roots of unity and $\R^\times$. There is a double cover $\widetilde{\R}^\times\to \R^\times$ given by $\mu\mapsto \mu^2$ for $\mu\in \widetilde{\R}^\times$, and we also have the isomorphism $\widetilde{\R}^\times\cong \Z/4\times \R_{>0}$. Consider the semidirect product~$\E^{1|1}\rtimes \widetilde{\R}^\times$ for the action $\mu\cdot (t,\theta)= (\mu^2 t,\mu \theta)$, where $(t,\theta)\in\E^{1|1}(S)$ and $\mu \in\widetilde{\R}^\times (S)$. The constructions below will make repeated use of the obvious left action of~$\E^{1|1}\rtimes \widetilde{\R}^\times$ on~$\R^{1|1}$; the notation~$\E^{1|1}$ and~$\R^{1|1}$ distinguishes between the Lie group and the supermanifold on which it acts, respectively. Similarly, let~$\E$ denote the Lie group whose underlying manifold is~$\R$ equipped with the usual additive structure. We observe that the canonical inclusion~$\R\hookrightarrow \R^{1|1}$ of the reduced manifold gives a homomorphism~$\E\hookrightarrow \E^{1|1}$. 

An \emph{$S$-family of (1-dimensional, based) lattices} is an $S$-family of monomorphisms $\ell \colon S\times \Z\to \E$. We often notationally identify a lattice~$\ell$ with its generator~$\ell\in \R^\times(S)\subset \E(S)$ gotten by restricting to $S\cong S\times\{1\}\hookrightarrow S\times \Z\to \E$. 
Given $\ell\in \R^\times(S)$, the canonical inclusion $\E\subset \E^{1|1}$ gives an action map
$$
 S\times\R^{1|1}\times \Z\stackrel{\ell}{\hookrightarrow} S\times \R^{1|1}\times \E\hookrightarrow S\times \R^{1|1}\times \E^{1|1}\stackrel{\rm act}{\to} S\times \R^{1|1}.
$$
This action is free, and the quotient supermanifold is a \emph{family of super circles}, denoted
$$
S^{1|1}_\ell:=(S\times \R^{1|1})/\Z.
$$ 
The quotient map defines a fiberwise (universal) cover, 
\beq
S\times \R^{1|1}\to S^{1|1}_\ell.\label{eq:circlecov}
\eeq
The following definition formalizes the notion of a map~$S^{1|1}_\ell\to S^{1|1}_{\ell'}$ that is locally determined by the action of $\E^{1|1}\rtimes \widetilde{\R}^\times$ on $\R^{1|1}$ using~\eqref{eq:circlecov}. 

\begin{defn} \label{rmk:bundleofgrps}A \emph{fiberwise rigid conformal map} between families of super circles is a map $S^{1|1}_\ell\to S^{1|1}_{\ell'}$ over a base change $S\to S'$  for which there exists a commutative diagram
\beq
\begin{tikzpicture}[baseline=(basepoint)];
\node (A) at (0,0) {$S\times \R^{1|1}$};
\node (B) at (4,0) {$S^{1|1}_\ell$};
\node (C) at (0,-1.5) {$S'\times \R^{1|1}$};
\node (D) at (4,-1.5) {$S^{1|1}_{\ell'}$}; 
\draw[->>] (A) to  (B);
\draw[->,dashed] (A) to  (C);
\draw[->>] (C) to (D);
\draw[->] (B) to (D);
\path (0,-.75) coordinate (basepoint);
\end{tikzpicture}\nonumber
\eeq
where the horizontal arrows are the quotient maps~\eqref{eq:circlecov}. We require that the dashed arrow be determined by a base change $S\to S'$ together with a map
$$
S\times \R^{1|1}\to \E^{1|1}\rtimes \widetilde{\R}^\times \times \R^{1|1}\to \R^{1|1}
$$
where the first arrow is given by an $S$-point of $\E^{1|1}\rtimes \widetilde{\R}^\times$, and the second arrow is the left action of this super Lie group on $\R^{1|1}$. We further require that the dashed arrow be a $\Z$-equivariant map for the $\Z$-action defining the families of super circles, relative to an $S$-family of homomorphisms $S\times \Z\to S\times \Z$ determined by an $S$-point of $\GL_1(\Z)$. 
\end{defn}


\begin{rmk} A super conformal map $\R^{1|1}\to \R^{1|1}$ is a diffeomorphism that preserves the distribution generated by the odd vector field $D:=\partial_\theta+i\theta\partial_t$, e.g., see~\cite[pg.~52]{5lectures}. The group of rigid conformal transformations $\E^{1|1}\rtimes \widetilde{\R}^\times$ is a strict subgroup of the super conformal transformation group. 
\end{rmk}

\begin{rmk}\label{rmk:11HST}
The above definition of a fiberwise rigid conformal map is a special case of (and was inspired by) fiberwise isometries between $S$-families with rigid conformal structure in the sense of Stolz and Teichner's rigid geometries~\cite[\S4.2]{ST11}; see also~\cite[\S6.3]{HST}. 
\end{rmk}

We observe that the existence of the diagram in Definition~\ref{rmk:bundleofgrps} means that the data of a fiberwise rigid conformal map over~$S$ is a section of the bundle of groups $(S\times \E^{1|1}\rtimes \widetilde{\R}^\times)/\Z\times \GL_1(\Z)\to S$, for the quotient by the fiberwise subgroup $S\times \Z\stackrel{\ell}{\hookrightarrow}S\times \E\hookrightarrow S\times \E^{1|1}\hookrightarrow S\times \E^{1|1}\rtimes \widetilde{\R}^\times$. In turn, such a section determines an $S$-point of $\widetilde{\R}^\times$, an $S$-point of $\GL_1(\Z)$, and a section of $(S\times \E^{1|1})/\Z\to S$. As in the previous section, we refer to the $S$-point of~$\widetilde{\R}^\times$ as the \emph{dilation} part of the fiberwise rigid conformal map and the section of $(S\times \E^{1|1})/\Z\to S$ as the \emph{super translation} part.

The following is the specialization of Definitions~\ref{defn:fields} and~\ref{defn:factor} to the $d=1$ case.

\begin{defn}\label{defn:superloop} The \emph{super loop stack of $M$}, denoted $\L^{1|1}(M)$, is the stack associated to the prestack whose objects over~$S$ are pairs $(\ell,\phi)$ where $\ell\in \R^\times(S)$ determines a family of super circles~$S^{1|1}_\ell$ and $\phi\colon S^{1|1}_\ell\to M$ is a map. Morphisms between these objects over a base change~$S\to S'$ consist of commuting triangles
\beq
\begin{tikzpicture}[baseline=(basepoint)];
\node (A) at (0,0) {$S^{1|1}_\ell$};
\node (B) at (3,0) {$S^{1|1}_{\ell'}$};
\node (C) at (1.5,-1.5) {$M$};
\draw[->] (A) to node [above=1pt] {$\cong$} (B);
\draw[->] (A) to node [left=1pt]{$\phi$} (C);
\draw[->] (B) to node [right=1pt]{$\phi'$} (C);
\path (0,-.75) coordinate (basepoint);
\end{tikzpicture}\label{11triangle}
\eeq
where the horizontal arrow is a fiberwise rigid conformal map. The stack of \emph{constant super loops}, denoted $\L^{1|1}_0(M)$, is the full substack for which~$\phi$ factors as
\beq
S^{1|1}_\ell\to S\times \R^{0|1}\to M. \label{eq:11factorization}
\eeq
For a map $M\to M'$, post-composition $S^{1|1}_\ell\to M\to M'$ defines a morphism of stacks $\L^{1|1}(M)\to \L^{1|1}(M')$ and $\L^{1|1}_0(M)\to \L^{1|1}_0(M')$. 
\end{defn}

\begin{rmk} The objects of $\L^{1|1}(M)$ are fields for supersymmetric mechanics with target~$M$, or equivalently, the $1|1$-dimensional $\sigma$-model from Example~\ref{ex:sigma}. 
\end{rmk}

\begin{rmk}
There is a more general notion of an $S$-family of super circles corresponding to $S$-families of based lattices $S\times \Z\hookrightarrow S\times \E^{1|1}$ inside the non-commutative super Lie group~$\E^{1|1}$. This moduli space is interesting, but doesn't fit nicely with the story from physics below; see Remark~\ref{eq:superlatices} for a discussion in the $2|1$-dimensional case.
\end{rmk}

\begin{rmk} \label{rmk:KUKO}We note that $i=\sqrt{-1}\in \widetilde{\R}^\times\subset \C^\times$ acts by orientation reversal on a family of super circles so that $\L^{1|1}(M)$ is an unoriented version of the super loop stack of~$M$. Replacing the rigid conformal group $\E^{1|1}\rtimes \widetilde{\R}^\times$ by the subgroup $\E^{1|1}\rtimes \R^\times$ results in an oriented version. Analogous constructions to the ones below using the oriented super loop stack yield a cocycle model for~${\rm KU}\otimes \C$ rather than $\KO\otimes \C$. \end{rmk}

Define a morphism of stacks $\L^{1|1}_0(\pt)\to [\pt\sq \C^\times]$ that is constant to~$\pt$ on objects over~$S$ and to a rigid conformal map $S^{1|1}_\ell\to S^{1|1}_{\ell'}$ assigns the $S$-point of $\widetilde{\R}^\times\subset \C^\times$ determined by the dilation part. As in the $d=0$ case, this description in terms of $S$-points of the prestack defining $\L^{1|1}_0(\pt)$ suffices to define a morphism of stacks; see Remark~\ref{rmk:stackad}. 

\begin{defn} Define line bundles $\omega^{\otimes k/2}$ on $\L^{1|1}_0(M)$ as the pullback of $k$th tensor power of the dual of the odd canonical line bundle over $[\pt\sq \C^\times]$ along the composition,
$$
\L^{1|1}_0(M)\to \L^{1|1}_0(\pt) \to [\pt\sq \C^\times],
$$
where the first arrow is the functor associated with the smooth map~$M\to \pt$. 
\end{defn}

Our next task is to establish a groupoid presentation of $\L^{1|1}_0(M)$ for the purposes of computing the space of sections of $\omega^{\otimes k/2}$. Consider the functor
\beq
\R^\times\times \Map(\R^{0|1},M)\to \L^{1|1}_0(M)\label{eq:11atlas}
\eeq
that sends an $S$-point $\ell\in \R^\times(S)$, $\phi_0\in \Map(\R^{0|1},M)(S)$ of the source to
\beq
\phi\colon S^{1|1}_\ell\stackrel{pr}{\to} S\times \R^{0|1}\stackrel{\phi_0}{\to} M\label{eq:proj11}
\eeq
where~$pr$ is determined by the standard projection $\R^{1|1}\to \R^{0|1}$. 

\begin{prop}
The map~\eqref{eq:11atlas} defines an atlas for the stack whose associated Lie groupoid presentation is 
$$
\left[ \begin{array}{c} (\E^{1|1}\rtimes \widetilde{\R}^\times\times  \R^\times)/\Z \times \GL_1(\Z)\times \Map(\R^{0|1},M) \\ \downarrow \downarrow \\ \R^\times\times \Map(\R^{0|1},M)\end{array} \right] \stackrel{\sim}{\to} \L^{1|1}_0(M).
$$
The source map is the projection. The target map is the composition
\beq
(\E^{1|1}\rtimes \widetilde{\R}^\times\times  \R^\times)/\Z \times \GL_1(\Z)\times \Map(\R^{0|1},M))&\to& \E^{0|1}\rtimes \widetilde{\R}^\times\times  \R^\times \times \GL_1(\Z)\times  \Map(\R^{0|1},M))\nonumber\\
&\to& \R^\times\times \Map(\R^{0|1},M)\nonumber
\eeq
where the first map comes from the projection $pr\colon (\E^{1|1}\rtimes \widetilde{\R}^\times \times \R^\times)/\Z\to \E^{0|1}\rtimes \widetilde{\R}^\times \times \R^\times$, and the second map is determined by the left action of $\E^{0|1}\rtimes \widetilde{\R}^\times<\E^{0|1}\rtimes \C^\times$ on $\Map(\R^{0|1},M)$ by precomposition, and the action of $\widetilde{\R}^\times\times \GL_1(\Z)$ on $\R^\times$ given by
$$
(\mu,\pm 1)\cdot \ell=\pm \mu^2\ell\qquad \mu\in \widetilde{\R}^\times(S), \ \pm 1\in \GL_1(\Z)(S), \ \ell\in \R^\times(S).
$$
\label{prop:11present}
\end{prop}

\bp
The stack $\L^{1|1}_0(M)$ is defined as the stackification of a prestack we shall denote by $\L^{1|1}_0(M)_{\rm pre}$ in this proof. We will show that $\L^{1|1}_0(M)_{\rm pre}$ is equivalent to the prestack defined by the Lie groupoid in the proposition. This will imply that their stackifications are equivalent. 

To show that a map of prestacks is an equivalence, it suffices to demonstrate an equivalence of groupoids at each $S$-point. 
To verify essential surjectivity at an $S$-point, the prestack $\L^{1|1}_0(M)_{\rm pre}$ has objects over~$S$ given by~$(\ell,\phi)$ for~$\ell\in \R^\times(S)$ and $\phi\colon S^{1|1}_\ell\to M$. By definition, there is a unique $\phi_0\colon S\times \R^{0|1}\to M$ such that $\phi=\phi_0\circ pr$ as in~\eqref{eq:proj11}. Hence, the map~\eqref{eq:11atlas} induces an essential surjection on objects over~$S$. We will extend this to a map from the appropriate Lie groupoid by characterizing morphisms in~$\L^{1|1}_0(M)_{\rm pre}$ over~$S$. These morphisms were described explicitly just before Definition~\ref{defn:superloop}. With $S=\R^\times\times \Map(\R^{0|1},M)$, the supermanifold of morphisms over~$S$ is~$(\E^{1|1}\rtimes \widetilde{\R}^\times\times  \R^\times)/\Z \times \GL_1(\Z)\times \Map(\R^{0|1},M)$. Therefore, morphisms over a general $S$ are $S$-points of $(\E^{1|1}\rtimes \widetilde{\R}^\times\times  \R^\times)/\Z \times \GL_1(\Z)\times \Map(\R^{0|1},M)$. This gives a map from the prestack underlying the Lie groupoid in the statement of the proposition to the prestack $\L^{1|1}_0(M)_{\rm pre}$, which (by construction) induces an equivalence of groupoids at each $S$-point. 

It remains to verify the claimed characterizations of the source and target maps in the Lie groupoid presentation. The source map is clearly the projection. For $(t,\theta,\mu,\pm 1)\in (\E^{1|1}\rtimes \widetilde{\R}^\times\times \GL_1(\Z))(S)$ determining a rigid conformal map with source $(\ell,\phi_0)\in (\R^\times\times \Map(\R^{0|1},M))(S)$, the target has $\pm \mu^2\ell\in \R^\times(S)$ and the rigid conformal map acts by precomposition on~$\phi=\phi_0\circ pr$. The latter action is through the projection to $\E^{0|1}\rtimes \widetilde{\R}^\times$, which is exactly the precomposition action on $\phi_0$ claimed in the statement. 
\ep

\begin{cor}
The map 
\beq
u\colon \Map(\R^{0|1},M)\hookrightarrow \R^\times\times \Map(\R^{0|1},M)\to \L^{1|1}_0(M)\label{eq:11atlas2}
\eeq
where the first arrow includes along $1\in \R^\times$ and the second is~\eqref{eq:11atlas} defines an atlas for the stack whose associated Lie groupoid is a quotient Lie groupoid,
$$
\left[\Map(\R^{0|1},M)\sq (\E^{1|1}/\Z)\rtimes \Z/4\right] \stackrel{\sim}{\to} \L^{1|1}_0(M),
$$
for the action given by the homomorphism
\beq
(\E^{1|1}/\Z)\rtimes \Z/4\to \E^{0|1}\rtimes \Z/4
\eeq
followed by the standard action of $\E^{0|1}\rtimes \Z/4<\E^{0|1}\rtimes \C^\times$ on $\Map(\R^{0|1},M)$ by precomposition. 
\label{cor:11present}
\end{cor}

\bp
The inclusion of the Lie groupoids
$$
\Map(\R^{0|1},M)\sq (\E^{1|1}/\Z)\rtimes \Z/4 \hookrightarrow \left\{ \begin{array}{c} (\E^{1|1}\rtimes \widetilde{\R}^\times\times  \R^\times)/\Z \times \GL_1(\Z) \times \Map(\R^{0|1},M) \\ \downarrow \downarrow \\ \R^\times\times \Map(\R^{0|1},M)\end{array} \right\}
$$
along $1\in \R^\times$ is fully faithful and essentially surjective. Essential surjectivity follows from the fact that the $\widetilde{\R}^\times$-action on $\R^\times$ is transitive. Fullness and faithfulness follow from the fact that automorphisms of a family of super circles with~$\ell=1$ are given precisely by $S$-points of~$(\E^{1|1}/\Z)\rtimes \Z/4$. Since the inclusion of Lie groupoids is an equivalence, the underlying stacks are equivalent. 
\ep

For future reference, we denote the (evaluation) map from the family of super circles mapping to~$M$ over the objects in the Lie groupoid of Corollary~\ref{cor:11present} by~$\ev$,
\beq
\ev\colon \Map(\R^{0|1},M)\times \R^{1|1}/\Z\to \Map(\R^{0|1},M)\times \R^{0|1} \to M\label{eq:11ev}
\eeq
where the first arrow is the projection, and the second arrow is the usual evaluation map. 

The homomorphism
$$
\E^{1|1}/\Z\rtimes \Z/4 \twoheadrightarrow \Z/4 \subset \C^\times
$$
determines a functor from the Lie groupoid in Corollary~\ref{cor:11present} to the Lie groupoid $\pt\sq \C^\times$. We can pullback the dual of the canonical odd line bundle along this functor, obtaining a line bundle over the groupoid that is isomorphic to the pullback of~$\omega^{1/2}$ from $\L^{1|1}_0(M)$. This allows us to compute sections in terms of functions on the cover with transformation properties.

\begin{proof}[Proof of Theorem~\ref{thm1}, $d=1$ case]  
By Corollary~\ref{cor:11present} and~\eqref{compsections}, the atlas~\eqref{eq:11atlas2} allows one to compute the sections $\Gamma(\L^{1|1}_0(M);\omega^{\otimes k/2})$ as the subspace of functions
$$
\Gamma(\L^{1|1}_0(M);\omega^{\otimes k/2})\subset C^\infty(\Map(\R^{0|1},M))\cong \Omega^\bullet(M)
$$ 
that are equivariant for the action by~$(\E^{1|1}/\Z)\rtimes \Z/4$. This action factors through the homomorphism
$$
\E^{1|1}\rtimes \Z/4\twoheadrightarrow \E^{0|1}\rtimes \Z/4.
$$
In view of the cocycle defining~$\omega^{\otimes k/2}$, sections are functions on $\Map(\R^{0|1},M)$ that are $\E^{0|1}$-invariant and for which the action by the generator $\sqrt{-1}\in \Z/4\subset \C$ is by~$(\sqrt{-1})^k$. Identifying functions on $\Map(\R^{0|1},M)$ with differential forms, this gives
$$
\Gamma(\L^{1|1}_0(M);\omega^{\otimes k/2})\cong \Omega^{4\bullet+k}_\cl(X).
$$
Assembling these vector spaces in each degree into a graded algebra we obtain the desired isomorphism of graded algebras
$$
\Gamma(\L^{1|1}_0(M);\omega^{\otimes \bullet/2})\cong \Omega^\bullet_{\cl}(M;\C[\beta,\beta^{-1}]),
$$
where~$\beta^{-1}$ is the trivializing section of $\omega^{\otimes 4/2}$ giving the (Bott) isomorphism 
$$
\Gamma(\L^{1|1}_0(M);\omega^{\otimes \bullet/2})\stackrel{\otimes \beta^{-1}}{\to} \Gamma(\L^{1|1}_0(M);\omega^{\otimes (\bullet+4)/2}).
$$
Naturality of~$\L^{1|1}_0(M)$ in~$M$ turns $\Gamma(\L^{1|1}_0(M);\omega^{\otimes \bullet/2})$ into a presheaf. We have shown that this presheaf is isomorphic to the sheaf of 4-periodic closed differential forms, so in fact $\Gamma(\L^{1|1}_0(M);\omega^{\otimes \bullet/2})$ is a sheaf, and the isomorphism is one of sheaves. 
\ep

\begin{rmk} We could also compute the space of sections $\Gamma(\L^{1|1}_0(M);\omega^{\otimes k/2})$ using the atlas from Proposition~\ref{prop:11present}. In this case, sections of~$\omega^{\otimes k/2}$ are generated by functions of the form (see~\eqref{eq:HST})
$$
\beta^{({\rm deg}(\alpha)+k)/4}\otimes \alpha \in C^\infty(\R^\times)\otimes \Omega^\bullet(M)\cong C^\infty(\R^\times\times \Map(\R^{0|1},M))
$$ 
where $\alpha\in \Omega^\bullet(M)$, the $k={\rm deg}(\alpha)\mod 4$, and $\beta=1/\ell^2\in C^\infty(\R^\times)$. In this case the Bott periodicity isomorphism is a bit easier to see, arising explicitly through the dependence on~$\ell$, or equivalently, $\beta$, viewed as a trivializing section of $\omega^{\otimes 4/2}$. When we restrict to the atlas in Corollary~\ref{cor:11present},~$\beta=1$ is the trivializing section and so it no longer enters the formulas explicitly. 
\end{rmk}

As is common in the study of loop spaces, an important aspect of the geometry of~$\L^{1|1}_0(M)$ is a circle action by loop rotation. Let $u\colon  \Map(\R^{0|1},M)\to \L^{1|1}_0(M)$ be the atlas from Corollary~\ref{cor:11present}. A vector bundle $W\to \L^{1|1}_0(M)$ when pulled back along $u$ will have a loop rotation action. Said differently, we have a factorization of~$u$ as
$$
u\colon \Map(\R^{0|1},M)\to [\Map(\R^{0|1},M)\sq S^1]\to [\Map(\R^{0|1},M)\sq \E^{1|1}/\Z\rtimes \Z/4]\simeq  \L^{1|1}_0(M)
$$
where the $S^1=\E/\Z$-action on $\Map(\R^{0|1},M)$ is the trivial action. So the pullback of a vector bundle over $\L^{1|1}_0(M)$ to $[\Map(\R^{0|1},M)\sq S^1]$ is a vector bundle on $\Map(\R^{0|1},M)$ equipped with an $S^1$-action on its fibers. This permits the following definition. 

\begin{defn} For $W\to \L^{1|1}_0(M)$ a vector bundle and $k\in \Z$, the \emph{$k$th weight space} is the subbundle of $u^*W$ on which the action by $S^1=\E^1/\Z<(\E^{1|1}/\Z)\rtimes \Z/4$ is multiplication by the function
$$
e^{2\pi i s k}\in C^\infty(\E^1/\Z\times \Map(\R^{0|1},M))
$$
for $s$ the standard coordinate on $\E^1$. The \emph{zero modes of $W$} are sections of weight zero, and the \emph{nonzero modes} lie in the span of sections with nonzero weight. 
\end{defn}

The utility of the above definition comes from the fact that in key examples the weight spaces break an infinite-rank vector bundle into a sum of finite-rank ones. 

\subsection{Vector bundles on $\L^{1|1}_0(M)$}

We refer to~\S\ref{sec:01MQ} for preliminaries on how to promote a vector bundle $V\to M$ with its connection and metric to a vector bundle in supermanifolds with a connection and a $C^\infty(M)$-bilinear pairing on sections. We recall that $\mathcal{V}$ denotes the sheaf of sections of~$V$ over~$M$ regarded as a supermanifold. 

\begin{rmk}
The following definitions implicitly use the adjunction~\eqref{eq:leftad} between the forgetful functor from stacks to prestacks and stackification: we will define a vector bundle on the prestack whose stackification is $\L^{1|1}_0(M)$, as this is equivalent to defining a vector bundle on the stack $\L^{1|1}_0(M)$. 
\end{rmk}

\begin{defn} \label{defn11ferbos}
For $V\to M$ a real oriented vector bundle, define a vector bundle $\Bos(V)\to \L^{1|1}_0(M)$ whose sections at an $S$-point $(\ell,\phi)$ is the $C^\infty(S)$-module 
$$
\Gamma(S^{1|1}_\ell,\phi^*\Pi V)=C^\infty(S^{1|1}_\ell)\otimes_{C^\infty(M)} \mathcal{V}
$$
where $\otimes$ is the projective tensor product of Fr\'echet spaces. Similarly define a vector bundle $\Fer(V)\to \L^{1|1}_0(M)$ whose sections at an $S$-point $(\ell,\phi)$ is the $C^\infty(S)$-module 
$$
\Gamma(S^{1|1}_\ell,\phi^*\Pi V)=C^\infty(S^{1|1}_\ell)\otimes_{C^\infty(M)}\Pi \mathcal{V},
$$
noting the parity reversal of the module $\mathcal{V}$ in this second definition. In terms of a sheaf of sections, these are the the direct image sheaves along the projection $S^{1|1}_\ell\to S$ of the sheaf of sections of the pullback $\phi^*V$ and $\phi^*\Pi V$ respectively. For an isomorphism in $\L^{1|1}_0(M)$ associated to a fiberwise rigid conformal map $f\colon S^{1|1}_\ell\to S^{1|1}_{\ell'}$ over a base change $S\to S'$, define maps of vector bundles from the pullback of sections along~$f$,
$$
\Gamma(S^{1|1}_{\ell'},\phi'^*V)\stackrel{f^*}{\to} \Gamma(S^{1|1}_\ell,\phi^*V)\stackrel{\mu\cdot }{\to} \Gamma(S^{1|1}_\ell,\phi^*V)
$$
$$
\Gamma(S^{1|1}_{\ell'},\phi'^*\Pi V)\stackrel{f^*}{\to} \Gamma(S^{1|1}_\ell,\phi^*\Pi V)
$$
where $\mu\in \R^\times(S)\subset C^\infty(S)^\ev$ is the dilation part of a fiberwise rigid conformal map. 
\end{defn}

\begin{rmk} Sections of $\Bos(V)$ and $\Fer(V)$ are fields for $V$-valued free bosons and $V$-valued free fermions, respectively, whence the notation. \end{rmk}

From Corollary~\ref{cor:11present}, the vector bundles $\Bos(V)$ and $\Fer(V)$ can be characterized in terms of vector bundles over the Lie groupoid that presents $\L^{1|1}_0(M)$. This consists of vector bundles $u^*\Bos(V)$ and $u^*\Fer(V)$ over $\Map(\R^{0|1},M)$ along with isomorphisms ${\sf s}^*u^*\Bos(V)\to {\sf t}^*u^*\Bos(V)$ and ${\sf s}^*u^*\Fer(V)\to {\sf t}^*u^*\Fer(V)$ of vector bundles, running in complete analogy to the $0|1$-dimensional case in~\S\ref{sec:01MQ}. Explicitly, sections of these vector bundles are the $C^\infty(\Map(\R^{0|1},M))$-modules
\beq
\Gamma(\Map(\R^{0|1},M);u^*\Bos(V))&=&\Gamma(\Map(\R^{0|1},M)\times \R^{1|1}/\Z;\ev^* V)\nonumber \\
\Gamma(\Map(\R^{0|1},M);u^*\Fer(V))&=&\Gamma(\Map(\R^{0|1},M)\times \R^{1|1}/\Z;\ev^*\Pi V)\nonumber 
\eeq
for $\ev$ the evaluation map~\eqref{eq:11ev}. In contrast to the $0|1$-dimensional case, $u^*\Fer(V)$ is no longer a finite-rank $C^\infty(\Map(\R^{0|1},M))$-module. 

In computations we will use \emph{component fields}, which for sections 
$$
X\in \Gamma(\Map(\R^{0|1},M),u^*\Bos(V))\qquad \psi\in \Gamma(\Map(\R^{0|1},M),u^*\Fer(V)),
$$ 
are defined using analogous formulas to before
\beq
\begin{array}{lll} 
a&:=&i_0^*X\in \Gamma(\Map(\R^{0|1},M)\times\R/\Z,i_0^*\ev^* V) \\ 
\eta&:=&i_0^* (\nabla_DX)\in \Gamma(\Map(\R^{0|1},M)\times\R/\Z,i_0^*\ev^*\Pi V)\\
\psi_1&:=&i_0^*\psi\in \Gamma(\Map(\R^{0|1},M)\times\R/\Z,i_0^*\ev^*\Pi V) \\ 
\psi_0&:=&i_0^* (\nabla_D\psi)\in \Gamma(\Map(\R^{0|1},M)\times\R/\Z,i_0^*\ev^*V).
\end{array} \label{eq:component11}
\eeq
In the above, $D=\partial_\theta+i\theta\partial_t$, $i_0\colon \Map(\R^{0|1},M)\times\R/\Z\hookrightarrow \Map(\R^{0|1},M)\times\R^{1|1}/\Z$ is the fiberwise inclusion of the reduced manifold, and $\ev$ is the evaluation map~\eqref{eq:11ev}. This gives isomorphisms of $C^\infty(\Map(\R^{0|1},M))$-modules 
\beq
&&\Gamma\left(\Map(\R^{0|1},M)\times \R^{1|1}/\Z,\ev^*\Pi V\right)\nonumber \\
&&\stackrel{\sim}{\to} \Gamma\left(\Map(\R^{0|1},M)\times\R/\Z,i_0^*\ev^*\Pi V\right)\oplus \Gamma\left(\Map(\R^{0|1},M)\times\R/\Z,i_0^*\ev^*V\right)\nonumber\\
&&\Gamma\left(\Map(\R^{0|1},M)\times \R^{1|1}/\Z,\ev^*V\right)\nonumber \\
&&\stackrel{\sim}{\to} \Gamma\left(\Map(\R^{0|1},M)\times\R/\Z,i_0^*\ev^* V\right)\oplus \Gamma\left(\Map(\R^{0|1},M)\times\R/\Z,i_0^*\ev^*\Pi V\right)\nonumber
\eeq
sending $\psi\mapsto (\psi_1,\psi_0)$ and $X\mapsto (a,\eta)$, respectively. The inverse isomorphisms are
$$
(\psi_1,\psi_0)\mapsto {\rm pr}^*\psi_1+\theta {\rm pr}^*\psi_0\qquad (a,\eta)\mapsto {\rm pr}^*a+\theta {\rm pr}^*\eta
$$
$$
{\rm pr}\colon \Map(\R^{0|1},M)\times\R^{1|1}/\Z\to \Map(\R^{0|1},M)\times\R/\Z
$$ 
where ${\rm pr}$ is the obvious projection. We caution that this projection is not invariant under super translations, and these Taylor components are similarly not invariant.  To simplify notation we often omit the pullback along the projection, i.e., writing~$\psi=\psi_1+\theta \psi_0$ and $X=a+\theta \eta$. For $p^*V\to V$ the pullback of~$V$ over itself, we have the canonical section~$\x$ that determines the section~$\ev^*\x\in \Gamma(\Map(\R^{0|1},V),u^*\Fer(p^*V))$ with components
\beq
&&\begin{array}{lll} \x_0&:=&i_0^*\ev^*\x \in \Gamma(\Map(\R^{0|1},V)\times\R/\Z,i_0^*\ev^*p^*V) \\ 
\x_1&:=&i_0^* (\nabla_D\ev^*\x)\in \Gamma(\Map(\R^{0|1},V)\times\R/\Z,i_0^*\ev^*p^*\Pi V),\end{array} \label{eq:component112}
\eeq
and we will similarly write $\x=\x_0+\theta\x_1$.

As before, component fields have a description in terms of differential forms since the composition
$$
\Map(\R^{0|1},M)\times\R/\Z\stackrel{i_0}{\hookrightarrow} \Map(\R^{0|1},M)\times\R^{1|1}/\Z\stackrel{\ev}{\to} M
$$
equals the composition of projections $\Map(\R^{0|1},M)\times\R/\Z\to \Map(\R^{0|1},M)\to M$. Hence, for $\psi\in \Gamma(\Map(\R^{0|1},M),u^*\Fer(V))$
\beq
&&\psi_0\in C^\infty(\R/\Z)\otimes \Omega^\bullet(M,V),\qquad \psi_1\in C^\infty(\R/\Z)\otimes \Omega^\bullet(M,\Pi V),\label{eq:comp11mod}
\eeq 
and for $X\in \Gamma(\Map(\R^{0|1},M),u^*\Bos(V))$
\beq
&&a\in C^\infty(\R/\Z)\otimes \Omega^\bullet(M,V),\qquad \eta\in C^\infty(\R/\Z)\otimes \Omega^\bullet(M,\Pi V),\label{eq:comp11modA}
\eeq
where the tensor product is the projective tensor product of Fr\'echet spaces. 

As a special case of Definition~\ref{defn11ferbos} above, we obtain vector bundles on $\L^{1|1}_0(M)$ from families of oriented manifolds $\pi\colon M\to B$. Indeed, consider~$\BosT \to \L^{1|1}_0(M)$ whose space of sections at an $S$-point $(\ell,\phi)$ is the $C^\infty(S)$-module
$$
\Gamma(S^{1|1}_\ell;\phi^*T(M/B))=C^\infty(S^{1|1}_\ell)\otimes_{C^\infty(M)} \Gamma(T(M/B)).
$$
Geometrically, when $B=\pt$ the nonzero modes of $\BosT$ are the normal bundle to the inclusion $\L^{1|1}_0(M)\hookrightarrow \L^{1|1}(M)$. For general~$B$ these sections are a relative version of such a normal bundle with respect to the projection $\L^{d|1}(M)\to \L^{d|1}(B)$.

\subsection{The analytic pushforward} \label{sec:familiesAhat}

We recall the action functional $\mathcal{S}_\sigma$ from~\eqref{eq:superLag} with quadratic part~\eqref{eq:sigmahess}. By applying this quadratic part to sections of $u^*\BosT$ we obtain the action functional
\beq
\begin{array}{c}\displaystyle
\HessS(X)=\int \langle X,\nabla_{i\partial_t}\nabla_DX\rangle [d\theta dt],\\ 
X\in \Gamma(\Map(\R^{0|1},M)\times \R^{1|1}/\Z,\ev^*T(M/B))\end{array}\label{eq:11HessS}
\eeq
where $\nabla$ is the pullback of the connection on~$T(M/B)$, $\langle-,-\rangle$ is the pullback of the metric, and the integral is with respect to the Berezinian measure $[d\theta dt]$ on the fibers of the projection
\beq
\Map(\R^{0|1},M)\times \R^{1|1}/\Z\to \Map(\R^{0|1},M).\label{eq:intfibers111}
\eeq


\begin{lem} The functional~\eqref{eq:11HessS} defined on objects descends to a functional defined on the Lie groupoid: 
$$
{\sf s}^*\HessS={\sf t}^*\HessS
$$
for ${\sf s},{\sf t}$ the source and target maps in the Lie groupoid presentation of $\L^{1|1}_0(M)$ in Corollary~\ref{cor:11present}, and where the equality uses the isomorphism of vector bundles~${\sf s}^*u^*\BosT\cong {\sf t}^*u^*\BosT$. Hence, $\HessS$ descends to the stack $\L^{1|1}_0(M)$.
\end{lem}

\bp 
This is a similar check as in the $0|1$-dimensional case. First, we observe that the integrand defining $\HessS$ is built from left-invariant vector fields, so is invariant under the left action by super translations (we recall that generators of the left action are \emph{right}-invariant vector fields, and these commute with the left-invariant vector fields). Since the Berezinian measure is also invariant, $\HessS$ is invariant under super translations. It remains to analyze behavior of~$\HessS$ under the action by the 4th roots of unity $\Z/4\subset \C^\times$. The pullbacks of the constituent pieces defining $\HessS$ along this action are: $X\mapsto \mu X$, $\nabla_D\mapsto \mu^{-1}\nabla_D$, $\nabla_{\partial_t}\mapsto \mu^{-2}\nabla_{\partial_t}$, and $[d\theta dt]\mapsto \mu [d\theta dt]$ for $\mu\in \Z/4$. These effects cancel, and so $\HessS$ is invariant and therefore descends to the Lie groupoid and hence the stack $\L^{1|1}_0(M)$.
\ep

\begin{lem} \label{lem:11Hesscomp}
The value of $\HessS$ on component fields $(a,\eta)$ is 
\beq
\HessS(a,\eta)&=&\int\Big(\langle a,(-\nabla_{\partial_t}^2-\nabla_{i\partial_t}R)a\rangle+\langle\eta,\nabla_{i\partial_t}\eta\rangle \Big)dt,\nonumber
\eeq
where $\nabla$ is the pullback of the connection on $T(M/B)$ along $\ev\circ i_0$,  $R\in \Omega^2(M,\End(T(M/B)))$ is the curvature 2-form of the fiberwise Levi-Civita connection acting on component fields using the identification~\eqref{eq:comp11modA}, and the integral is over the fibers
$$
\Map(\R^{0|1},M)\times\R/\Z\to \Map(\R^{0|1},M).
$$

\end{lem}

\bp Let $\nabla$ denote the pullback of the connection on $T(M/B)$ along $\ev$; note that $\nabla$ in the statement of the lemma is $i_0^*\nabla$ in the notation in the proof. We compute the Taylor components as in (\ref{eq:comp11modA}) of the section~$\nabla_{i\partial_t}\nabla_DX$:
\beq
i^*_0(\nabla_{i\partial_t}\nabla_DX)&=&(i_0^*\nabla)_{i\partial_t} i^*_0 (\nabla_DX)=i_0^*(\nabla)_{i\partial_t} \eta\nonumber\\
i^* _0\nabla_D(\nabla_{i\partial_t}\nabla_DX)&=&i^*_0(\nabla_{i\partial_t}\nabla_D\nabla_D X)=\frac{1}{2}i^*_0(\nabla_{i\partial_t} (\nabla_{[D,D]}+\ev^*R(D,D))X)\nonumber\\
&=&i_0^*(-\nabla_{\partial_t}^2+\nabla_{i\partial_t}\frac{1}{2}\ev^*R(D,D)))i^*_0X\nonumber\\
&=&(-i_0^*(\nabla)_{\partial_t}^2+\frac{1}{2}i_0^*(\nabla)_{i\partial_t}i_0^*(\ev^*R(D,D)))a\nonumber
\eeq
where we used 
$$
\nabla_D^2=\frac{1}{2}(\nabla_D\nabla_D+\nabla_D\nabla_D)=\frac{1}{2}(\nabla_{[D,D]}+\ev^*R(D,D))=\nabla_{i\partial_t}+\frac{1}{2}\ev^*R(D,D),
$$
$\iota_{\partial_t} \ev^*R=0$ and $[\partial_t,D]=0$.  
The result now follows from the same arguments as in the proof of Lemma~\ref{lem:calccomp01}: compute the Berezinian integral using $\int \theta [d\theta]=1$, observe that the evaluation map~\eqref{eq:11ev} factors as
$$
\Map(\R^{0|1},M)\times\R^{1|1}/\Z\to \Map(\R^{0|1},M)\times \R^{0|1}\to M.
$$
where the first arrow is projection and the second is the evaluation map from the $0|1$-dimensional case, and apply Lemma~\ref{lem:curvpull}. 
\ep

Our regularization procedure requires we also consider a complementary bundle $\Fer(n):=\Fer(\underline{\R}^n)$ from Definiton~\ref{defn11ferbos} for $\underline{\R}^n$ the trivial bundle. Equipping it with the trivial connection, we get a family of operators~$D$ acting on sections of $u^*\Fer(n)$ and an associated functional on sections
$$
\mathcal{S}_\Fer(\psi)=\frac{1}{2}\int \langle \psi,D\psi\rangle [d\theta dt]
$$
$$
\psi\in \Gamma(\Map(\R^{0|1},M)\times \R^{1|1}/\Z,\Pi \underline{\R}^n),
$$
where $\langle-,-\rangle$ is the standard metric on $\R^n$ and the integral is with respect to the Berezinian measure $[d\theta dt]$ on the fibers of the usual projection. Writing a section in terms of component fields, $\psi=\psi_1+\theta \psi_0$ and integrating with respect to~$\theta$, we find
$$
\mathcal{S}_\Fer(\psi_1,\psi_0)=\frac{1}{2}\int \left(\langle \psi_0,\psi_0\rangle-\langle \psi_1,i\partial_t\psi_1\rangle\right) dt.
$$
This formula is an easy case of Proposition~\ref{prop:SMQ11} below. 

\begin{notation} Let $\Bos_l(T(M/B))$ denote the $l$th weight space of the vector bundle $u^*\BosT$ on $\Map(\R^{0|1},M)$. Denote the decomposition into even and odd parts of this vector bundle as~$\Bos_l(T(M/B))= \Bos_l^\ev(T(M/B))\oplus \Bos_l^\odd(T(M/B))$.  Similarly, let $\Fer_l(n)$ denote the $l$th weight space of the vector bundle $u^*\Fer(n)$ on $\Map(\R^{0|1},M)$. Denote the decomposition into even and odd parts of this vector bundle as~$\Fer_l(n)=\Fer_l^\ev(n)\oplus\Fer_l^\odd(n)$.
\end{notation}
 
\begin{proof}[Proof of Proposition \ref{prop:Aclass}]
By Definition~\ref{defn:21sigmareg}, we need to calculate the finite-dimensional super determinants
\beq
\sdet(\nabla_{i\partial_t}\nabla_D|_{\Bos_l(T(M/B))})&=&\frac{\det(\nabla_{i\partial_t}\nabla_D|_{\Bos_l^\odd(T(M/B))})}{\det(\nabla_{i\partial_t}\nabla_D|_{\Bos_l^\ev(T(M/B))}}
=\frac{\det(i\nabla_{\partial_t}|_{\Bos_l^\odd(T(M/B))})}{\det((-\nabla_{\partial_t}^2-\nabla_{i\partial_t}R)|_{\Bos_l^\ev(T(M/B))})}\nonumber\\
\sdet(D|_{\Fer_l(n)})&=&\frac{\det(i\partial_t|_{\Fer_l^\odd(n)})}{\det(\id |_{\Fer_l^\ev(n)})}\nonumber
\eeq
where the last equality in the first line is from Lemma~\ref{lem:11Hesscomp}. We compute these determinants locally on~$U\subset M$ using a local basis for weight~$l$ component fields
$$
e^{-2\pi i lt}\otimes a_j\in C^\infty(\R/\Z)\otimes \Omega^\bullet(U,T(M/B)),
$$
$$ e^{-2\pi i lt}\otimes \eta_j\in C^\infty(\R/\Z)\otimes \Omega^\bullet(U,\Pi T(M/B)),
$$
where $\{a_j\}$ (respectively, $\{\eta_j\}$) form a local basis of sections of~$T(M/B)$ (respectively, $\Pi T(M/B)$) on~$U\subset M$. Similarly we have a global basis 
$$
e^{-2\pi i lt}\otimes \psi_j^0\in C^\infty(\R/\Z)\otimes \Omega^\bullet(M,\R^n),$$
$$ e^{-2\pi i lt}\otimes \psi_j^1\in C^\infty(\R/\Z)\otimes \Omega^\bullet(M,\Pi \R^n),
$$ 
where $\psi_j^0$ is a basis for $\R^n$ and $\psi_j^1$ is a basis for $\Pi \R^n$. 
In terms of the former basis, $\nabla_{\partial_t} =d/dt\otimes {\rm id}$, so we have
$$
(\nabla_{i\partial_t}\nabla_D)|_{\Bos_l^\ev(T(M/B))}=-\frac{d^2}{d t^2}\otimes \id-i\frac{d}{dt} \otimes R,\quad\quad (\nabla_{i\partial_t}\nabla_D)|_{\Bos_l^\odd(T(M/B))}=i\frac{d}{dt}\otimes \id.
$$
So then on the $l$th weight space, 
\beq
\det\left(i\frac{d}{d t}\otimes \id\right)&=&\det\left(2\pi l\otimes \id\right)=\left(2\pi l\right)^{n}\nonumber 
\eeq
and 
\beq
\det\left(-\frac{d^2}{d t^2}\otimes \id-i\frac{d}{dt} \otimes R\right)&=&\det\left(4\pi^2l^2\otimes \id+2\pi  l\otimes R\right)\nonumber\\
&=&\det\left(4\pi^2l^2\otimes \id\right)\det\left(\id\otimes \id+\frac{1}{2\pi l}\otimes R\right)\nonumber \\
&=&\left(2\pi l\right)^{2n}\det\left(\id\otimes \id+\frac{1}{2\pi l}\otimes R\right)\nonumber
\eeq
where $n={\rm dim}(M)-{\rm dim}(B)$. On $\Fer_l(n)$ we have 
$$
\frac{\det(i\partial_t|_{\Fer_l^\odd(n)})}{\det(\id |_{\Fer_l^\ev(n)})}=(2\pi l)^n.
$$
Putting this together, the contribution to the regularized determinant from the $l$th weight space is 
\beq
\det(\id\otimes \id+\frac{1}{2\pi l}\otimes R)&=&\exp({\rm Tr}(\log(\id\otimes \id+\frac{1}{2\pi l}\otimes R)))\nonumber\\
&=&\exp\left({\rm Tr}\left(\sum_k \frac{(-1)^{k+1}}{k}\frac{1}{(2\pi l)^k}R^k\right)\right)\nonumber \\
&=&\exp\left({\rm Tr}\left(-\sum_k \frac{1}{2k}\frac{1}{(2\pi l)^{2k}}R^{2k}\right)\right)\nonumber
\eeq
where in the last equality we used that traces of odd powers of the curvature vanish. 
The product over $l>0$ is the exponential of the absolutely convergent series
$$
\sum_{l>0 }{\rm Tr}\left(\sum_k  \frac{1}{2k}\frac{1}{(2\pi l)^{2k}}R^{2k}\right)=\sum_{k=1}^\infty\frac{{\rm Tr}(R^{2k}) }{2k (2\pi)^{2k}}\zeta(2k)
$$
where $\zeta(k)$ denotes the value of the Riemann $\zeta$-function at~$k$. Putting this together we get
$$
\sdet_{\rm reg}'(\nabla_{i\partial_t}\nabla_D)^{-1/2}=\exp\left(\sum_{k=1}^\infty\frac{{\rm Tr}(R^{2k})}{2k (2\pi)^{2k}}\zeta(2k)\right).
$$
It remains to compare with the $\hat{A}$-form. For a real vector bundle $V$, let ${\rm ph}_{k}(V)$ denote the degree $4k$ part of the Pontryagin character of $V$, i.e., ${\rm ph}_k(V)={\rm ch}_{2k}(V\otimes \C)$ where ${\rm ch}_k$ is the degree~$2k$ piece of the Chern character. In our cocycle model, 
\beq
{\rm Tr}(R^{2k})=(-1)^k2(2k)! {\rm ph}_{k}(V),\label{eq:charclass}
\eeq
where ${\rm ph}_{k}$ denotes the $4k^{\rm th}$ component of the Pontryagin character as a function on~$\L^{1|1}_0(M)$. Hence
\beq
\sdet_{\rm reg}'(\nabla_{i\partial_t}\nabla_D)^{-1/2}&=&\exp\left(\sum_{k=1}^\infty \frac{(2k)! {\rm ph}_{k}(T(M/B))}{2k} \frac{2\zeta(2k)}{(2\pi i)^{2k}}\right)\nonumber
\eeq
which we identify as a cocycle representative of the $\hat{A}$-class of the family $\pi \colon M\to B$ as a function on~$\L^{1|1}_0(M)$. The function descends to the stack. We conclude that the analytic pushforward in this case is given by the formula
$$
\widehat{\pi}_!^{\rm an}(f)=\int_{M/B} f\cdot\sdet_{\rm reg}'(\nabla_{i\partial_t}\nabla_D)^{-1/2}= \int_{M/B} f\cdot \hat{A}(M/B),
$$
which is a cocycle refinement of the complexification of the pushforward in~$\KO$.
\ep

\subsection{The Mathai--Quillen form}\label{sec:KQM}

We specialize the vector bundle~$\Fer(V)$ to the vector bundle $p^*V\to V$ gotten from pulling $V$ back over itself. There is a functional on sections of $u^*\Fer(p^*V)$ given by essentially the same formula as in the $0|1$-dimensional case
\beq
\SMQ(\psi)=\int \left(\frac{1}{2}\langle \psi,\nabla_D\psi\rangle+i\langle \psi,\phi^*\x\rangle\right)[d\theta dt]\label{eq:11MQaction}
\eeq
except now the odd vector field is $D=\partial_\theta+i\theta\partial_t$ and the integral is over the fibers of the bundle of super circles,
\beq
\Map(\R^{0|1},V)\times\R^{1|1}/\Z\to \Map(\R^{0|1},V)\label{eq:intfibers11}
\eeq
using the fiberwise Berezinian measure descending from~$[d\theta dt]$ on $\R^{1|1}$. 

\begin{rmk}
The Mathai--Quillen action functional~\eqref{eq:11MQaction} is invariant under super translations, but its behavior under the action of $\Z/4\subset \C^\times$ is more subtle. On the physics side, this has to do with various subtleties related to time-reversal symmetry (and lack thereof). 
Our convention for the vector bundle $\Fer(V)$ is set up so that the quadratic part of the classical action~\eqref{eq:11MQaction} descends to the stack; see Lemma~\ref{lem:11perp}. 
\end{rmk}

The value of the Mathai--Quillen classical action on component fields is an element of~$\Omega^\bullet(V)$, which we calculate in terms of more familiar quantities.

\begin{prop} \label{prop:SMQ11}
The value of $\SMQ$ on component fields $(\psi_1,\psi_0)$ is
$$
\SMQ(\psi_1,\psi_0)=\int \left( \langle \psi_0,\frac{1}{2}\psi_0+i\x \rangle +\langle \psi_1,\frac{1}{2}(-\nabla_{i\partial_t}+p^*F)\psi_1+i\nabla \x \rangle \right)dt\in  \Omega^\bullet(V)
$$
where $\nabla$ is the pullback of the connection on $p^*V$ along $\ev\circ i_0$, $F\in \Omega^2(M;\End(V))$ is the curvature 2-form of~$V$, $\x\in \Omega^0(V,p^*V)$ is the canonical section, $\nabla\x\in \Omega^1(V,p^*V)$ is its covariant derivative, and the integral is over the fibers,
\beq
\Map(\R^{0|1},V)\times\R/\Z\to \Map(\R^{0|1},V),\label{eq:11compintegral}
\eeq
using the fiberwise measure~$dt$.
\end{prop}

\begin{proof} Let $\nabla$ denote the pullback of the connection on $T(M/B)$ along $\ev$; note that $\nabla$ in the statement of the proposition is $i_0^*\nabla$ in the notation in the proof. We calculate the components as 
\beq\label{eq:components}
\begin{array}{llll}
\phantom{BLA}&(\nabla_D\psi)_0&=&i_0^*\nabla_D\psi=\psi_0\\
&(\nabla_D\psi)_1&=&i_0^*(\nabla_D\nabla_D\psi)=i_0^*(\frac{1}{2}(\nabla_D\nabla_D+\nabla_D\nabla_D)\psi)\\
&&=&i_0^*(\frac{1}{2}(\ev^*F(D,D)+\nabla_{[D,D]})\psi) =(\frac{1}{2}i_0^*(\ev^*F(D,D))+(i_0^*\nabla)_{i\partial_t})\psi_1. \\
&(\phi^*\x)_0&=&i_0^*\ev^*\x\\
&(\phi^*\x)_1&=&i_0^*(\nabla_D \ev^*\x).
\end{array}
\eeq
Then the result follows from the same arguments as in the proof of Lemma~\ref{lem:calccomp01}: compute the Berezinian integral using $\int \theta [d\theta]=1$, observe that the evaluation map~\eqref{eq:11ev} factors as
$$
\Map(\R^{0|1},M)\times\R^{1|1}/\Z\to \Map(\R^{0|1},M)\times \R^{0|1}\to M,
$$
where the first arrow is projection and the second is the evaluation map from the $0|1$-dimensional case, and finally apply Lemma~\ref{lem:curvpull}. 
\ep

Now we set to work on providing rigorous meaning for the integral $\int e^{-\SMQ(\psi)}[d\psi]$ over sections, as per~\eqref{MQsplit}. 
We start with defining the integral over the zero modes. In terms of component fields, these zero modes are the subspaces
$$
\psi_0\in \Omega^\bullet(V,p^*V)\hookrightarrow C^\infty(\R/\Z)\otimes \Omega^\bullet(V,p^*V),$$$$\psi_1\in \Omega^\bullet(V,p^*\Pi V)\hookrightarrow C^\infty(\R/\Z)\otimes \Omega^\bullet(V,p^*\Pi V)
$$
of sections that are constant in the fiber direction. Hence, we define the integral over zero modes as the one in~\S\ref{sec:01MQmeasure} for the $0|1$-dimensional case: the integral of the functional over the sections~$\psi_0$ is with respect to the fiberwise (real) volume form on the fibers~$p^*V\to V$ coming from the metric and orientation on~$V$, and the integral of the functional over the sections~$\psi_1$ is the Berezinian integral 
$$
\Omega^\bullet(V,\Lambda^\bullet V^\vee)\to \Omega^\bullet(V,\Lambda^{\rm top} V^\vee)\cong \Omega^\bullet(V),
$$
which uses the orientation on~$V$.

\begin{proof}[Proof of Lemma~\ref{lem:d1zeromodes} when $d=1$]
Evaluating the Mathai--Quillen classical action on a zero mode, and applying Proposition~\ref{prop:SMQ11}, we obtain
\beq
\SMQ(\psi_1,\psi_0)&=&\int \left(\langle \psi_0,\frac{1}{2}\psi_0+i \x\rangle +\langle \psi_1,\frac{1}{2}(-i\nabla_{\partial_t}+p^*F)\psi_1+i \nabla\x\rangle\right)dt\nonumber\\
&=&\left(\langle \psi_0,\frac{1}{2}\psi_0+i \x\rangle+\langle \psi_1,\frac{1}{2} p^*F\psi_1+i \nabla\x\rangle\right)\in C^\infty(S)\nonumber
\eeq
where we used $\nabla_{\partial_t}\psi_1=0$ on zero modes (since $\partial_t$ generates the circle action defining the weight spaces) and the integral is over the fibers~\eqref{eq:intfibers11} so the integral evaluates simply as the fiberwise volume (namely,~$1$) multiplied by the integrand. So we are left to evaluate 
$$
\int e^{-\SMQ(\psi)}[d\psi]_0=\int \exp(-\langle \psi_0,\frac{1}{2}\psi_0+i \x\rangle)d \psi_0\int \exp(- \langle \psi_1,\frac{1}{2}p^*F\psi_1+i \nabla\x\rangle)d\psi_1
$$
as a function on $\Map(\R^{0|1},V)$. The claim now follows immediately from the case where~$d=0$ proved in~\S\ref{d0case}, and the normalizing factor is $\frac{\epsilon({\rm dim}(V))}{(2\pi)^{{\rm dim}(V)}}$ as before. 
\ep

Next we unpack $\SMQ$ on the nonzero modes.

\begin{lem}\label{lem:11perp} For a nonzero mode $\psi$ of $u^*\Fer(p^*V)$, we have 
\beq
\SMQ(\psi)=\int\frac{1}{2} \langle \psi,\nabla_D\psi\rangle [d\theta dt],\label{eq:11perp1}
\eeq
where the integral is over the fibers~\eqref{eq:intfibers11} and the operator $\nabla_D$ is invertible. Furthermore, the functional~\eqref{eq:11perp1} pulls back along the map
$$\Map(\R^{0|1},V)\to \Map(\R^{0|1},M)\to \L^{1|1}_0(M),
$$ 
meaning the formula on the atlas descends to a functional on sections over the stack. 

\end{lem}

\bp 
From~\eqref{eq:comp11mod}, weight $n$ component fields are spanned by ones of the form
$$
\psi_0=e^{2\pi i nt}\otimes v_0\in C^\infty(\R/\Z)\otimes \Omega^\bullet(V,p^*V), 
$$$$
\psi_1=e^{2\pi i nt}\otimes v_1\in C^\infty(\R/\Z)\otimes \Omega^\bullet(V,p^*\Pi V). 
$$
Similarly, $(\phi^*\x)_0$ and $(\phi^*\x)_1$ are $1\otimes \x$ and $1\otimes \nabla \x$ in this description. Then we evaluate the term $\langle \psi,\phi^*\x\rangle$ in~$\SMQ$ on these component fields as 
\beq
\int \langle \psi_0,(\phi^*\x)_0\rangle dt&=&\int e^{2\pi i nt}  dt \cdot \langle v_0,\x\rangle=0\nonumber\\
\int \langle \psi_1,(\phi^*\x)_1\rangle dt &=&\int e^{2\pi i nt}  dt\cdot \langle v_1, \nabla \x\rangle=0\nonumber
\eeq
where the integrals are over the fibers~\eqref{eq:11compintegral}. Hence, for nonzero modes this term in the classical action vanishes, verifying~\eqref{eq:11perp1}. We observe that $-i\nabla_{\partial_t}$ is an invertible operator on the nonzero modes. So by Proposition~\ref{prop:SMQ11}, $\nabla_D=-i\nabla_{\partial_t}+F$ is a nilpotent modification of this invertible operator (the curvature 2-form is nilpotent), so is invertible. 

To verify the remainder of the lemma, since $\nabla$ in~\eqref{eq:11perp1} pulls back from $V\to M$ along $\ev^*\circ p^*$, we get that $\SMQ$ on nonzero modes does indeed pullback from a functional on $u^*\Fer(V)$. Finally, we need to check that the functional $\int\frac{1}{2} \langle \psi,\nabla_D\psi\rangle [d\theta dt]$ on sections of $u^*\Fer(V)$ descends to a functional on sections of $\Fer(V)$ over the stack $\L^{1|1}_0(M)$. This is a standard check of invariance under isomorphisms in the Lie groupoid presentation from Corollary~\ref{cor:11present}. The integrand is built from left-invariant vector fields, so is invariant under the left action by super translations. It remains to analyze behavior under the action by the 4th roots of unity $\Z/4\subset \C^\times$. The pullbacks of the constituent pieces along the $\Z/4$-action map are: $[d\theta dt]\mapsto \frac{1}{\mu} [d\theta] \mu^2 [dt]=\mu [d\theta dt]$ and $\nabla_D\mapsto \mu^{-1}\nabla_D$. Assembling the functional from these pieces, we see that it is indeed invariant. 
\ep

Our regularization procedure requires we also consider a complementary bundle $\Bos(m):=\Bos(\underline{\R}^m)$ from Definiton~\ref{defn11ferbos} for $\underline{\R}^m$ the trivial bundle. Equipping it with the trivial connection, we get a family of operators~$\Delta=Di\partial_t$ acting on sections of $u^*\Bos(m)$ and an associated functional on sections
$$
\mathcal{S}_\Bos(X)=\frac{1}{2}\int \langle X,\Delta X\rangle [d\theta dt]
$$
$$
X\in \Gamma(\Map(\R^{0|1},M)\times \R^{1|1}/\Z,\underline{\R}^m),
$$
where $\langle-,-\rangle$ is the standard metric on $\R^m$ and the integral is with respect to the Berezinian measure $[d\theta dt]$ on the fibers of the usual projection. Writing a section in terms of component fields, $X=a+\theta \eta$ and integrating with respect to~$\theta$, as a special case of Lemma~\ref{lem:11Hesscomp} we find
$$
\mathcal{S}_\Bos(a,\eta)=\frac{1}{2}\int \left(-\langle a,\partial_t^2 a\rangle+\langle \eta,i\partial_t\eta\rangle\right)  dt.
$$

\begin{notation} Let $\Fer_l(V)$ denote the $l$th weight space of the vector bundle $u^*\Fer(V)$ on $\Map(\R^{0|1},M)$. Denote the decomposition into even and odd parts of this vector bundle as~$\Fer_l(V)\cong \Fer^\ev_l(V)\oplus \Fer^\odd_l(V)$.  Similarly, let $\Bos_l(m)$ denote the $l$th weight space of the vector bundle $u^*\Bos(m)$ on $\Map(\R^{0|1},M)$. Denote the decomposition into even and odd parts of this vector bundle as~$\Bos_l(m)=\Bos_l^\ev(m)\oplus\Bos_l^\odd(m)$.
\end{notation}

\begin{proof}[Proof of Proposition~\ref{thm:thomkthy}] 
By Lemma~\ref{lem:d1zeromodes} and~\eqref{MQsplit}, it remains to identify the regularized super determinant of $\nabla_D$ (from Definition~\ref{defn:21MQreg}) with the Riemann--Roch factor $p^*\hat{A}(V)^{-1}$. 
So we need to calculate the finite-dimensional super determinants
$$
\sdet(\nabla_D|_{\Fer_l(V)})= \frac{\det(\nabla_D|_{\Fer_l^\odd(V)})}{\det(\nabla_D|_{\Fer_l^\ev(V)})}=\frac{\det(-i\nabla_{\partial_t}+F)}{\det(\id)}
$$
$$
\sdet(\Delta|_{\Bos_l(m)})=\frac{\det(i\partial_t|_{\Bos_l^\odd(m)})}{\det(\partial_t^2|_{\Bos_l^\ev(m)})}
$$
where the last equality in the top equation is from Lemma~\ref{lem:11perp}. 

We begin by expressing the finite-dimensional determinants locally on~$M$ in a local basis for weight~$l$ component fields
$$
e^{2\pi i lt}\otimes v_j\in C^\infty(\R/\Z)\otimes \Omega^\bullet(M,\Pi V),
$$
where $\{v_j\}$ form a local basis of sections of~$\Pi V$ over~$M$ and
$$
e^{2\pi i lt}\otimes a_j\in C^\infty(\R/\Z)\otimes \Omega^\bullet(M,\R^m)\qquad e^{2\pi i lt}\otimes \eta_j\in C^\infty(\R/\Z)\otimes \Omega^\bullet(M,\Pi \R^n)
$$
for $\{a_j\}$ a basis of $\R^m$ and $\{\eta_j\}$ a basis of $\Pi \R^m$. In terms of the basis $\{v_j\}$, we have~$\nabla_D|_{\Fer_l^\odd(V)}=-i\partial_t\otimes \id+\id\otimes F$ and hence
\beq
\det(-i\nabla_{\partial_t}+F)&=&\det(2\pi l \otimes \id+\id\otimes F)=\det(2\pi l\otimes \id)\det(\id\otimes \id+\frac{1}{2\pi l}\otimes F)\nonumber \\
&=&\left(2\pi l\right)^m\det(\id\otimes \id+\frac{1}{2\pi l}\otimes F)\nonumber
\eeq
where $m={\rm dim}(V)$. We also have
$$
\sdet(\Delta|_{\Bos_l(m)})=\frac{\left(2\pi l\right)^m}{\left(2\pi l\right)^{2m}}=\frac{1}{\left(2\pi l\right)^m}.
$$
So the contribution to the regularized determinant from the $l$th weight space is 
\beq
\det(\id\otimes \id+\frac{1}{2\pi l}\otimes F)&=&\exp({\rm Tr}(\log(\id\otimes \id+\frac{1}{2\pi l}\otimes F)))\nonumber\\
&=&\exp\left({\rm Tr}\left(\sum_k \frac{(-1)^{k+1}}{k}\frac{1}{(2\pi l)^k}F^k\right)\right)\nonumber \\
&=&\exp\left({\rm Tr}\left(-\sum_k \frac{1}{2k}\frac{1}{(2\pi l)^{2k}}F^{2k}\right)\right)\nonumber
\eeq
where in the last equality we used that traces of odd powers of the curvature vanish. The product over weight spaces is the exponential of the absolutely convergent series
$$
\sum_{l> 0}{\rm Tr}\left(-\sum_k  \frac{1}{2k}\frac{1}{(2\pi l)^{2k}}F^{2k}\right)=-\sum_{k=1}^\infty\frac{{\rm Tr}(F^{2k}) }{2k (2\pi)^{2k}}\zeta(2k)
$$
where $\zeta(k)$ denotes the value of the Riemann $\zeta$-function at~$k$. So together we have
$$
\sdet_{\rm reg}'(\nabla_D)^{-1/2}=p^*\left(\exp\left(-\sum_{k=1}^\infty\frac{{\rm Tr}(F^{2k}) }{2k (2\pi )^{2k}}\zeta(2k)\right)\right).
$$
Using~\eqref{eq:charclass}, we get
$$
\sdet_{\rm reg}'(\nabla_D)^{-1/2}=p^*\left(\exp\left(-\sum_{k=1}^\infty \frac{(2k)! {\rm ph}_{k}(V)}{2k} \frac{2\zeta(2k)}{(2\pi i)^{2k}}\right)\right)\nonumber
$$
which we identify with $\hat{A}(V)^{-1}$ as an element of $C^\infty(\L^{1|1}_0(M))$.

This regularized super determinant combined with Lemma~\ref{lem:d1zeromodes} shows that 
$$
\frac{\epsilon({\rm dim}(V))}{(2\pi)^{{\rm dim}(V)}} \sdet_{\rm reg}'(\nabla_D)^{-1/2} \cdot \int \exp(-\SMQ)[d\psi]_0  \in \Gamma_{\rm c}(\L^{1|1}_0(V);\omega^{\otimes{\rm dim}(V)/2})
$$ 
is the Mathai--Quillen cocycle representative for the K-theoretic Thom class.
\ep

\subsection{An index theorem}\label{sec:Kindex}

\begin{proof}[Proof of Theorem~\ref{thm:3}, $d=1$ case]
As discussed in Remark~\ref{rmk:complexifiedindex11}, the equality of analytic and topological pushforwards boils down to an equality of functions
\beq
&&\hat{A}(\nu)^{-1}=\hat{A}(M/B)\in C^\infty(\L^{1|1}_0(M)),\nonumber
\eeq
where the left-hand side is the Riemann--Roch factor that modifies the Thom cocycle in de~Rham cohomology, and the right-hand side is the Riemann--Roch factor that modifies integration of differential forms. The equalities of the components of the Pontryagin character~${\rm ph}_k(\nu)=-{\rm ph}_k(T(M/B))$ verify the equality of these functions. This implies that the pushforwards are equal. 
\ep

\section{Pushforwards in complexified TMF from $2|1$-dimensional field theories}\label{sec:TMF}

\subsection{Super tori and fields for the $2|1$-dimensional $\sigma$-model}\label{sec:TMFC}

Let $\E^{2|1}$ denote the super Lie group with underlying supermanifold $\R^{2|1}$ and multiplication 
$$
(z,\bar z,\theta)\cdot (z',\bar z',\theta')=(z+z',\bar z+\bar z'+\theta\theta',\theta+\theta'), \quad (z,\bar z,\theta),(z',\bar z',\theta')\in \R^{2|1}(S),
$$
where we have identified $S$-points of $\R^{2|1}$ with 
\beq
\R^{2|1}(S)\simeq \{z,\bar z\in C^\infty(S)^\ev, \theta\in C^\infty(S)^\odd \mid \overline{(z)}_{\rm red}=(\bar z)_{\rm red}\}. \label{eq:R21Spts}
\eeq
We emphasize that $z$ and $\bar z$ are only conjugates after restriction to the reduced manifold $S_{\rm red}\hookrightarrow S$; in other words, on $S$ the functions $z$ and $\bar z$ are only complex conjugates up to nilpotents. We refer to Example~\ref{eq:CSpoints} for further discussion. The Lie algebra of left-invariant vector fields on~$\E^{2|1}$ has an even generator $\partial_z$ and an odd generator $D=\partial_\theta-\theta\partial_{\bar z}$ with $[\partial_z,D]=0$ and $\frac{1}{2}[D,D]=-\partial_{\bar z}$. Form the semidirect product $\E^{2|1}\rtimes \C^\times$ using the action $(\mu,\bar\mu)\cdot (z,\bar z,\theta)=(\mu^2 z,\bar\mu^2\bar z, \bar\mu \theta)$ for $(\mu,\bar\mu)\in \C^\times(S)$ and $(z,\bar z,\theta)\in \E^{2|1}(S)$, where we have used notation similar to~\eqref{eq:R21Spts} for $S$-points of $\C^\times\subset \C\simeq \R^2$. The constructions below make repeated use of the obvious left action of~$\E^{2|1}\rtimes \C^\times$ on~$\R^{2|1}$; the notation~$\R^{2|1}$ and~$\E^{2|1}$ distinguishes between the Lie group and the supermanifold on which it acts. Similarly, let $\E^2$ denote the Lie group whose underlying manifold is $\R^2$ equipped with the usual additive structure. We observe that the canonical inclusion $\E^2\hookrightarrow \E^{2|1}$ is a group homomorphism. 

A \emph{family of $2$-dimensional based lattices} is an $S$-family of monomorphisms $\ell \colon S\times \Z^2\to S\times \E^2$ so that the ratio of the images of $S\times\{1,0\}$ and $S\times \{0,1\}$ under $\ell\colon S\times \Z^2\to S\times \E^2\cong S\times \C$ are in $S\times (\C\setminus \R)\subset S\times \C$; equivalently, these generators are linearly independent over~$\R$. A \emph{family of $2$-dimensional based, oriented lattices} is an $S$-family of monomorphisms $\ell \colon S\times \Z^2\to S\times \E^2$ such that the ratio of the images of $S\times\{1,0\}$ and $S\times \{0,1\}$ under $\ell\colon S\times \Z^2\to S\times \E^2\cong S\times \C$ are in $S\times \mathbb{H}\subset S\times \C$ for $\mathbb{H}\subset \C$ the upper half plane. Let $\widetilde{\Lat}$ denote the presheaf on supermanifolds whose $S$-points are based lattices and $\Lat$ the presheaf on supermanifolds whose $S$-points are based, oriented lattices. We denote an $S$-point of $\Lat$ or $\widetilde{\Lat}$ by $\ell$ when emphasizing the map $S\times \Z^2\to S\times \E^2$, or $(\ell_1,\bar\ell_1,\ell_2,\bar\ell_2)$ when emphasizing the images of the generators $S\times\{1,0\}$ and $S\times\{0,1\}$ in~$\C$. The complex conjugate notation for the $S$-points follows the same discussion as~\eqref{eq:R21Spts}, e.g., $\ell_1,\bar\ell_1\in C^\infty(S)^\ev$ are complex conjugates only after restriction to $S_{\rm red}$. We observe that $\Lat$ is representable, $\Lat\stackrel{\sim}{\to} \mathbb{H}\times \C^\times$ given on $S$-points by $(\ell_1,\bar\ell_1,\ell_2,\bar\ell_2)\mapsto (\ell_2/\ell_1,\bar\ell_2/\bar\ell_1,\ell_1,\bar\ell_1)$. 

An $S$-family of based lattices defines a free action
$$
 S\times\R^{2|1}\times \Z^2\stackrel{\ell}{\hookrightarrow} S\times \R^{2|1}\times \E^2\hookrightarrow S\times \R^{2|1}\times \E^{2|1}\stackrel{\rm act}{\to} S\times \R^{2|1}.
$$
whose quotient supermanifold is a \emph{family of super tori}
$$
T^{2|1}_\ell:=(S\times \R^{2|1})/\Z^2,
$$
The quotient map defines a cover, 
\beq
S\times \R^{2|1}\to T^{2|1}_\ell.\label{eq:toruscov}
\eeq
 The following definition formalizes the notion of a map~$T^{2|1}_\ell\to T^{2|1}_{\ell'}$ that is locally determined by the action of $\E^{2|1}\rtimes \C^\times$ on $\R^{2|1}$ using~\eqref{eq:toruscov}. 

\begin{defn} \label{defn:21fibconf} A \emph{fiberwise rigid conformal map} between families of super tori is a map $T^{2|1}_\ell\to T^{2|1}_{\ell'}$ over a base change $S\to S'$ for which there exists a commutative  square
\beq
\begin{tikzpicture}[baseline=(basepoint)];
\node (A) at (0,0) {$S\times \R^{2|1}$};
\node (B) at (4,0) {$T^{2|1}_\ell$};
\node (C) at (0,-1.5) {$S'\times \R^{2|1}$};
\node (D) at (4,-1.5) {$T^{2|1}_{\ell'}$}; 
\draw[->] (A) to  (B);
\draw[->,dashed] (A) to  (C);
\draw[->] (C) to (D);
\draw[->] (B) to (D);
\path (0,-.75) coordinate (basepoint);
\end{tikzpicture}\nonumber
\eeq
where the horizontal arrows are the covers~\eqref{eq:toruscov}, and the dashed arrow is determined by a base change $S\to S'$ together with a map
$$
S\times \R^{2|1}\to \E^{2|1}\rtimes \C^\times \times \R^{2|1}\to \R^{2|1}
$$
where the first arrow is given by an $S$-point of $\E^{2|1}\rtimes \C^\times$, and the second arrow is the left action of this super Lie group on $\R^{2|1}$. We further require that the dashed arrow be a $\Z^2$-equivariant map for the $\Z^2$-action defining the families of super tori, relative to an $S$-family of homomorphisms $S\times \Z^2\to S\times \Z^2$ determined by an $S$-point of $\GL_2(\Z)$. 
\end{defn}

We observe that the data of a fiberwise rigid conformal map is a section of the bundle of groups $(S\times \E^{2|1}\rtimes \C^\times\times \GL_2(\Z))/\Z^2\to S$, for the quotient by the fiberwise subgroup $S\times \Z^2\stackrel{\ell}{\hookrightarrow}S\times \E^2\hookrightarrow S\times \E^{2|1}\hookrightarrow S\times \E^{2|1}\rtimes \C^\times\times \GL_2(\Z)$. In turn, such a section determines an $S$-point of $\C^\times$, a section of $(S\times \E^{2|1})/\Z^2\to S$, and an $S$-point of $\GL_2(\Z)$. 
In analogy to the previous cases, we refer to $(\mu,\bar\mu)\in \C^\times(S)$ as the \emph{dilation} part of the fiberwise rigid conformal map and the section of $(S\times \E^{2|1})/\Z^2\to S$ as the \emph{super translation} part of the fiberwise rigid conformal map.

The following is the specialization of Definitions~\ref{defn:fields} and~\ref{defn:factor} to the $d=2$ case.

\begin{defn} \label{defn:supertori} The stack of \emph{super tori in~$M$}, denoted $\L^{2|1}(M)$ is the stack associated to the prestack whose objects over~$S$ are pairs $(\ell,\phi)$ where $\ell\in \widetilde{\Lat}(S)$ determines a family of super tori~$T^{2|1}_\ell$ and $\phi\colon T^{2|1}_\ell\to M$ is a map. Morphisms over a base change consist of commuting triangles
\beq
\begin{tikzpicture}[baseline=(basepoint)];
\node (A) at (0,0) {$T^{2|1}_\ell$};
\node (B) at (3,0) {$T^{2|1}_{\ell'}$};
\node (C) at (1.5,-1.5) {$M$};
\draw[->] (A) to node [above=1pt] {$\cong$} (B);
\draw[->] (A) to node [left=1pt]{$\phi$} (C);
\draw[->] (B) to node [right=1pt]{$\phi'$} (C);
\path (0,-.75) coordinate (basepoint);
\end{tikzpicture}\label{21triangle}
\eeq
where the horizontal arrow is a fiberwise rigid conformal map. The stack of \emph{constant super tori}, denoted $\L{}^{2|1}_0(M)$, is the full substack with objects $(\ell,\phi)$ where $\phi$ factors as
$$
T^{2|1}_\ell\to S\times \R^{0|1}\to M
$$
where the first arrow is induced by the projection $\R^{2|1}\to \R^{0|1}$. 
For a map $M\to M'$, post-composition $T^{2|1}_\ell\to M\to M'$ defines a morphism of stacks $\L^{2|1}(M)\to \L^{2|1}(M')$ and $\L^{2|1}_0(M)\to \L^{2|1}_0(M')$. 
\end{defn}

\begin{rmk} The objects of $\L^{2|1}(M)$ are fields for the $2|1$-dimensional $\sigma$-model with target~$M$ from Example~\ref{ex:sigma}. \end{rmk}

\begin{rmk}\label{eq:superlatices}
In choosing lattices for super tori in Definition~\ref{defn:supertori}, it is also possible to consider $S$-families of homomorphisms $S\times \Z^2\to \E^{2|1}$ into the (non-commutative) super Lie group~$\E^{2|1}$ rather than just $\E^2<\E^{2|1}$. Although mathematically well-defined, this doesn't square well with the physics: the $\Z^2$-quotient is by the right action, and so the resulting family of super tori has a residual left action of~$\E^{2|1}$. Typically one defines Lagrangians for classical field theories (e.g., see~\eqref{eq:superLag}) by the left-invariant vector fields which commute with the left action. However, for these more general lattices such left-invariant vector fields can fail to be invariant under the $S\times \Z^2$-action, and so the associated Lagrangians fail to descend to the families of super tori. On could instead define super tori in terms of a $\Z^2$-quotient by the left action, so that the Lagrangian does descend. However, generic such families of super tori only have a residual left action by the subgroup~$\E^2<\E^{2|1}$, so that symmetries for this Lagrangian field theory are a strict subgroup of the rigid conformal group. In physical jargon, \emph{supersymmetry is broken} in this compactification. In our applications to cohomology, the de~Rham operator comes from the generator of supersymmetry. As such, breaking supersymmetry is undesirable for us and so we work with the more restrictive notion of lattices and the associated vacua in which supersymmetry is preserved. 
\end{rmk}

\begin{rmk} 
In parallel to Remark~\ref{rmk:KUKO}, there are variants of the stacks $\L^{2|1}(M)$ where one restricts the isomorphisms between families of super tori from Definition~\ref{defn:21fibconf} to have isomorphisms between based lattices lying in a subgroup $\Gamma<\GL_2(\Z)$. Choosing the level~$N$ subgroups results in a cocycle model for $\TMF\otimes \C$ with level structure. 
\end{rmk}

Define a morphism of stacks $\L^{2|1}_0(\pt)\to [\pt\sq \C^\times]$ that sends all objects over~$S$ to $\pt$ and to a fiberwise rigid conformal map $T^{2|1}_\ell\to T^{2|1}_{\ell'}$ associates the dilation part $S\to \C^\times$. As in the previous cases, this description in terms of $S$-points of the prestack defining $\L^{2|1}_0(\pt)$ suffices to define a morphism of stacks; see Remark~\ref{rmk:stackad}. 

\begin{defn} Define line bundles $\omega^{\otimes k/2}$ over $\L^{2|1}_0(M)$ as the pullback of the $k$th tensor power of the canonical odd line over $[\pt\sq \C^\times]$ along the composition
$$
\L^{2|1}_0(M)\to \L^{2|1}_0(\pt)\to [\pt\sq \C^\times],
$$
where the first arrow is the induced functor from $M\to \pt$. 
\end{defn}

Next we find a Lie groupoid presentation of $\L^{2|1}_0(M)$ in which we can compute sections of $\omega^{\otimes k/2}$. Consider the functor 
\beq
&&u\colon \Lat\times \Map(\R^{0|1},M)\to \L^{2|1}_0(M), \label{eq:21atlas}
\eeq
that sends an $S$-point $\ell\in \Lat(S)$ and $\phi_0\colon S\times \R^{0|1}\to M$ to
$$
\phi\colon T^{2|1}_\ell\stackrel{pr}{\to}  S\times \R^{0|1}\stackrel{\phi_0}{\to} M
$$
where $pr$ is determined by the standard projection $\R^{2|1}\to \R^{0|1}$. 

\begin{prop} 
The map~\eqref{eq:21atlas} defines an atlas for the stack whose associated groupoid presentation is 
$$
\left( \begin{array}{c} (\E^{2|1}\rtimes \C^\times\times\SL_2(\Z)\times \Lat)/\Z^2 \times \Map(\R^{0|1},M) \\ \downarrow \downarrow \\ \Lat\times \Map(\R^{0|1},M)\end{array} \right) \stackrel{\sim}{\to} \L^{2|1}_0(M).
$$
The source map is the projection. The target map is the composition
\beq
(\E^{2|1}\rtimes \C^\times\times\SL_2(\Z)\times  \Lat)/\Z \times \Map(\R^{0|1},M))&\to& \E^{0|1}\rtimes \C^\times\times\SL_2(\Z)\times  \Lat \times \Map(\R^{0|1},M))\nonumber\\
&\to& \Lat\times \Map(\R^{0|1},M)\nonumber
\eeq
where the first map is determined by the projection $pr\colon (\Lat\times \E^{2|1})/\Z\to \Lat\times \E^{0|1}$, and the second map is determined by the left action of $\E^{0|1}\rtimes \C^\times$ on $\Map(\R^{0|1},M)$ and the action of $\C^\times \times \SL_2(\Z)$ on $\Lat$. 
\label{prop:21present}
\end{prop}

\bp The stack $\L^{2|1}_0(M)$ is defined as the stackification of a prestack denoted $\L^{2|1}_0(M)_{\rm pre}$ in this proof. We will show that $\L^{2|1}_0(M)_{\rm pre}$ is equivalent to the prestack defined by the Lie groupoid in the proposition. This implies that their stackifications are equivalent. To show that a map of prestacks is an equivalence, it suffices to demonstrate an equivalence of groupoids at each $S$-point. 

First we verify essential surjectivity. The prestack $\L^{2|1}_0(M)_{\rm pre}$ has objects over~$S$ given by~$(\ell,\phi)$ for~$\ell\in \widetilde{\Lat}(S)$ and $\phi\colon T^{2|1}_\ell\to M$. For any given family of super tori $T^{2|1}_\ell$, there is a family $T^{2|1}_{\ell'}$ for which the associated $S$-family of based lattices is oriented $\ell'\in \Lat(S)\subset \widetilde{\Lat}(S)$ and a fiberwise rigid conformal map $T^{2|1}_\ell\to T^{2|1}_{\ell'}$. To see this, first we observe that on a connected component of~$S$, the ratio of the generators of a based lattice either defines an $S$-point of the upper half plane or the lower half plane in~$\C$. In the former case, the family is already oriented over that connected component of~$S$. In the latter case, modifying the generators by an orientation-reversing element of $\GL_2(\Z)$ (e.g., the element that exchanges the order of the generators) defines a fiberwise conformal map whose target is a family of super tori associated with an oriented lattice. Finally, by definition there is a unique $\phi_0\colon S\times \R^{0|1}\to M$ such that $\phi=\phi_0\circ pr$ for $pr\colon T^{2|1}_\ell\to S\times \R^{0|1}$ the projection. Hence, the map~\eqref{eq:21atlas} induces an essential surjection on $S$-points. 

We will extend~\eqref{eq:21atlas} to a map from the appropriate Lie groupoid by characterizing morphisms in~$\L^{2|1}_0(M)_{\rm pre}$ over~$\Lat\times\Map(\R^{0|1},M)$. These morphisms are determined by families of rigid conformal maps. We characterized these immediately after Definition~\ref{defn:21fibconf}. If we restrict to $S$-families of super tori with $\ell\in \Lat(S)\subset\widetilde{\Lat}(S)$ an oriented lattice, this restricts the fiberwise rigid conformal maps to sections of~$(S\times \E^{2|1}\rtimes \C^\times\times \SL_2(\Z))/\Z^2\to S$, i.e., $\GL_2(\Z)$ gets replaced by the subgroup~$\SL_2(\Z)$.  With $S=\Lat\times \Map(\R^{0|1},M)$, morphisms over $S$ are $(\E^{2|1}\rtimes \C^\times\times \SL_2(\Z)\times  \Lat)/\Z^2 \times \Map(\R^{0|1},M)$. Hence, to an $S$-point of the Lie groupoid in the statement we obtain an $S$-family of rigid conformal maps, and hence a morphism of $\L^{2|1}_0(M)_{\rm pre}$ over~$S$. By construction, this map gives an equivalence of groupoids at each $S$-point, and hence and equivalence of prestacks. 

It remains to verify that the source and target maps are as claimed. The identity rigid conformal map leaves~$(\ell,\phi)$ unchanged, so the source map is the projection. The action of~$\SL_2(\Z)$ on~$\Lat$ is the standard one, as the associated rigid conformal maps are change of basis maps $S\times \Z^2\to S\times \Z^2$. This leaves the map~$\phi$ unchanged, so the $\SL_2(\Z)$-action on $\Map(\R^{0|1},M)$ is trivial. The action of~$\E^{2|1}\rtimes \C^\times$ on~$\ell$ is through the projection to~$\C^\times$ which acts by rotating and dilating based lattices, and on the map~$\phi$ through precomposition with a rigid conformal map of a family of super tori. Since $\phi$ factors through $pr$, the associated action of $\E^{2|1}\rtimes \C^\times$ on $\phi_0$ factors through the homomorphism $\E^{2|1}\rtimes \C^\times\to \E^{0|1}\rtimes \C^\times$ and is precisely the claimed composition in the statement. This verifies the target map.
\ep

For future reference, we denote the map from the family of super tori to~$M$ associated with the atlas~\eqref{eq:21atlas} by $\ev$, 
\beq
\ev\colon (\Lat\times \Map(\R^{0|1},M)\times \R^{2|1})/\Z^2\to \Map(\R^{0|1},M)\times \R^{0|1} \to M\label{eq:21ev}
\eeq
where the first arrow is the projection, and the second arrow is the usual evaluation map.

The projection 
$$
\E^{2|1}\rtimes \C^\times\times \SL_2(\Z) \twoheadrightarrow \C^\times
$$
is a homomorphism which gives a functor from the groupoid in Proposition~\ref{prop:21present} to the Lie groupoid $\pt\sq \C^\times$. The pullback of the canonical odd line bundle on $\pt\sq \C^\times$ gives a line bundle over the Lie groupoid whose line bundle over the underlying stack is isomorphic to~$\omega^{1/2}$. This allows us to compute sections of $\omega^{\otimes k/2}$ in terms of functions on $\Lat\times \Map(\R^{0|1},M)$ with transformation properties.

Before turning to this computation, we define what it means for a section of $\omega^{\otimes k/2}$ to be holomorphic. This definition is set up so that the field-theoretic constructions in the remainder of this section give rise to holomorphic sections. Let $\vol\in C^\infty(\Lat)$ denote the smooth function on lattices that reads of the volume of a torus associated to a given lattice. In terms of lattice generators, $\vol:=\frac{1}{2i}(\ell\bar\ell'-\bar\ell\ell')$. Consider the composition 
\beq
\Lat\times \Map(\R^{0|1},M)\to \Lat\times \Map(\R^{0|1},M)\to \L^{2|1}_0(M)\label{eq:atlasmod}
\eeq
where the first arrow on $S$-points of $\Lat\times \Map(\R^{0|1},M)\cong \Lat\times \Pi TM$ is $(\ell,x,\psi)\mapsto (\ell,x,\vol^{-1/2}(\ell)\psi)$, i.e., rescale the fibers of $\Pi TM\to M$ by $\vol^{1/2}$. We observe that the pullback of a section of~$\omega^{\otimes k/2}$ under the map above is an $\SL_2(\Z)$-invariant function on~$\Lat\times \Map(\R^{0|1},M)$, and so its restriction to $\mathbb{H}\times \Map(\R^{0|1},M)$ has a well-defined $q$-expansion.

\begin{defn} \label{defn:21dilinvt} 
Let $\bar\partial$ be the restriction of the $\bar\partial$ operator on $\C\times \C$ to $\Lat\subset \C\times \C$. A section of $\omega^{\otimes k/2}$ is \emph{holomorphic} if
\begin{enumerate}
\item its pullback along~\eqref{eq:atlasmod} is in the kernel of~$\bar\partial$, and 
\item for compact subsets $K\subset M$, the further pullback to $\mathbb{H}\times \Map(\R^{0|1},K)\subset \Lat\times \Map(\R^{0|1},M)$ defines an element of $\mathcal{O}(\mathbb{H};C^\infty(\Map(\R^{0|1},K)))\hookrightarrow C^\infty(\mathbb{H}\times \Map(\R^{0|1},K))$ that is meromorphic as~$\tau\to i\infty$ for~$\tau\in \mathbb{H}$ (compare Definition~\ref{defn:MF}).
\end{enumerate}
We denote the space of holomorphic sections by
$$
\mathcal{O}(\L^{2|1}_0(M);\omega^{\otimes k/2})\subset \Gamma(\L^{2|1}_0(M);\omega^{\otimes k/2}).
$$
\end{defn}

\begin{rmk} \label{rmk:holo}
The above holomorphic condition on sections might seem a bit arbitrary: why not study \emph{all} smooth sections of $\omega^{\otimes k/2}$? The restriction is motivated by the fact that the sections constructed by field theory techniques (namely, partition functions of supersymmetric field theories) have these additional properties. Holomorphic dependence on~$\Lat$ is a consequence of \emph{chiral} supersymmetry because the anti-holomorphic contributions to the partition function from bosons and fermions precisely cancel, e.g., see~\cite[pg.~375]{Pilch}. The condition as $\tau\to i\infty$ is equivalent to having a $q$-expansion with finitely many negative powers, which in turn is a consequence of energy in quantum field theory being bounded below. 
See~\cite[4.4]{ST11} for a more details. More practically---and for the the purposes of this paper---the pushforwards we construct \emph{do} apply to arbitrary smooth sections of~$\omega^{\otimes k/2}$, but they preserve the subspace of holomorphic sections as defined above. 
Of course, this holomorphic condition is also required to obtain cocycles for cohomology with coefficients in modular forms. 
\end{rmk}

\begin{proof}[Proof of Theorem~\ref{thm1}, $d=2$ case]
By Proposition~\ref{prop:21present} and~\eqref{compsections}, (smooth) sections of $\omega^{\otimes k/2}$ are functions on the objects, $\Lat\times \Map(\R^{0|1},M)$, that are equivariant for the groupoid action. We will further impose the holomorphic condition from Definition~\ref{defn:21dilinvt} on these smooth sections. 

To start, we have isomorphisms
$$
C^\infty(\Map(\R^{0|1},M)\times \Lat)\cong C^\infty(\Map(\R^{0|1},M);C^\infty(\Lat))\cong \Omega^\bullet(M;C^\infty(\Lat)).
$$
By virtue of the definition of~$\omega^{\otimes k/2}$, sections are invariant under the action of $\E^{2|1}$ and $\SL_2(\Z)$. We have that $\E^{2|1}$ acts purely through $C^\infty(\Map(\R^{0|1},M))\cong \Omega^\bullet(M)$ as the de~Rham operator, and $\SL_2(\Z)$ acts purely through $C^\infty(\Lat)$ by changing the basis of the lattice. So we have
$$
\Omega^\bullet(M;C^\infty(\Lat))^{\E^{2|1}\times \SL_2(\Z)}\cong \Omega^\bullet_\cl(M;C^\infty(\Lat)^{\SL_2(\Z)})\cong \bigoplus_j \Omega^j_\cl(M;C^\infty(\Lat)^{\SL_2(\Z)}).
$$
We can decompose an element $\alpha\in \Omega^\bullet_\cl(M;C^\infty(\Lat)^{\SL_2(\Z)})$ according to its differential form degree and we write this as $\alpha=\sum_j \vol^{j/2}\alpha_j$ where $\alpha_j\in \Omega^j_\cl(M;C^\infty(\Lat)^{\SL_2(\Z)})$. Holomorphy in the sense of Definition~\ref{defn:21dilinvt} then demands that $\alpha_j\in \Omega^j_\cl(M;\mathcal{O}(\Lat)^{\SL_2(\Z)})\subset \Omega^j(M;C^\infty(\Lat)^{\SL_2(\Z)})$, and that the $q$-expansion has only finitely many negative powers of~$q$ on compact subsets of~$M$. 

It remains to understand $\C^\times$-equivariance. The action on $\mathcal{O}(\Lat)^{\SL_2(\Z)}$ is the action associated to the weight of modular forms, so we may assume that $\vol^{j/2}\alpha_j\in \Omega^j(M;\MF^l)$, using the notation of the previous paragraph. The action on degree $j$-forms is by $\bar\mu^{j/2}$. So the $\C^\times$-action on $\vol^{j/2}\alpha_j$ is
$$
\vol^{j/2}\alpha_j\mapsto (\mu\bar\mu\vol^{1/2})^j\bar\mu^{-j/2}\alpha_j\mu^{l/2}=\mu^{(j+l)/2}\vol^{j/2}\alpha_j.
$$ 
To be a section of $\omega^{\otimes k/2}$, we require that $j+l=k$. This means we have 
$$
\mathcal{O}(\L^{2|1}_0(M);\omega^{\otimes k/2})\cong \bigoplus_j \Omega^j_{\rm cl}(M;\MF^{k-j}).
$$

Naturality of~$\L^{2|1}_0(M)$ in~$M$ turns $\mathcal{O}(\L^{2|1}_0(M);\omega^{\otimes \bullet/2})$ into a presheaf on manifolds. We have shown that this presheaf is isomorphic to the sheaf of closed differential forms with values in modular forms, so in fact $\mathcal{O}(\L^{2|1}_0(M);\omega^{\otimes \bullet/2})$ is a sheaf, and the isomorphism is one of sheaves. 
\ep

\label{sec:21zeromode}

Just as in the $1|1$-dimensional case, the pullback of a vector bundle along $u\colon \Lat\times \Map(\R^{0|1},M)\to \L^{2|1}_0(M)$ has some extra structure coming from automorphisms in the stack. In the $1|1$-dimensional case, we had a decomposition into weight spaces from an $S^1$-action by loop rotation. In the $2|1$-dimensional case, this is a bit more subtle as there is ``more than one torus" acting by rotation. Instead, there is a Lie groupoid action by
$$
\Rot(M):=\left\{ \begin{array}{c} (\E^2\times  \Lat)/\Z^2 \times \Map(\R^{0|1},M) \\ \downarrow \downarrow \\ \Lat\times \Map(\R^{0|1},M)\end{array} \right\}.
$$
To see this, we observe that $\Rot(M)$ is a subgroupoid of the presentation from Proposition~\ref{prop:21present}. This gives a factorization of the atlas~$u$,
$$
u\colon \Lat\times \Map(\R^{0|1},M)\to [\Rot(M)]\to \L^{2|1}_0(M). 
$$
Given a vector bundle $W$ over $\L^{2|1}_0(M)$ the above defines weight spaces for the action of $\Rot(M)$ as follows. 
Let $u^*W$ be the pullback along $u\colon \Lat\times \Map(\R^{0|1},M)\to \L^{2|1}_0(M)$. The data of a vector bundle over the stack also gives an isomorphism of vector bundles ${\sf s}^*u^*W\stackrel{\sim}{\to}{\sf t}^*u^*W$ over the supermanifold of morphisms in the Lie groupoid presentation from Proposition~\ref{prop:21present}, where~${\sf s}$ and~${\sf t}$ the source and target maps. When restricted to the morphisms of the subgroupoid~${\rm Rot}(M)$, ${\sf s}={\sf t}$ is the projection, so the isomorphism between these vector bundles is a vector bundle \emph{automorphism}. 

\begin{defn} For $W\to \L^{1|1}_0(M)$ a vector bundle and $(k,l)\in \Z\times \Z$, the \emph{$(k,l)$th weight space} is the subbundle of $u^*W$ on which the action by $\Rot(M)$ is the vector bundle automorphism gotten by multiplication by the function
$$
\exp\left(\frac{\pi}{\vol}(w(k \bar\ell_1 +l\bar{\ell_1}')-\bar{w}(k\ell_1+l\ell_2))\right)\in C^\infty((\E^2\times \Lat)/\Z^2\times \Map(\R^{0|1},M))
$$
for $(w,\bar w)$ the standard coordinates on $\E^2\cong \C$ and $(\ell_1,\bar\ell_1,\ell_2,\bar\ell_2)$ the standard coordinates on~$\Lat\subset \C\times \C$. 
\end{defn}

As in the $1|1$-dimensional case, we introduce some related terminology. The \emph{zero modes of $W$} are sections of the $(0,0)$ weight space and the nonzero modes are the span of sections in the weight spaces with $(k,l)\ne (0,0)$. 

\subsection{Vector bundles on $\L^{2|1}_0(M)$}

As in the previous case from~\S\ref{sec:KQM}, to a real oriented vector bundle $V\to M$ we define vector bundles $\Bos(V)\to \L^{2|1}_0(M)$ and $\Fer(V)\to \L^{2|1}_0(M)$. See~\S\ref{sec:01MQ} for the preliminaries on how to promote $V\to M$ with its connection and metric to a vector bundle in supermanifolds with a connection and a $C^\infty(M)$-bilinear pairing on sections. We recall that $\mathcal{V}$ denotes the sheaf of sections of~$V$ over~$M$ regarded as a supermanifold.

\begin{defn} \label{defn21ferbos}
For $V\to M$ a real oriented vector bundle, define a vector bundle $\Bos(V)\to \L^{2|1}_0(M)$ whose sections at an $S$-point $(\ell,\phi)$ is the $C^\infty(S)$-module 
$$
\Gamma(T^{2|1}_\ell,\phi^* V)=C^\infty(T^{2|1}_\ell)\otimes_{C^\infty(M)} \mathcal{V}.
$$
Define a similar vector bundle $\Fer(V)\to \L^{2|1}_0(M)$ whose sections at an $S$-point $(\ell,\phi)$ is the $C^\infty(S)$-module 
$$
\Gamma(T^{2|1}_\ell,\phi^*\Pi V)=C^\infty(T^{2|1}_\ell)\otimes_{C^\infty(M)}\Pi \mathcal{V},
$$
noting the parity reversal in this second definition. In terms of a sheaf of sections, these are the the direct image sheaves along the projection $T^{2|1}_\ell\to S$ of the sheaf of sections of the pullback $\phi^*V$ and $\phi^*\Pi V$ respectively. For an isomorphism in $\L^{2|1}_0(M)$ associated to a fiberwise rigid conformal map $f\colon T^{2|1}_\ell\to T^{2|1}_{\ell'}$ over a base change $S\to S'$, define a map of vector bundles from the pullback of sections along~$f$ followed by rescaling
$$
\Gamma(T^{2|1}_{\ell'},\phi'^* V)\stackrel{f^*}{\to} \Gamma(T^{2|1}_\ell,\phi^* V)$$$$
\Gamma(T^{2|1}_{\ell'},\phi'^*\Pi V)\stackrel{f^*}{\to} \Gamma(T^{2|1}_\ell,\phi^*\Pi V)\stackrel{\mu^{-1}}{\to} \Gamma(T^{2|1}_\ell,\phi^*\Pi V)
$$
where $\mu^{-1}\in C^\infty(S)^\ev$ is part of the data $(\mu,\bar\mu)\in \C^\times(S)$ of a rigid conformal map. 
\end{defn}

\begin{rmk} Sections of $\Bos(V)$ and $\Fer(V)$ are fields for $V$-valued free bosons and $V$-valued free fermions, respectively, whence the notation. \end{rmk}

From Proposition~\ref{prop:21present}, the vector bundles $\Bos(V)$ and $\Fer(V)$ can be characterized in terms of a vector bundle over the presenting Lie groupoid, i.e., vector bundles $u^*\Bos(V)$ and $u^*\Fer(V)$ on the objects $\Lat\times \Map(\R^{0|1},M)$ and isomorphisms of vector bundles ${\sf s}^*u^*\Bos(V)\stackrel{\sim}{\to} {\sf t}^*u^*\Bos(V)$ and ${\sf s}^*u^*\Fer(V)\stackrel{\sim}{\to} {\sf t}^*u^*\Fer(V)$ over the morphisms. Explicitly, these vector bundles over objects are the $C^\infty(\Lat\times \Map(\R^{0|1},M))$-modules
\beq
\Gamma(\Lat\times \Map(\R^{0|1},M);u^*\Bos(V))&=&\Gamma((\Lat\times \Map(\R^{0|1},M)\times \R^{2|1})/\Z^2,\ev^* V)\nonumber\\
\Gamma(\Lat\times \Map(\R^{0|1},M);u^*\Fer(V))&=&\Gamma((\Lat\times \Map(\R^{0|1},M)\times \R^{2|1})/\Z^2,\ev^*\Pi V)\nonumber
\eeq
for $\ev$ the map~\eqref{eq:21ev}. 

We define component fields for sections 
$$
X\in \Gamma(\Lat\times \Map(\R^{0|1},M);u^*\Bos(V))\qquad \psi\in \Gamma(\Lat\times \Map(\R^{0|1},M);u^*\Fer(V))
$$ 
by the same formulas~\eqref{eq:component11} as in the $1|1$-dimensional case, 
\beq
&&\begin{array}{lll} 
a&:=&i_0^*A\in \Gamma((\Lat\times \Map(\R^{0|1},M)\times \R^2)/\Z^2,i_0^*\ev^* V) \\ 
\eta&:=&i_0^* (\nabla_DA)\in \Gamma((\Lat\times \Map(\R^{0|1},M)\times \R^2)/\Z^2,i_0^*\ev^*\Pi V)\\
\psi_1&:=&i_0^*\psi\in \Gamma((\Lat\times \Map(\R^{0|1},M)\times \R^2)/\Z^2,i_0^*\ev^*\Pi V)\\ 
\psi_0&:=&i_0^* (\nabla_D\psi)\in \Gamma((\Lat\times \Map(\R^{0|1},M)\times \R^2)/\Z^2,i_0^*\ev^*V)\end{array} \label{eq:component21}
\eeq
where now $D=\partial_\theta-\theta\partial_{\bar z}$, $i_0\colon (\Lat \times \Map(\R^{0|1},M))\times \R^2)/\Z^2\hookrightarrow (\Lat\times \Map(\R^{0|1},M)\times \R^{2|1})/\Z^2$ is the fiberwise inclusion of the reduced manifold, and $\ev$ is the relevant evaluation map~\eqref{eq:21ev}. This gives an isomorphism of $C^\infty(\Lat\times \Map(\R^{0|1},M))$-modules 
\beq
&&\Gamma\left((\Lat \times \Map(\R^{0|1},M)\times \R^{2|1})/\Z^2,\ev^*\Pi V\right)\nonumber \\
&&\stackrel{\sim}{\to} \Gamma\left((\Lat\times \Map(\R^{0|1},M)\times \R^2)/\Z^2,i_0^*\ev^*\Pi V\right)\oplus \Gamma\left((\Lat\times \Map(\R^{0|1},M)\times \R^2)/\Z^2,i_0^*\ev^*V\right)\nonumber\\
&&\Gamma\left((\Lat \times \Map(\R^{0|1},M)\times \R^{2|1})/\Z^2,\ev^* V\right)\nonumber \\
&&\stackrel{\sim}{\to} \Gamma\left((\Lat\times \Map(\R^{0|1},M)\times \R^2)/\Z^2,i_0^*\ev^* V\right)\oplus \Gamma\left((\Lat\times \Map(\R^{0|1},M)\times \R^2)/\Z^2,i_0^*\ev^*\Pi V\right)\nonumber
\eeq
sending $\psi\mapsto (\psi_1,\psi_0)$ and $X\mapsto (a,\eta)$. The inverse isomorphisms are
$$
(\psi_1,\psi_0)\mapsto {\rm pr}^*\psi_1+\theta {\rm pr}^*\psi_0,\quad (a,\eta)\mapsto {\rm pr}^*a+\theta {\rm pr}^*\eta
$$$$
{\rm pr}\colon (\Lat \times \Map(\R^{0|1},M)\times \R^{2|1})/\Z^2\to (\Lat\times \Map(\R^{0|1},M)\times \R^2)/\Z^2$$ 
where ${\rm pr}$ is the obvious projection. We caution that this projection is not invariant under super translations, and these Taylor components are similarly not invariant. To simplify notation we often omit the pullback along the projection, i.e., writing~$\psi=\psi_1+\theta \psi_0$ and $X=a+\theta \eta$. For $p^*V\to V$ the pullback of a vector bundle over itself, we have the canonical section $\x$ that determines a section $\ev^*\x\in \Gamma(\Lat\times \Map(\R^{0|1},V),u^*\Fer(p^*V))$ with component fields
\beq
&&\begin{array}{lll} \x_0&:=&i_0^*\ev^*\x \in \Gamma((\Lat\times \Map(\R^{0|1},V)\times\R)/\Z,i_0^*\ev^*p^*V), \\ 
\x_1&:=&i_0^* (\nabla_D\ev^*\x)\in \Gamma((\Lat\times \Map(\R^{0|1},V)\times\R)/\Z,i_0^*\ev^*p^*\Pi V)\end{array} \label{eq:component212}
\eeq
and we similarly write $\x=\x_0+\theta\x_1$.

As before, component fields have a description in terms of differential forms since the composition
$$
(\Lat\times \Map(\R^{0|1},M)\times \R^2)/\Z^2\stackrel{i_0}{\hookrightarrow} \Lat\times \Map(\R^{0|1},M))\times_\ell \R^{2|1}\stackrel{\ev}{\to} M
$$
equals the composition of projections $(\Lat\times \Map(\R^{0|1},M)\times \R^2)/\Z^2\to \Map(\R^{0|1},M)\to M$. Hence, for $\psi\in \Gamma(\Lat\times \Map(\R^{0|1},M),u^*\Fer(V))$, 
\beq
&&\psi_0\in C^\infty((\Lat\times \R^2)/\Z^2)\otimes \Omega^\bullet(M,V),\qquad \psi_1\in C^\infty((\Lat\times \R^2)/\Z^2)\otimes \Omega^\bullet(M,\Pi V), \label{eq:comp21mod}
\eeq
and for $X\in \Gamma(\Lat\times \Map(\R^{0|1},M),u^*\Bos(V))$
\beq
&&a\in C^\infty((\Lat\times \R^2)/\Z^2)\otimes \Omega^\bullet(M,V),\qquad \eta\in C^\infty((\Lat\times \R^2)/\Z^2)\otimes \Omega^\bullet(M,\Pi V), \label{eq:comp21modA}
\eeq
where the tensor product is the projective tensor product of Fr\'echet spaces.

As a special case of Definition~\ref{defn21ferbos} above, we obtain vector bundles on $\L^{2|1}_0(M)$ from families of oriented manifolds $\pi\colon M\to B$. Namely,~$\BosT \to \L^{2|1}_0(M)$ whose space of sections at an $S$-point $(\ell,\phi)$ is the $C^\infty(S)$-module
$$
\Gamma(T^{2|1}_\ell;\phi^*T(M/B))=C^\infty(T^{2|1}_\ell)\otimes_{C^\infty(M)} \Gamma(T(M/B)).
$$
\begin{rmk}
Geometrically, when $B=\pt$ the nonzero modes of $\BosT$ are the normal bundle to the inclusion $\L^{2|1}_0(M)\hookrightarrow \L^{2|1}(M)$. For general~$B$ these sections are a relative version of such a normal bundle. 
\end{rmk}

\subsection{The analytic pushforward} \label{sec:WG}

We recall the action functional $\mathcal{S}_\sigma$ from~\eqref{eq:superLag} with quadratic part~\eqref{eq:sigmahess}. Applying this quadratic part to sections of $u^*\BosT$ we obtain the action functional
\beq
\begin{array}{c}\displaystyle
\phantom{BLAH} \HessS(X):=\int \langle X,\nabla_{\partial_z}\nabla_DX\rangle [d\theta d\bar z dz],\\ X\in \Gamma((\Lat \times \Map(\R^{0|1},M)\times \R^{2|1})/\Z^2;\ev^*T(M/B))\end{array} \label{eq:21HessS}
\eeq
where $\nabla$ is the pullback of the connection on~$T(M/B)$, $\langle-,-\rangle$ is the pullback of the metric, and the integral is with respect to the Berezinian measure $[d\theta d\bar z dz]$ on the fibers of the projection
\beq
(\Lat \times \Map(\R^{0|1},M)\times \R^{2|1})/\Z^2\to \Lat\times \Map(\R^{0|1},M).\label{eq:intfibers211}
\eeq

\begin{lem} The functional~\eqref{eq:21HessS} defined on objects descends to a functional defined on the stack: 
$$
{\sf s}^*\HessS={\sf t}^*\HessS
$$
for ${\sf s},{\sf t}$ the source and target maps in the groupoid presentation of $\L^{2|1}_0(M)$, and where the equality implicitly uses isomorphism of vector bundles~${\sf s}^*u^*\BosT\cong {\sf t}^*u^*\BosT$. 
\end{lem}

\bp The argument is very similar to the previous cases. Since it is built out of left-invariant vector fields, the integrand is automatically invariant under the left action by super translation. Together with invariance of the Berezinian under super translation, we find $\mathcal{S}_\sigma$ is super translation invariant. To study the effect of super dilations, we have formulas for the pullbacks $\nabla_z\mapsto \frac{1}{\mu^2}\nabla_z$, $\nabla_D\mapsto \frac{1}{\bar \mu} \nabla_D$ and $[d\theta d\bar z dz]\mapsto \bar \mu\mu^2[ d\theta d\bar z dz]$. Recombining these ingredients, we see that ${\sf s}^*\HessS={\sf t}^*\HessS$, proving the lemma. 
\ep

\begin{lem} \label{lem:21Hesscomp}
The value of $\HessS$ on component fields $(a,\eta)$ is 
\beq
\HessS(a,\eta)&=&\int\left(\langle(a,(-\nabla_{\partial_z}\nabla_{\partial_{\bar z}}-\nabla_{\partial_z}R)a\rangle+\langle\eta,\nabla_{\partial_z}\eta\rangle\right) d\bar z dz,\nonumber
\eeq
where $\nabla$ is the pullback of the connection on $T(M/B)$ along $\ev\circ i_0$,  $R\in \Omega^2(M,\End(T(M/B)))$ is the curvature 2-form of the fiberwise Levi-Civita connection acting on component fields using the identification~\eqref{eq:comp21modA}, and the integral is over the fibers
$$
(\Lat\times\Map(\R^{0|1},M)\times\R^2)/\Z^2\to \Lat\times \Map(\R^{0|1},M).
$$

\end{lem}

The proof is completely analogous to the $1|1$-dimensional case, Lemma~\ref{lem:11Hesscomp}. 

Our regularization procedure requires we also consider a complementary bundle $\Fer(n):=\Fer(\underline{\R}^n)$ from Definiton~\ref{defn11ferbos} for $\underline{\R}^n$ the trivial bundle. Equipping it with the trivial connection, we get a family of operators~$D$ acting on sections of $u^*\Fer(n)$ and an associated functional on sections
$$
\mathcal{S}_\Fer(\psi)=\frac{1}{2}\int \langle \psi,D\psi\rangle [d\theta d\bar z dz]
$$
$$
\psi\in \Gamma((\Lat\times \Map(\R^{0|1},M)\times \R^{2|1})/\Z^2,\Pi \underline{\R}^n),
$$
where $\langle-,-\rangle$ is the standard metric on $\R^n$ and the integral is with respect to the Berezinian measure $[d\theta d\bar z dz]$ on the fibers of the usual projection. Writing a section in terms of component fields, $\psi=\psi_1+\theta \psi_0$ and integrating with respect to~$\theta$, we find
$$
\mathcal{S}_\Fer(\psi_1,\psi_0)=\frac{1}{2}\int \left(\langle \psi_0,\psi_0\rangle+\langle \psi_1,\partial_{\bar z}\psi_1\rangle\right) d\bar z dz.
$$
This formula is a special case of Proposition~\ref{prop:SMQ21} below.

\begin{notation} Let $\Bos_{k,l}(T(M/B))$ denote the $(k,l)$th weight space of the vector bundle $u^*\BosT$ on $\Lat\times \Map(\R^{0|1},M)$. Denote the decomposition into even and odd parts of this vector bundle as~$\Bos_{k,l}(T(M/B))= \Bos_{k,l}^\ev(T(M/B))\oplus \Bos_{k,l}^\odd(T(M/B))$.  Similarly, let $\Fer_{k,l}(n)$ denote the $(k,l)$th weight space of the vector bundle $u^*\Fer(n)$ on $\Lat\times \Map(\R^{0|1},M)$. Denote the decomposition into even and odd parts of this vector bundle as~$\Fer_{k,l}(n)=\Fer_{k,l}^\ev(n)\oplus\Fer_{k,l}^\odd(n)$.
\end{notation}

\begin{proof}[Proof of Theorem \ref{thm:1}]

By Definition~\ref{defn:21sigmareg}, we need to calculate the finite-dimensional super determinants 
\beq
 \sdet(\nabla_{\partial_z}\nabla_D|_{\Bos_{k,l}(T(M/B))})&=&\left(\frac{\det(\nabla_{\partial_z}\nabla_D|_{\Bos_{k,l}^\odd(T(M/B))})}{\det(-\nabla_{\partial_z}\nabla_D|_{\Bos^\ev_{k,l}(T(M/B))})}\right)= \left(\frac{\det(\nabla_{\partial_z})}{\det(-\nabla_{\partial_{\bar z}}\nabla_{\partial_z}- \nabla_{\partial_z}R)}\right)\nonumber\\
 \sdet(D|_{\Fer_{k,l}(n)})&=&\frac{\det(\partial_{\bar z}|_{\Fer_l^\odd(n)})}{\det(\id |_{\Fer_l^\ev(n)})}\nonumber
\eeq
where the last equality in the first line is from Lemma~\ref{lem:21Hesscomp}. As before, we compute these determinants locally in~$M$ using a local basis for the $(k,l)$th weight space component fields
$$
f_{k,l}\otimes a_j\in C^\infty((\Lat\times \R^2)/\Z^2)\otimes \Omega^\bullet(M,T(M/B)),\quad f_{k,l}\otimes \eta_j\in C^\infty((\Lat\times \R^2)/\Z^2)\otimes \Omega^\bullet(M,\Pi T(M/B)),
$$
where $\{a_j\}$ (respectively, $\{\eta_j\}$) are local generators for the $C^\infty(M)$-module of sections of~$T(M/B)$ (respectively, $\Pi T(M/B)$) and 
$$
f_{k,l}(z,\bar{z}):=\exp\left(\frac{\pi}{\vol}(\bar{z}(k\ell_1+l\ell_2)-z(k \bar\ell_1 +l\bar\ell_2))\right)\in C^\infty((\Lat\times \R^2)/\Z^2). 
$$
In terms of these bases, $\nabla_{\partial_z} =\partial_z\otimes {\rm id}$ and $\nabla_{\partial_{\bar z}} =\partial_{\bar z}\otimes {\rm id}$ so we have
$$
(\nabla_{\partial_z}\nabla_D)|_{\Bos_{k,l}^\ev(T(M/B))}=-\frac{\partial^2}{\partial_z\partial_{\bar z}}\otimes \id-\frac{\partial}{\partial_z} \otimes R,\quad\quad (\nabla_{\partial_z}\nabla_D)|_{\Bos_{k,l}^\odd(T(M/B))}=\frac{\partial}{\partial_z}\otimes \id.
$$
Similarly we have a global basis 
$$
f_{k,l}\otimes \psi_j^0\in C^\infty(\R/\Z)\otimes \Omega^\bullet(M,\R^n),$$
$$ 
f_{k,l}\otimes \psi_j^1\in C^\infty(\R/\Z)\otimes \Omega^\bullet(M,\Pi \R^n),
$$ 
where $\psi_j^0$ is a basis for $\R^n$ and $\psi_j^1$ is a basis for $\Pi \R^n$. 

On the $(k,l)$th weight space, 
\beq
\det\left(\frac{\partial}{\partial z}\otimes \id\right)&=&\det\left(-\frac{\pi (k\bar\ell_1+l\bar\ell_2)}{\vol}\otimes \id\right)=\left(-\frac{\pi (k\bar\ell_1+l\bar\ell_2)}{\vol}\right)^{n}\nonumber 
\eeq
and 
\beq
\det\left(-\frac{\partial^2}{\partial_z\partial_{\bar z}}\otimes \id-i\frac{\partial}{\partial_z} \otimes R\right)&=&\det\left(\frac{\pi^2}{\vol^2}|k\ell_1+l\ell_2|^2\otimes \id+\frac{\pi}{\vol}(k\bar\ell_1+l\bar\ell_2)\otimes R\right)\nonumber\\
&=&\det\left(\frac{\pi^2}{\vol^2}|k\ell_1+l\ell_2|^2\otimes \id\right)\det\left(\id\otimes \id+\frac{\vol}{\pi(k\ell_1+l\ell_2)}\otimes R\right)\nonumber \\
&=&\left(\frac{\pi^2}{\vol^2}|k\ell_1+l\ell_2|^2\right)^{n}\det\left(\id\otimes \id+\frac{\vol}{\pi(k\ell_1+l\ell_2)}\otimes R\right)\nonumber
\eeq
where $n={\rm dim}(M)-{\rm dim}(B)$. On $\Fer_{k,l}(n)$ we have 
$$
\frac{\det(\partial_{\bar z}|_{\Fer_{k,l}^\odd(n)})}{\det(\id |_{\Fer_{k,l}^\ev(n)})}=\left(-\frac{\pi (k\ell_1+l\ell_2)}{\vol}\right)^{n}.
$$
Putting this together, the contribution to the regularized determinant from the $(k,l)$th weight space is 
\beq
\det\left(\id\otimes \id+\frac{\vol}{\pi(k\ell_1+l\ell_2)}\otimes R\right)&=&\exp({\rm Tr}(\log(\id\otimes \id+\frac{\vol}{\pi(k\ell_1+l\ell_2)}\otimes R)))\nonumber\\
&=&\exp\left({\rm Tr}\left(\sum_m \frac{(-1)^{m+1}}{m}\frac{\vol^m}{(\pi(k\ell_1+l\ell_2))^m}R^m\right)\right)\nonumber \\
&=&\exp\left({\rm Tr}\left(-\sum_m \frac{1}{2m}\frac{\vol^{2m}}{(\pi(k\ell_1+l\ell_2))^{2m}}R^{2m}\right)\right)\nonumber
\eeq
where in the last equality we used that traces of odd powers of the curvature vanish. 

Next we take the product over weight spaces. For $m>1$, the contribution to the regularized determinant is the exponential of the absolutely convergent sum involving the Eisenstein series~$E_{2m}$
\beq
\sum_{(k,l)\in \Z^2_+} {\rm Tr}\left(\sum_m \frac{1}{2m}\frac{\vol^{2m}}{(\pi (k\ell_1+l\ell_2))^{2m}}R^{2m}\right)&=&\sum_{(k,l)\in \Z^2_*} {\rm Tr}\left(\sum_m \frac{1}{4m}\frac{\vol^{2m}}{(\pi (k\ell_1+l\ell_2))^{2m}}R^{2m}\right)\nonumber\\
&=&\sum_{m=1}^\infty\frac{\vol^{2m} E_{2m}}{4m \pi^{2m} }{\rm Tr}(R^{2m}),\nonumber
\eeq
where we recall that
$$
\Z^2_+=\{(k,l)\in \Z^2\mid k>0 \ {\rm or} \ k=0, l>0\},\quad \Z^2_*=\{(k,l)\in \Z^2\mid (k,l)\ne (0,0)\}.
$$
For $m=1$ we get a contribution from the exponential of the conditionally convergent sum,
$$
\sum_{(k,l)\in \Z^2_+} {\rm Tr}\left(\frac{\vol^{2}}{2(\pi (k\ell_1+l\ell_2))^{2}}R^{2}\right)=\frac{\vol^{2} E_{2}^{\rm ren}}{4 \pi^{2} }{\rm Tr}(R^{2})
$$
where $E_2^{\rm ren}$ is a version of the 2nd Eisenstein series that depends on a choice of ordering the sum. There is no modular form of weight~2, so there is no choice so that~$E_2^{\rm ren}$ is a modular form. 

For a real vector bundle $V$, let ${\rm ph}_{k}(V)$ denote the degree $4k$ part of the Pontryagin character of $V$, i.e., ${\rm ph}_k(V)={\rm ch}_{2k}(V\otimes \C)$ where ${\rm ch}_k$ is the degree~$2k$ piece of the Chern character. In our cocycle model, we have the description 
$$
(2\vol)^{2k}{\rm Tr}(R^{2k})=(-1)^k2(2k)! {\rm ph}_{k}(T(M/B))
$$
and so we obtain the function on $\Lat\times \Map(\R^{0|1},M)$
$$
\sdet_{\rm reg}'\left(\nabla_{\partial_z}\nabla_D\right)^{-1/2}=\exp\left(\frac{p_1(T(M/B))}{(2\pi i)^2}E_2^{\rm ren}+\sum_{k\ge 2}^\infty \frac{(2k)!{\rm ph}_k(T(M/B))}{2k(2\pi i)^{2k}} E_{2k}\right)
$$
using that ${\rm ph}_1(T(M/B))=p_1(T(M/B))$, i.e., the degree~4 part of the Pontryagin character is the 1st Pontryagin form. Since $E_2^{\rm ren}$ is not a modular form, 
$$
\sdet_{\rm reg}'\left(\nabla_{\partial_z}\nabla_D\right)^{-1/2}\in C^\infty(\Lat\times \Map(\R^{0|1},M))
$$
descends to a holomorphic function on the stack~$\L^{2|1}_0(M)$ to the stack if and only if $p_1(T(M/B))=0$ as a differential form, proving part (i) of the theorem. 

We observe that for any choice of ordering, $\sdet_{\rm reg}'\left(\nabla_{\partial_z}\nabla_D\right)^{-1/2}$ is invariant under the action of super translations, so the failure of descent is encoded by the $\C^\times\times \SL_2(\Z)$-action. The ratio of $\sdet_{\rm reg}'\left(\nabla_{\partial_z}\nabla_D\right)^{-1/2}$ when pulled back along the projection and action maps,
$$
{\rm pr},{\rm act}\colon \C^\times\times \SL_2(\Z)\times \Lat\times \Map(\R^{0|1},M)
$$
is the nonvanishing function ${\rm pr}^*\sdet_{\rm reg}'\left(\nabla_{\partial_z}\nabla_D\right)^{-1/2}/{\rm act}^*\sdet_{\rm reg}'\left(\nabla_{\partial_z}\nabla_D\right)^{-1/2}$ that defines a cocycle for a (super translation invariant) line bundle on~$\L^{2|1}(M)$ we denote by $\mathcal{A}(p_1)$. This proves part (ii) of the theorem. 

Given a rational string structure $H$ with $dH=p_1(T(M/B))$, we can consider the concordance of sections
\beq
&&\exp\left(\frac{d(\lambda H)}{(2\pi i)^2}E_2^{\rm ren}+\sum_{k\ge 2}^\infty \frac{(2k)!{\rm ph}_k(T(M/B))}{2k(2\pi i)^{2k}} E_{2k}\right)
\in \Gamma(\L^{2|1}_0(M\times \R);\mathcal{A}(d(\lambda H))).\label{eq:ellThomconc}
\eeq
In the target of this concordance, the term involving the $2^{\rm nd}$ Eisenstein series is eliminated so the line bundle is trivialized. Said differently, the result is a function on the atlas $\Lat\times \Map(\R^{0|1},M)$ that descends to a holomorphic function on~$\L^{2|1}_0(M)$. By construction, this function is a cocycle representative of the (modular) Witten class of the family $\pi\colon M\to B$. By inspection, the resulting quantity is independent of the choice of ordering defining $E^{\rm ren}_2$ and the rational string structure. This proves parts (iii) and (iv) of the theorem. 
\ep

\subsection{The elliptic Mathai--Quillen form}\label{sec:ellMQ}

We specialize the free fermion vector bundle to the one associated with $p^*V\to V$ gotten by pulling $V$ back along itself. There is a function on sections of $u^*\Fer(p^*V)$ given by the usual formula
\beq
\SMQ(\psi)=\int \left(\frac{1}{2}\langle \psi,\nabla_D\psi\rangle+\frac{i}{\sqrt{\vol}}\langle \psi,\phi^*\x\rangle\right)[d\theta d\bar z dz] \label{eq:21MQaction}
\eeq
where $D=\partial_\theta -\theta\partial_{\bar z}$, and the integral is over the fibers of the bundle of super tori 
\beq
(\Lat \times \Map(\R^{0|1},V)\times \R^{2|1})/\Z^2\to \Lat\times \Map(\R^{0|1},V)\label{eq:21integral}
\eeq
using the fiberwise Berezinian measure $[d\theta d\bar z dz]$ that descends from $\R^{2|1}$. 

\begin{lem} The functional~\eqref{eq:21MQaction} defined on objects descends to a functional defined on the stack: 
$$
{\sf s}^*\SMQ={\sf t}^*\SMQ
$$
for ${\sf s},{\sf t}$ the source and target maps in the groupoid presentation of $\L^{2|1}_0(V)$, and where the equality implicitly uses the isomorphism of vector bundles~${\sf s}^*u^*\Fer(V)\cong {\sf t}^*u^*\Fer(V)$. 
\end{lem}

\bp
The argument proceeds analogously to the previous versions: invariance under super translations follows from the fact that the integrand is build from left-invariant vector fields (the left action of super translations is generated by \emph{right}-invariant vector fields) and that the Berezinian measure is invariant under super translations. Invariance under dilations comes from pulling back the constituent pieces defining $\SMQ$:~$\psi\mapsto \mu^{-1} \psi$, $[d\theta  d\bar z d z]\mapsto \frac{1}{\bar \mu} [d\theta] \mu^2 [dz] \bar \mu^2[ d\bar z] =\bar \mu \mu^2 [d\theta d\bar z dz]$, $\nabla_D\mapsto \bar \mu^{-1}\nabla_D$, and $\sqrt{\vol}\mapsto \mu\bar\mu\sqrt{\vol}$. Assembling ${\sf t}^*\SMQ$ from these pieces, we find it agrees with~${\sf s}^*\SMQ$, and so~$\SMQ$ descends to the stack. 
\ep

The value of the Mathai--Quillen classical action on component fields is an element of~$C^\infty(\Lat)\otimes \Omega^\bullet(V)$, which we calculate in terms of more familiar quantities.

\begin{prop} \label{prop:SMQ21}
The value of $\SMQ$ on component fields $(\psi_1,\psi_0)$ is 
$$
\SMQ(\psi_1,\psi_0)=\int \left( \langle \psi_0,\frac{1}{2}\psi_0+\frac{i}{\sqrt{\vol}}\x \rangle +\langle \psi_1,\frac{1}{2}(\nabla_{\partial_{\bar z}}+p^*F)\psi_1+\frac{i}{\sqrt{\vol}}\nabla \x \rangle \right)d\bar z dz
$$
where $\nabla$ is the pullback of the connection on $p^*V$ along $\ev\circ i_0$, $F\in \Omega^2(M;\End(V))$ is the curvature 2-form of~$V$, $\x\in \Omega^0(V,p^*V)$ is the canonical section, $\nabla\x\in \Omega^1(V,p^*V)$ is its covariant derivative, and the integral is over the fibers,
\beq
(\Lat \times \Map(\R^{0|1},V)\times \R^2)/\Z^2\to \Lat\times \Map(\R^{0|1},V),\label{eq:21integral2}
\eeq
using the fiberwise measure~$d\bar z dz$.
\end{prop}

The proof is virtually identical to the proof of Proposition~\ref{prop:SMQ11}.

The next task is to give rigorous meaning to the integral $\int e^{-\SMQ(\psi)}[d\psi]_0$ as per~\eqref{MQsplit}. We start with the contribution from the zero modes. On component fields, these are the subspaces 
$$
\psi_0\in C^\infty(\Lat)\otimes \Omega^\bullet(V,p^*V)\hookrightarrow C^\infty((\Lat\times \R^2)/\Z^2)\otimes \Omega^\bullet(V,p^*V),$$$$\psi_1\in C^\infty(\Lat)\otimes \Omega^\bullet(V,p^*\Pi V)\hookrightarrow C^\infty((\Lat\times \R^2)/\Z^2)\otimes \Omega^\bullet(V,p^*\Pi V)
$$
of sections that are constant on the fibers of $(\Lat\times \R^2)/\Z^2\to \Lat$. This allows us to define the integral over zero modes as the $C^\infty(\Lat)$-linear extension of the one in~\S\ref{sec:01MQmeasure} for the $0|1$-dimensional case: the integral of the functional over the sections~$\psi_0$ is with respect to the fiberwise (real) volume form on~$p^*V\to V$ from the metric and orientation on~$V$, and the integral of the functional over the sections~$\psi_1$ is the Berezinian integral 
$$
C^\infty(\Lat)\otimes \Omega^\bullet(V,\Lambda^\bullet V^\vee)\to C^\infty(\Lat)\otimes \Omega^\bullet(V,\Lambda^{\rm top} V^\vee)\cong  C^\infty(\Lat)\otimes \Omega^\bullet(V),
$$
which uses the orientation on~$V$.

\begin{proof}[Proof of Lemma~\ref{lem:d1zeromodes} when $d=2$]
Evaluating the Mathai--Quillen classical action on a zero mode and applying Proposition~\ref{prop:SMQ21}, we obtain
\beq
\SMQ(\psi_1,\psi_0)&=&\int \left(\langle \psi_0,\frac{1}{2}\psi_0+\frac{i}{\sqrt{\vol}} \x\rangle +\langle \psi_1,\frac{1}{2}(-i\nabla_{\partial_{\bar z}}+p^*F)\psi_1+\frac{i}{\sqrt{\vol}} \nabla\x\rangle\right)d\bar z dz\nonumber\\
&=&\vol \left(\langle \psi_0,\frac{1}{2}\psi_0+\frac{i}{\sqrt{\vol}} \x\rangle+\langle \psi_1,\frac{1}{2} p^*F\psi_1+\frac{i}{\sqrt{\vol}} \nabla\x\rangle\right)\in C^\infty(S)\nonumber
\eeq
where we used $\nabla_{\partial_{\bar z}}\psi_1=0$ on zero modes and the integral is over the fibers~\eqref{eq:21integral2} so the integral evaluates simply as the fiberwise volume (namely,~$\vol$) multiplied by the integrand. So we are left to evaluate 
$$
\int e^{-\SMQ(\psi)}[d\psi]_0=\int \exp(-\vol\langle \psi_0,\frac{1}{2}\psi_0+\frac{i}{\sqrt{\vol}} \x\rangle)d \psi_0\int \exp(-\vol \langle \psi_1,\frac{1}{2}p^*F\psi_1+\frac{i}{\sqrt{\vol}} \nabla\x\rangle)d\psi_1.
$$
The first factor is a standard Gaussian integrated over the fibers of $p^*V\to V$, which we integrate with respect to the (real) volume form, obtaining
$$
\int \exp(-\vol\langle \psi_0,\frac{1}{2}\psi_0+\frac{i}{\sqrt{\vol}} \x\rangle)d \psi_0=\left(\frac{2\pi}{\vol}\right)^{{\rm dim}(V)/2} \exp(-\langle \x,\x\rangle/2)\in C^\infty(\Lat)\otimes \Omega^0_c(V).
$$
Writing the integral over $\psi_1$ as
$$
\int \exp(-\vol \langle \psi_1,\frac{1}{2}p^*F\psi_1+\frac{i}{\sqrt{\vol}} \nabla\x\rangle)d\psi_1=\int \exp(-\langle \psi_1,\frac{\vol}{2}p^* F\psi_1+i \sqrt{\vol}\nabla\x\rangle)d\psi_1,
$$
we identify the Berezinian integral
$$
\left(\frac{2\pi}{\vol}\right)^{{\rm dim}(V)/2} \exp(-\langle \x,\x\rangle/2)\int \exp(-\langle \psi_1,\frac{\vol}{2} p^*F\psi_1+i \sqrt{\vol}\nabla\x\rangle)d\psi_1
$$
with the Mathai--Quillen form~\eqref{eq:MQThom} in de~Rham cohomology under the inclusion of closed $n$-forms into $\Gamma(\L^{2|1}_0(V);\omega^{\otimes{\rm dim}(V)/2})$ in our cocycle model, up to the claimed factor of $2\pi $ and~$i$ in the statement of the lemma. 
\end{proof}

Next we unpack $\SMQ$ on the nonzero modes. 

\begin{lem}\label{lem:21perp} 

 For a nonzero mode $\psi$ of $u^*\Fer(p^*V)$, we have 
\beq
\SMQ(\psi)=\int\frac{1}{2} \langle \psi,\nabla_D\psi\rangle,\label{eq:21perp1}
\eeq
where the integral is over the fibers~\eqref{eq:21integral} and the operator $\nabla_D$ is invertible. Furthermore, $\nabla_D$ pulls back along the map~$\Lat\times \Map(\R^{0|1},V)\to \Lat\times \Map(\R^{0|1},M)$ induced by the projection $p\colon V\to M$ from a family of operators acting on the nonzero modes of $u^*\Fer(V)$. 
\end{lem}

The argument is totally analogous to the proof of Lemma~\ref{lem:11perp}.

Our regularization procedure requires we also consider a complementary bundle $\Bos(n):=\Bos(\underline{\R}^n)$ from Definiton~\ref{defn21ferbos} for $\underline{\R}^n$ the trivial bundle. Equipping it with the trivial connection, we get a family of operators~$\Delta_n=D \partial_z$ acting on sections of $u^*\Bos(n)$ and an associated functional on sections
$$
\mathcal{S}_\Bos(X)=\frac{1}{2}\int \langle X,\Delta_n X\rangle [d\theta d\bar z dz]
$$
$$
X\in \Gamma(\Map(\R^{0|1},M)\times \R^{1|1}/\Z,\underline{\R}^n),
$$
where $\langle-,-\rangle$ is the standard metric on $\R^n$ and the integral is with respect to the Berezinian measure $[d\theta d\bar z dz]$ on the fibers of the usual projection. Writing a section in terms of component fields, $X=a+\theta \eta$ and integrating with respect to~$\theta$, as a special case of Lemma~\ref{lem:21Hesscomp} we find
$$
\mathcal{S}_\Bos(a,\eta)=\frac{1}{2}\int \left(-\langle a,\partial_{\bar z}\partial_z a\rangle-\langle \eta,\partial_{\bar z}\eta\rangle\right) d\bar z dz.
$$

\begin{notation} Let $\Fer_{k,l}(V)$ denote the $(k,l)$th weight space of the vector bundle $u^*\Fer(V)$ on $\Lat\times \Map(\R^{0|1},M)$. Denote the decomposition into even and odd parts of this vector bundle as~$\Fer_{k,l}(V)\cong \Fer^\ev_{k,l}(V)\oplus \Fer^\odd_{k,l}(V)$.  Similarly, let $\Bos_{k,l}(m)$ denote the $(k,l)$th weight space of the vector bundle $u^*\Bos(m)$ on $\Lat\times \Map(\R^{0|1},M)$. Denote the decomposition into even and odd parts of this vector bundle as~$\Bos_{k,l}(m)=\Bos_{k,l}^\ev(m)\oplus\Bos_{k,l}^\odd(m)$.
\end{notation}

\begin{notation} Let $\Fer_{k,l}(V)$ denote the $(k,l)$th weight space of the vector bundle $u^*\Fer(V)$ on $\Lat\times \Map(\R^{0|1},M)$. Denote the decomposition into even and odd parts of this vector bundle as~$\Fer_{k,l}(V)= \Fer^\ev_{k,l}(V)\oplus \Fer^\odd_{k,l}(V)$.  
\end{notation}

\begin{proof}[Proof of Theorem~\ref{thm:2}] By Lemma~\ref{lem:d1zeromodes} and~\eqref{MQsplit}, it remains to compute the regularized super determinant of $\nabla_D$ and compare with the Riemann--Roch factor $\Wit(V)^{-1}$. By Definition~\ref{defn:21MQreg}, we start by computing the finite-dimensional super determinants
$$
\sdet(\nabla_D|_{\Fer_{k,l}(p^*V)})=p^*\left(\sdet(\nabla_D|_{\Fer_{k,l}(V)})\right),\qquad \sdet(\Delta|_{\Bos_{k,l}(m)})
$$
for $k,l\in \Z$ not both zero. Expanding, 
\beq
\sdet(\nabla_D|_{\Fer_{k,l}(V)})&=& \frac{\det(\nabla_D|_{\Fer_{k,l}^\odd(V)})}{\det(\nabla_D|_{\Fer_{k,l}^\ev(V)})}=\frac{\det(\nabla_{\partial_{\bar z}}+F)}{\det(\id)}\nonumber\\
\sdet(\Delta|_{\Bos_{k,l}(m)})&=&\frac{\det(\partial_{\bar z}|_{\Bos_{k,l}(m)})}{\det(\partial_{\bar z}\partial_z|_{\Bos_{k,l}(m)})}\nonumber
\eeq
where the last equality in the first line is from Lemma~\ref{lem:21perp}.

We begin by expressing the finite-dimensional determinants locally on~$M$ in a local basis for weight~$(k,l)$ component fields
$$
f_{k,l}\otimes v_j\in C^\infty((\Lat\times \R^2)/\Z^2)\otimes \Omega^\bullet(M,\Pi V),
$$
$$
f_{k,l}(z,\bar{z}):=\exp\left(\frac{\pi}{\vol}(\bar{z}(k\ell_1+l\ell_2)-z(k\bar\ell_1 +l\bar\ell_2))\right)
$$
where $\{v_j\}$ form a local basis of sections of~$\Pi V$ over an open subset~$U\subset M$ and
$$
f_{k,l}\otimes a_j\in C^\infty(\Lat)\otimes \Omega^\bullet(M,\R^m)\qquad f_{k,l}\otimes \eta_j\in C^\infty(\Lat)\otimes \Omega^\bullet(M,\Pi \R^n)
$$
for $\{a_j\}$ a basis of $\R^m$ and $\{\eta_j\}$ a basis of $\Pi \R^m$.  In terms of the basis $\{v_j\}$,~$\nabla_D|_{\Fer_l^\odd(V)}=\partial_{\bar z}\otimes \id+\id\otimes F$ and so 
\beq
\det(\nabla_{\partial_{\bar z}}+F)&=&\det(\frac{\pi}{\vol}(k\ell_1+l\ell_2)\otimes \id+\id\otimes F)\nonumber \\
&=&\det(\frac{\pi}{\vol}(k\ell_1+l\ell_2)\otimes \id)\det(\id\otimes \id+\frac{\vol}{\pi (k\ell_1+l\ell_2)}\otimes F)\nonumber \\
&=&\left(\frac{\pi}{\vol}(k\ell_1+l\ell_2)\right)^{{\rm dim}(V)}\det(\id\otimes \id+\frac{\vol}{\pi (k\ell_1+l\ell_2)}\otimes F).\nonumber
\eeq
The determinant calculation above is clearly independent of the choice of local basis $\{v_j\}$, and gives a formula that is globally defined on~$M$. We also have
$$
\sdet(\Delta|_{\Bos_{k,l}(m)})=\frac{\left(\frac{\pi}{\vol}(k\bar \ell_1+l\bar \ell_2)\right)^{{\rm dim}(V)}}{\left|\frac{\pi}{\vol}(k\ell_1+l\ell_2)\right|^{2{{\rm dim}(V)}}}=\frac{1}{\left(\frac{\pi}{\vol}(k\ell_1+l\ell_2)\right)^{{\rm dim}(V)}}.
$$
So the contribution to the regularized determinant from the $(k,l)$th weight space is 
\beq
\det(\id\otimes \id+\frac{\vol}{\pi (k\ell_1+l\ell_2)}\otimes F)&=&\exp({\rm Tr}(\log(\id\otimes \id+\frac{\vol}{\pi (k\ell_1+l\ell_2)}\otimes F)))\nonumber \\
&=&\exp\left({\rm Tr}\left(\sum_j \frac{(-1)^{j+1}}{j}\frac{\vol^j}{(\pi (k\ell_1+l\ell_2))^j}F^j\right)\right)\nonumber \\
&=&\exp\left({\rm Tr}\left(-\sum_j \frac{1}{2j}\frac{\vol^{2j}}{(\pi (k\ell_1+l\ell_2))^{2j}}F^{2j}\right)\right)\nonumber
\eeq
where in the last equality we used that traces of odd powers of the curvature vanish. 

For $j>1$, the contribution to the regularized determinant is the exponential of the absolutely convergent sum involving the Eisenstein series~$E_{2j}$
$$
\sum_{(k,l)\in \Z^2_+} {\rm Tr}\left(-\sum_j \frac{1}{2j}\frac{\vol^{2j}}{(\pi (k\ell_1+l\ell_2))^{2j}}F^{2j}\right)=-\sum_{j=1}^\infty\frac{\vol^{2j} E_{2j}}{4j \pi^{2j} }{\rm Tr}(F^{2j}),
$$
and for $j=1$ we get a contribution that is the exponential of the promised conditionally convergent sum,
$$
\sum_{(k,l)\in \Z^2_+} {\rm Tr}\left(-\frac{\vol^{2}}{2(\pi (k\ell_1+l\ell_2))^{2}}F^{2}\right)= {\rm Tr}\left(\frac{\vol^2}{2\pi^2}F^2\right) \sum_{(k,l)\in \Z^2_+} \frac{1}{((k\ell_1+l\ell_2))^{2}}=\frac{\vol^{2} E_{2}^{\rm ren}}{4 \pi^{2} }{\rm Tr}(F^{2})
$$
where $E_2^{\rm ren}$ is some choice of ordering defining a version of the 2nd Eisenstein series. There is no modular form of weight~2, so there is no choice so that~$E_2^{\rm ren}$ is a modular form.

In our cocycle model, we have
$$
(2\vol)^{2j}{\rm Tr}(F^{2j})=(-1)^j2(2j)! {\rm ph}_{j}(V)
$$
and so we obtain the function on $\Lat\times \Map(\R^{0|1},M)$
$$
\sdet_{\rm reg}'(\nabla_D)^{-1/2}=\left(\exp\left(-\frac{p_1(V)}{(2\pi i)^2}E_2^{\rm ren}-\sum_{j\ge 2}^\infty \frac{(2j)!{\rm ph}_j(V)}{2j(2\pi i)^{2j}} E_{2j}\right)\right), \nonumber
$$
using that ${\rm ph}_1(V)=p_1(V)$ and where we choose $Z^{{\rm dim}(V)}$ as the normalizing factor in the definition of $\sdet_{\rm reg}'(\nabla_D)^{-1/2}$. Since $E_2^{\rm ren}$ is not a modular form, 
$$
\sdet_{\rm reg}'(\nabla_D)^{-1/2}\in C^\infty(\Lat\times \Map(\R^{0|1},M))
$$
descends to a holomorphic function on the stack~$\L^{2|1}_0(M)$ to the stack if and only if $p_1(T(M/B))=0$ as a differential form, proving part (i) of the theorem. 

We observe that $\sdet_{\rm reg}'(\nabla_D)^{-1/2}$ is invariant under the action of super translations, so the failure of descent is encoded by the $\C^\times\times \SL_2(\Z)$-action on $\sdet_{\rm reg}'(\nabla_D)^{-1/2}$. The ratio of $\sdet_{\rm reg}'(\nabla_D)^{-1/2}$ when pulled back along the projection and action maps,
$$
{\rm pr},{\rm act}\colon \C^\times\times \SL_2(\Z)\times \Lat\times \Map(\R^{0|1},M)
$$
is the nonvanishing function ${\rm pr}^*\sdet_{\rm reg}'(\nabla_D)^{-1/2}/{\rm act}^*\sdet_{\rm reg}'(\nabla_D)^{-1/2}$ that defines a cocycle for a (super translation invariant) line bundle on~$\L^{2|1}_0(M)$ we denote by $\mathcal{A}(p_1)$. This proves part (ii) of the theorem. 

Given a rational string structure $H$ with $dH=p_1(V)$, we can consider the concordance of sections
\beq
&&\exp\left(-\frac{d(\lambda H)}{(2\pi i)^2}E_2^{\rm ren}-\sum_{k\ge 2}^\infty \frac{(2k)!{\rm ph}_k(V)}{2k(2\pi i)^{2k}} E_{2k}\right)
\in \Gamma(\L^{2|1}_0(M\times \R);\mathcal{A}(-d(\lambda H)))\label{eq:ellThomconc2}
\eeq
where $\lambda$~is the standard coordinate on~$\R$. In the target of this concordance (where $\lambda=0$), the term involving the $2^{\rm nd}$ Eisenstein series is eliminated, thereby giving a function on $\Lat\times \Map(\R^{0|1},M)$ that does descend to $\L^{2|1}_0(M)$. By construction, the section of the trivial bundle is $\Wit(V)^{-1}\in \mathcal{O}(\L^{2|1}_0(M))$. Pulling this back to $\L^{2|1}_0(V)$ and taking the product with the contribution from the zero modes (i.e., the Thom form in de~Rham cohomology) we obtain a representative of the Thom class in $\TMF^{{\rm dim}(V)}_{\rm c}(V)\otimes \C$; this proves parts (iii) and (iv) of the theorem. 
\ep

\subsection{An equality of differential pushforwards}

\begin{proof}[Proof of Theorem~\ref{thm:3}, $d=2$ case]
As discussed in Remark~\ref{rmk:complexifiedindex21}, the equality of analytic and topological pushforwards boils down to an equality of functions
\beq
&&\mathcal{O}(\L^{2|1}_0(M))\ni \Wit(\nu)^{-1}=\Wit(M/B)\in \mathcal{O}(\L^{2|1}_0(M)),\nonumber
\eeq
where the left-hand side is the Riemann--Roch factor that modifies the Thom cocycle in de~Rham cohomology, and the right-hand side is the Riemann--Roch factor that modifies integration of differential forms. The equalities of the components of the Pontryagin character~${\rm ph}_k(\nu)=-{\rm ph}_k(T(M/B))$ verify the equality of these functions. This implies that the pushforwards are equal. 
\ep

\appendix

\section{Background miscellany} \label{appen}

\subsection{Mathai--Quillen forms} \label{appen:MQ}\label{appen:geo}
The pullback of a real vector bundle $p\colon V\to M$ over itself
\beq
\begin{tikzpicture}[baseline=(basepoint)];
\node (A) at (0,0) {$p^* V$};
\node (B) at (3,0) {$V$};
\node (C) at (0,-1.5) {$V$};
\node (D) at (3,-1.5) {$M$}; 
\draw[->] (A) to (B);
\draw[->, bend left] (C) to node [left=1pt] {$\x$} (A);
\draw[->] (A) to (C);
\draw[->] (C) to node [above=1pt] {$p$} (D);
\draw[->] (B) to (D);
\path (0,-.75) coordinate (basepoint);
\end{tikzpicture}\nonumber
\eeq
has a tautological section~$\x$. If we equip $V$ with an orientation, metric $\langle-,-\rangle$, and compatible connection~$\nabla$ with curvature $F$, the \emph{Mathai--Quillen Thom form} is the Berezinian integral (c.f.~\cite{BGV}~\S1.6\footnote{Our formula differs by a factor of $1/2$ in the curvature term because contracting with the metric differs from the convention~\cite[pg.~55]{BGV} that converts the curvature into a quadratic function on sections valued in 2-forms.})
\beq
&&u_V=\frac{\epsilon({\rm dim}(V))}{(2\pi)^{{\rm dim}(V)/2}}\int^{\rm Ber} \exp\left(-\frac{1}{2}\langle \x,\x\rangle-i\langle \nabla \x,-\rangle-\frac{1}{2}\langle -,p^*F-\rangle\right)\in \Omega_{\rm cl,c}^{{\rm dim}(V)}(V)\label{eq:MQThom}
\eeq
where here and throughout the subscript ${\rm c}$ denotes \emph{compact support} (or really, rapidly decreasing support; see Remark~\ref{rmk:compactsupp}) for sections of a sheaf, and $\Omega_{\rm cl}^k$ denotes the sheaf of \emph{closed} differential $k$-forms. We set $\epsilon(n)=1$ when $n$ is even and $i$ when $n$ is odd. The Berezinian integral is the map $\Omega^\bullet(V,\Lambda V^\vee)\to \Omega^\bullet(V)$ that projects to the top power of~$\Lambda V^\vee$ and uses the orientation on~$V$ to trivialize this top power. We use square brackets, e.g., $[u_V]$, to denote the cohomology class associated with the cocycle representative~$u_V$. 

The complexification of the Thom classes in $\KO$ and $\TMF$ are modifications of usual Thom class by the inverse $\hat{A}$-class and the Witten class, respectively; see~\eqref{diag:RRThom} and~\eqref{diag:WitRR1}. These characteristic classes are determined by the power series
$$
\hat{A}(z)=\frac{z/2}{\sinh(z/2)}=\exp\left(\sum_{k=1}^\infty\frac{z^{2k}}{2k(2\pi i)^{2k}}2\zeta(2k)\right)\in \R\llbracket z\rrbracket,
$$
and 
$$
\Wit(q,z)=\frac{z/2}{\sinh(z/2)}\prod_{n\ge 1} \frac{(1-q^n)^2}{(1-q^n e^{z/2})(1-q^ne^{-z/2})}=\exp\left(\sum_{k\ge 1}\frac{E_{2k}(q)}{2k(2\pi i)^{2k}} z^{2k}\right)\in \C\llbracket z,q\rrbracket
$$
where $\zeta$ is the Riemann zeta function and $E_{2k}$ is the $2k$th Eisenstein series. 
One can build differential form representatives of~$[\hat{A}(V)^{-1}]$ and $[\Wit(V)^{-1}]$ using the Pontryagin forms associated with a chosen connection on~$V$. 

\begin{rmk} Mathai and Quillen~\cite{MathaiQuillen} constructed a differential form representative of the complexified $\KO$-theory Thom class, $[u_V]\smallsmile[\hat{A}(V)^{-1}]$, by a careful study of the Chern character applied to the Clifford modules that define Thom classes in $\KO$. 
\end{rmk}

Let $i\colon M\hookrightarrow \R^N$ be an embedding with normal bundle~$\nu$. Then the characteristic classes $[\hat{A}(\nu)^{-1}]$ and $[\Wit(\nu)^{-1}]$ coincide with $[\hat{A}(TM)]$ and $[\Wit(TM)]$, respectively. In this case, a choice of metric on~$M$ gives differential form representatives of~$[\hat{A}(TM)]$ and $[\Wit(TM)]$ using the Levi-Civita connection. When studying such differential form representatives in families, we adopt the following from~\cite[\S~\!\!a]{BismutFreed1} and~\cite[\S1.1]{Freed_Det}. 
\begin{defn} \label{defn:geofam}
A \emph{family of Riemannian manifolds over $B$} is 
\begin{enumerate}
\item a smooth fibration of oriented manifolds $\pi \colon M\to B$,
\item a Riemannian metric on the fibers, i.e., a metric on the vertical tangent bundle~$T(M/B)={\rm ker}(d\pi)\subset TM$, and
\item a projection $P\colon TM\to T(M/B)$. 
\end{enumerate}
\end{defn}
A family of Riemannian manifolds over $M$ can be endowed with a connection on $T(M/B)$ by first fixing an arbitrary metric $g_B$ on the base~$B$. The metric on~$T(M/B)$ from (2) above together with the horizontal lift of $g_B$ using $P$ from (3) gives a metric on $M$ that has a Levi-Civita connection. Then, \cite[Theorem 1.9]{BismutAS} shows that the projection of this Levi-Civita connection to $T(M/B)$ is independent of the choice of~$g_B$. Let $p_k(T(M/B))$ denote the relative Pontryagin forms associated with the connection on the vertical tangent bundle of a family of Riemannian manifolds constructed in the above way. Similarly, let $\hat{A}(M/B)$ and $\Wit(M/B)$ denote the differential form representatives of the~$\hat{A}$- and Witten classes of the family $\pi\colon M\to B$ using relative Pontryagin forms.

\begin{defn}
A \emph{rational string structure} on a family of Riemannian manifolds is a 3-form $H\in \Omega^3(M)$ such that $dH=p_1(T(M/B))$ is a trivialization of the relative 1st Pontryagin form.
\end{defn}

\subsection{Supermanifolds}\label{appen:super}

\begin{defn} A \emph{$k|l$-dimensional supermanifold} is a locally ringed space whose underlying topological space is Hausdorff and second countable, and whose structure sheaf is locally isomorphic to $C^\infty(U)\otimes_\C\Lambda^\bullet(\C^l)$ as a super algebra over~$\C$ for $U\subset \R^k$ an open submanifold. Supermanifolds and maps between them (as locally ringed spaces) form a category we denote by~${\sf SMfld}$. The set of morphisms between supermanifolds~$N$ and~$N'$ is denoted ${\sf SMfld}(N,N')$. \end{defn}

We follow the standard notation wherein~$C^\infty(N)$ denotes the global sections of the structure sheaf of a supermanifold~$N$, and refer to these sections as \emph{functions} on~$N$. The existence of partitions of unity imply that maps between supermanifolds $N\to N'$ are equivalent to the data of a $\Z/2$-graded algebra map~$C^\infty(N')\to C^\infty(N)$, i.e., supermanifolds are affine. 

\begin{rmk} The objects called supermanifolds in this paper are called $cs$-manifolds in~\cite{DM}. Indeed, the more common definition of a supermanifold has a structure sheaf defined over~$\R$. However, working with Wick-rotated field theories necessitates the category of $cs$-manifolds (see~\cite[Example 4.9.3, pg.~55]{strings1}), and the field theories we consider in this paper are Wick rotated. 
\end{rmk}

\begin{ex}
For natural numbers $n$ and $m$, the supermanifold~$\R^{n|m}$ is given by the manifold $\R^n$ equipped with the structure sheaf $C^\infty(U)\otimes_\C\Lambda^\bullet(\C^m)$ for $U\subset \R^n$. Following the notation above, $C^\infty(\R^{n|m})\cong C^\infty(\R^n)\otimes_\C\Lambda^\bullet(\C^m)$. Just as any $n$-dimensional manifold is locally isomorphic to $\R^n$, any $n|m$-dimensional supermanifold is locally isomorphic to~$\R^{n|m}$. 
\end{ex}

\begin{ex} \label{ex:mansup}Ordinary $n$-dimensional manifolds with their sheaf of complex-valued smooth functions define supermanifolds of dimension~$n|0$. This assignment gives a faithful embedding of manifolds into supermanifolds. 
\end{ex}

Given a supermanifold~$N$, let $N_{\rm red}$ denote the supermanifold gotten by taking the quotient of the structure sheaf by its nilpotent ideal; then $N_{\rm red}$ is an ordinary manifold that we have regarded as a supermanifold, and there is an evident map $N_{\rm red}\hookrightarrow N$. 

\begin{rmk} 
Ignoring gradings, the above allows one to think of the relationship between manifolds and supermanifolds as being analogous to the relationship between reduced schemes and non-reduced schemes. Indeed, the inclusion $N_{\rm red}\hookrightarrow N$ witnesses $N$ as a nilpotent thickening of its reduced manifold.
\end{rmk}
\begin{rmk}\label{rmk:real}
Note that the sheaf of functions on $N_{\rm red}$ has a real structure (i.e., a complex-conjugation map) but we do not demand that this structure extend to the structure sheaf of~$N$. In particular, the complex conjugate of a function on a supermanifold has no meaning in general. 
\end{rmk}

\begin{thm}[\cite{batchelor}] Any supermanifold~$N$ is isomorphic to $(N_{\rm red},\Gamma(\Lambda^\bullet E^\vee))$ for $E\to N_{\rm red}$ a complex vector bundle over a smooth manifold~$N_{\rm red}$. 
\end{thm}

The supermanifold $(N_{\rm red},\Gamma(\Lambda^\bullet E^\vee))$ is often denoted $\Pi E\cong N$, though we use this notation sparingly to avoid confusion with the parity reversal of a super vector bundle (see below). A map of vector bundles~$E\to E'$ defines a map between supermanifolds $\Pi E\to \Pi E'$, but the crucial fact (that makes supermanifolds interesting) is that there additional maps of supermanifolds $\Pi E\to \Pi E'$ that do not come from vector bundle maps. 

\begin{ex} The supermanifold $\R^{n|m}$ is $\Pi \underline{\C}^m$ for $\underline{\C}^m\to \R^n$ the trivial rank $m$ vector bundle on $\R^n$. \end{ex}

It is often convenient to study a supermanifold in terms of its functor of points. This uses the Yoneda embedding into presheaves, wherein a supermanifold $N$ defines the (set-valued) presheaf $S\mapsto {\sf SMfld}(S,N)$. The value of this presheaf on $S$ is called the \emph{$S$-points of~$N$} and denoted~$N(S)$. By the Yoneda lemma, morphisms between representable presheaves are in bijection with morphisms between their underlying supermanifolds. There are also interesting presheaves on supermanifolds that might fail to be representable. One example is the mapping presheaf $\Map(T,N)$ whose $S$-points are defined as
$$
\Map(T,N)(S):={\sf SMfld}(S\times T,N).
$$ 
This presheaf typically fails to be representable, just as mapping spaces between ordinary manifolds aren't usually (finite-dimensional) manifolds. 

\begin{ex}
The $S$-points of $\R^{n|m}$ are
\beq
&&\R^{n|m}(S)= \{x_1,\dots,x_n\in C^\infty(S)^\ev,\ \theta_1,\dots,\theta_m\in C^\infty(S)^\odd\mid (x_i)_{\rm red}=\overline{(x_i)}_{\rm red}\}\label{RnmSpot}
\eeq
where $(x_i)_{\rm red}$ denotes the restriction of a function $x_i$ to the reduced manifold $S_{\rm red}\hookrightarrow S$. We emphasize that the reality condition on the $x_i$ is only imposed on this restriction; see Remark~\ref{rmk:real}. 
\end{ex}

\begin{ex}\label{eq:CSpoints}
By the description above, we can describe $S$-points of $\C^2\cong \R^2$ as
\beq
\C(S)=\{z,w \in C^\infty(S)^{\ev}\mid z_{\rm red}=\overline{w}_{\rm red}\}\label{eq:C(S)}
\eeq
where $z_{\rm red}$ and $w_{\rm red}$ denote the restriction of $z$ and $w$ to the reduced manifold~$S_{\rm red}$. In a common abuse of notation, we write this pair of functions~$(z,w)$ as~$(z,\bar z)$, though we emphasize that $z$ and $\bar z$ are \emph{not} complex conjugates in~$C^\infty(S)$: indeed, this algebra need not have a real structure. In particular, when working with $S$-points of $\R^{2|1}$, we implicitly use the identification $\C\cong \R^2\subset \R^{2|1}$ and the above convention to write an $S$-point as $(z,\bar z,\theta)\in \R^{2|1}(S)$. Similarly, we write $S$-points of $\C^\times$ as $(\mu,\bar \mu)\in \C^\times(S)\subset \C(S)$. 
\end{ex}

\begin{ex}\label{ex:GL1}
Another example of a non-representable presheaf is $\GL(\C)$, the presheaf on supermanifolds whose value on~$S$ is the set of even invertible functions on~$S$; equivalently, this is the presheaf that assigns to a supermanifold~$S$ the automorphisms of the free rank one module over~$C^\infty(S)$. Multiplication of invertible functions promotes $\GL(\C)$ to a group object in presheaves. To see that $\GL(\C)$ is not representable as a supermanifold, observe that if it were representable its underlying manifold would be $\GL(\C)(\pt)\simeq C^\infty(\pt)^\times \simeq \C^\times$. On the other hand, there is an evident surjective homomorphism of group objects in presheaves
\beq
\C^\times\to \GL(\C)\qquad \C^\times (S)\ni (z,\bar z) \mapsto z\in (C^\infty(S)^\times)^{\ev}\label{eq:GL1map}
\eeq
where 
$$
\C^\times(S)\simeq \{z,\bar z \in (C^\infty(S)^\times)^\ev\mid z_{\rm red}=\overline{(\bar z_{\rm red})}\}
$$
uses the notation from the previous example. However,~\eqref{eq:GL1map} is not a monomorphism of presheaves, as $z\in (C^\infty(S)^\times)^\ev$ does not uniquely determine $\bar z\in (C^\infty(S)^\times)^\ev$. From this we see that $\GL(\C)$ is not representable. To the best of our knowledge, this argument is due to Lars Borutzky. There is a similarly defined non-representable presheaf $\GL(\C^{0|1})\simeq \GL(\C)$ with the same $S$-points as $\GL(\C)$, but viewed as automorphisms of the free rank one \emph{odd} module over~$C^\infty(S)$. 
\end{ex}

\begin{defn}
A \emph{vector bundle} over a supermanifold is a locally free module over its structure sheaf. We call elements of these modules \emph{sections} of the vector bundle. 
\end{defn}

\begin{ex} For a supermanifold $N$, consider the sheaf of modules of $\C$-linear derivations of its sheaf of functions. This is a vector bundle, namely the tangent bundle of $N$, denoted~$TN$. We observe that maps between supermanifolds~$f\colon N\to N'$ determine morphisms of tangent bundles~$df\colon TN\to TN'$. This allows us to import many standard ideas and theorems from ordinary manifolds, e.g., a map $f\colon N\to N'$ is a local isomorphism (or \emph{local embedding}) if and only if~$df$ induces an isomorphism on tangent spaces~\cite[\S3.4]{DM}. 
\end{ex}

We often abuse notation, writing a vector bundle over $N$ as $V\to N$ in spite of $V$ being defined as a sheaf of modules rather than a total space over~$N$. We use $\Pi V$ to denote the parity reversal of the module associated with~$V$. 

\begin{rmk} We caution that promoting a real vector bundle over an ordinary manifold to a vector bundle over the associated supermanifold necessitates complexification of the sheaf of sections, as it must be a module over the sheaf of complex-valued functions. We discuss this further in~\S\ref{sec:01MQ}. 
\end{rmk}

There is a well-developed theory of fiberwise integration on supermanifolds, namely fiberwise Berezinian integration; for example see~\cite[\S3]{Wittenintegration} and~\cite[pg.~94-96]{strings1}. In short, sections of a relative Berezinian line bundle can be integrated. We spell this out in the cases of interest. A fibration of supermanifolds is locally of the form~$S\times \R^{k|l}\to S$, and in this paper we will only need the case $l=1$ and $k=0,1,2$. Then a choice of coordinates $x_1,\dots,x_k,\theta$ on $\R^{k|1}$ defines a trivialization of the fiberwise Berezinian line (usually this trivialization is denoted~$[d\theta dx_1\cdots dx_k]$), giving a map
$$
\int\colon C^\infty_c(S\times \R^{k|1})\to C^\infty(S),
$$
defined on functions on $S\times \R^{k|1}$ with compact vertical support. 
On a function~$f$, this integral first projects to the Taylor coefficient of $\theta$ to get a function on~$S\times \R^k$, and applies the $C^\infty(S)$-linear extension of the usual integral on~$\R^k$ with respect to $dx_1\cdots dx_k$. To nail down the signs, for $a+\beta\theta\in C^\infty(S\times \R^k)[\theta]\cong C^\infty(S\times \R^{k|1})$, we take
$$
\int(a+\beta\theta)[d\theta dx]=\int \beta dx\in C^\infty(S).
$$
So the integral over~$\theta$ is determined by the formula~$\int \theta [d\theta]=1$. In particular, the integral over odd variables is a totally algebraic operation. We observe that total derivatives are sent to zero, e.g., for~$f=f(x,\theta)=f_0(x)+\theta f_1(x)$ we have
$$
\int (\partial_\theta f)[d\theta]=\int f_1 [d\theta]= 0. 
$$

\subsection{Lie groupoids and stacks}\label{appen:stacks}

 We give a brief introduction to super Lie groupoids and super stacks; good references for this material include~\cite{BlohmannStacks,Hollander_Stacks,Lerman} and (especially for our purposes) the appendix of~\cite{HKST}. 
 
 \begin{defn}
 A \emph{super Lie groupoid} is a groupoid internal to supermanifolds. 
 \end{defn}
 
In more detail, a super Lie groupoid~$\mathcal{G}$ has as data a pair of supermanifolds $\mathcal{G}_0$ (the objects) and $\mathcal{G}_1$ (the morphisms) together with maps
$$
s,t\colon \mathcal{G}_1\to \mathcal{G}_0\quad u\colon \mathcal{G}_0\to \mathcal{G}_1\quad c\colon \mathcal{G}_1\times_{\mathcal{G}_0}\mathcal{G}_1\to \mathcal{G}_1\quad (-)^{-1}\colon \mathcal{G}_1\to \mathcal{G}_1
$$
with the notation standing for source, target, unit, composition, and inversion, respectively. These data are required to satisfy the usual axioms for a groupoid, phrased in terms of commutative diagrams. For example, the diagrams
\begin{equation}
\begin{array}{c}
\begin{tikzpicture}
  \node (A) {$\mathcal{G}_1$};
  \node (B) [node distance= 2cm, right of=A] {$\mathcal{G}_1\times_{\mathcal{G}_0}\mathcal{G}_1$};
  \node (C) [node distance = 2cm, right of=B] {$\mathcal{G}_1$};
  \node (D) [node distance = 1.5cm, below of=A] {$\mathcal{G}_0$};
  \node (E) [node distance = 1.5cm, below of=B] {$\mathcal{G}_1$};
  \node (F) [node distance = 1.5cm, below of=C] {$\mathcal{G}_0$};
  \draw[->] (B) to node [above] {$p_1$} (A);
  \draw[->] (B) to node [above] {$p_2$} (C);
  \draw[->] (A) to node [left] {$t$} (D);
  \draw[->] (B) to node [right] {$c$} (E);
  \draw[->] (C) to node [right] {$s$} (F);
  \draw[->] (E) to node [below] {$t$} (D);
  \draw[->] (E) to node [below] {$s$} (F);
\end{tikzpicture}\end{array}\nonumber\qquad
\begin{array}{c}
\begin{tikzpicture}
  \node (A) {$\mathcal{G}_1\times_{\mathcal{G}_0}\mathcal{G}_1\times_{\mathcal{G}_0}\mathcal{G}_1$};
  \node (B) [node distance= 4cm, right of=A] {$\mathcal{G}_1\times_{\mathcal{G}_0}\mathcal{G}_1$};
  \node (C) [node distance = 1.5cm, below of=A] {$\mathcal{G}_1\times_{\mathcal{G}_0}\mathcal{G}_1$};
  \node (D) [node distance = 1.5cm, below of=B] {$\mathcal{G}_1$};
  \draw[->] (A) to node [above] {$c\times \id_{\mathcal{G}_1}$} (B);
  \draw[->] (A) to node [left] {$\id_{\mathcal{G}_1}\times c$} (C);
  \draw[->] (B) to node [right] {$c$} (D);
  \draw[->] (C) to node [below] {$c$} (D);
\end{tikzpicture}\end{array}\nonumber
\end{equation}
require that identity arrows act as the identity and that composition is associative. In all example we consider, the source map is a submersion so that the fibered products (e.g.,~$\mathcal{G}_1\times_{\mathcal{G}_0} \mathcal{G}_1$) automatically exist in supermanifolds. We will often drop the modifier ``super," simply referring to these as Lie groupoids. 

\begin{ex} Let $G$ be a super Lie group acting on a supermanifold~$N$. Then the \emph{quotient groupoid}, denoted $N\sq G$ has as object supermanifold $N:=(N\sq G)_0$ and morphism supermanifold $G\times N:=(N\sq G)_1$. The source map is the projection $G\times N\to N$ and target map is the action $G\times N\to N$. The identity map includes $N$ along the identity element of~$G$, composition is determined by multiplication in the group, and inversion is determined by inversion on the group. 
\end{ex}
 
 \begin{defn}
 A \emph{(smooth) functor} $F\colon \mathcal{G}\to \mathcal{H}$ between Lie groupoids $\mathcal{G}=\{\mathcal{G}_1\rightrightarrows \mathcal{G}_0\}$ and $\mathcal{H}=\{H_1\rightrightarrows H_0\}$ has as data a pair of maps $F_0\colon \mathcal{G}_0\to H_0$ and $F_1\colon \mathcal{G}_1\to H_1$, and we require that these data make the diagrams defining a functor between groupoids commute. A \emph{smooth natural transformation} between functors $F,G\colon\mathcal{G}\to \mathcal{H}$ is a smooth map $\mathcal{G}_0\to\mathcal{H}_1$ satisfying the usual axioms of a natural transformation. 
 \end{defn}

Lie groupoids, smooth functors, and smooth natural transformations form a strict 2-category, ${\sf LieGrpd}$. We caution that a fully faithful and essentially surjective map between Lie groupoids need not have a smooth inverse. A prototypical example comes from regarding a principal $G$-bundle~$P\to M$ as a groupoid $P\sq G$. Then the quotient map to~$M$ (regarded as a discrete Lie groupoid) is fully faithful and essentially surjective, but has a smooth inverse if and only if $P$ is a trivial $G$-bundle. One way to enhance our 2-category of Lie groupoids to include these inverses is to pass to smooth (super) stacks.
 
 \begin{defn} A \emph{cover} of a supermanifold $S$ is a local embedding $p\colon U\to S$ so that the smooth map on underlying manifolds $p_{\rm red}\colon U_{\rm red} \to S_{\rm red}$ is an ordinary open cover. 
 \end{defn}
 
 \begin{defn}
A \emph{prestack} $\X$ is a weak 2-functor from supermanifolds to groupoids, $\X\colon {\sf SMfld}^\op\to{\sf Grpd}$. A \emph{stack} is a prestack satisfying descent with respect to open covers of supermanifolds. 
\end{defn}

\begin{rmk} Fibered categories give an equivalent formulation of stacks, e.g., see~\cite[Definition~7.9]{HKST}. \end{rmk}

Natural transformations of 2-functors define morphisms between prestacks or stacks, which assemble into a groupoid we denote by~${\sf PreSt}(\X,\Y)$ or~${\sf St}(\X,\Y)$. Prestacks and stacks form bicategories, denoted ${\sf PreSt}$ and ${\sf St}$ respectively. 

\begin{ex} Any sheaf of sets on supermanifolds defines a stack by viewing sets as discrete categories. In particular, the sheaf $S\mapsto {\sf SMfld}(S,N)$ defined for any supermanifold~$N$ is a stack that we denote simply by~$N$, and the assignment $S\mapsto C^\infty(S)$ is a stack denoted~$C^\infty(-)$. \end{ex}

\begin{defn} A \emph{function} on a stack $\X$ is a morphism of stacks~$\X\to C^\infty(-)$. The collection of all such morphisms is denoted $C^\infty(\X)$, and has the structure of an algebra. 
\end{defn}

\begin{ex} 
Consider ${\sf Vect}\colon {\sf SMfld}^\op\to {\sf Grpd}$ that assigns to a supermanifold $S$ the groupoid of vector bundles and vector bundle isomorphisms over~$S$. Since vector bundles can be pulled back along maps $S\to S'$, ${\sf Vect}$ defines a prestack. Since vector bundles and vector bundle maps can be assembled from local data on an open cover, ${\sf Vect}$ is in fact a stack. 
\end{ex}

\begin{defn} A vector bundle on a stack $\X$ is a morphism of stacks $\X\to {\sf Vect}$. \end{defn}

Equivalently, a vector bundle $\mathcal{V}\to \X$ assigns to each $S$-point of $\X$ a vector bundle~$\mathcal{V}_S\to S$ over $S$, and for each 2-commuting triangle 
\beq
\begin{tikzpicture}[baseline=(basepoint)];
\node (A) at (0,0) {$S$};
\node (B) at (3,0) {$S'$};
\node (C) at (1.5,-1) {$\X$};
\draw[->] (A) to node [above] {$\varphi$} (B);
\draw[->] (A) to (C);
\draw[->] (B) to (C);
\path (0,-.75) coordinate (basepoint);
\end{tikzpicture}\nonumber
\eeq
an isomorphism of vector bundles $\varphi^*\mathcal{V}_{S'}\simeq \mathcal{V}_S$, and these isomorphisms are required to be compatible with compositions of base changes. 

We recall that there is a functor from prestacks to stacks called \emph{stackifification}, which is the left adjoint to the forgetful functor $u\colon {\sf St}\to {\sf PreSt}$ from stacks to prestacks, analogous to sheafification of presheaves. In particular, we have an equivalence of groupoids
\beq
{\sf PreSt}(\X_{\rm pre},u(\Y))\simeq {\sf St}(\X,\Y)\label{eq:leftad}
\eeq
where $\X$ is the stackification of a prestack $\X_{\rm pre}$. The above is quite useful in practice, allowing one to compute the groupoid of maps between a pair of stacks in terms of maps between associated prestacks. For example, the functions on~$\X$ coincide with the functions on $\X_{\rm pre}$, and the groupoid of vector bundles on $\X$ coincides with the groupoid of vector bundles on $\X_{\rm pre}$. 
 
There is a 2-functor ${\sf LieGrpd}\to{\sf PreSt}$ that sends a Lie groupoid to the prestack whose value on a supermanifold $S$ is the groupoid $\mathcal{G}(S):=\{\mathcal{G}_1(S)\rightrightarrows \mathcal{G}_0(S)\}$. 
 
 \begin{defn}\label{defn:grpdstack} The \emph{stack underlying} a Lie groupoid $\mathcal{G}$, denoted $[\mathcal{G}]$, is the stackification of the prestack associated to $\mathcal{G}$. If a stack $\X$ is equivalent to $[\mathcal{G}]$ for some Lie groupoid $\mathcal{G}$, then $\X$ is \emph{a geometric stack} and $\mathcal{G}$ is a \emph{groupoid presentation} of $\X$. 
 \end{defn}
 
 \begin{rmk} The stackification of the prestack associated to a Lie groupoid has an explicit description in terms of bibundles, e.g., see~\cite{BlohmannStacks,Lerman}. However, for the computations in this paper we only use formal properties of stackification. \end{rmk}

\begin{rmk}\label{rmk:localization}
Fully faithful and essentially surjective smooth functors between Lie groupoids are sent to equivalences of stacks under the map~$\mathcal{G}\mapsto [\mathcal{G}]$. Hence, one point of view on geometric stacks is as a localization of the category of Lie groupoids so that fully faithful and essentially surjective smooth functors become invertible. 
\end{rmk}

\begin{defn}
Given a diagram of stacks
\beq
\begin{tikzpicture}
  \node (A) {$\null$};
  \node (B) [node distance= 1.7cm, right of=A] {$\X$};
  \node (C) [node distance = 1.5cm, below of=A] {$\Y$};
  \node (D) [node distance = 1.5cm, below of=B] {$\mathcal{Z}$};
  \draw[->] (B) to node [right] {$F$} (D);
  \draw[->] (C) to node [below] {$F'$} (D);
\end{tikzpicture}\nonumber
\eeq
the \emph{2-pullback} $\X\times_\mathcal{Z} \Y$ is the stack whose objects over $S$ are
$$
\X\times_\mathcal{Z} \Y(S):=\{f,f',\eta\mid f\colon S\to \X, \ f'\colon S\to \Y, \ \eta\colon F\circ f\Rightarrow F'\circ f'\},
$$
and morphisms $(f,f',\eta)\to (g,g',\omega)$ in $\X\times_\mathcal{Z} \Y(S)$ are given by morphisms $\phi\colon f\to g$, $\phi'\colon f'\to g'$ so that $\omega\circ F(\phi)=F'(\phi')\circ \eta$. 
\end{defn}

If $\X$ is a geometric stack, there is a map $\mathcal{G}_0\to \X$ and an isomorphism between the morphisms of $\mathcal{G}$ and the 2-pullback of objects over themselves,~$\mathcal{G}_1\simeq \mathcal{G}_0\times_\X \mathcal{G}_0$. Hence, the groupoid presenting~$\X$ can be recovered from the map~$\mathcal{G}_0\to \X$. This point of view can be useful in practice, and is formalized as follows. 

\begin{defn}
An \emph{atlas} for a stack $\X$ is a map $p\colon U\to \X$ with source a supermanifold~$U$ so that for any other map $q\colon V\to \X$ with source a supermanifold $V$, the 2-fibered product $U\times_\X V$ is representable, and the map $U\times_\X V\to V$ is a submersion. 
\end{defn}

An atlas defines a groupoid presentation, $\{U\times_\X U\rightrightarrows U\}$, and so a stack is geometric if and only if it admits an atlas. 

\begin{ex} When $\X$ is geometric, we can compute its functions in terms of a groupoid presentation. Indeed,
$$
C^\infty(\X)\cong C^\infty([\mathcal{G}])\cong \{f\in C^\infty(\mathcal{G}_0)\mid s^*f=t^*f\}
$$
for $\mathcal{G}=\{\mathcal{G}_1\rightrightarrows \mathcal{G}_0\}$ the presenting groupoid with source and target maps $s$ and $t$, respectively. 
\end{ex}

\begin{ex} Recall the non-representable group object in presheaves $\GL(\C^{0|1})$ from Example~\ref{ex:GL1}. Let~$[\pt\sq \GL(\C^{0|1})]$ denote the stack that classifies odd line bundles, i.e., the stackification of the prestack that sends a supermanifold $S$ to the groupoid with a single object (the trivial odd line bundle on~$S$) and morphisms $\GL(\C^{0|1})(S)$. The notation~$[\pt\sq \GL(\C^{0|1})]$ is slightly misleading in that $\GL(\C^{0|1})$ is not a representable presheaf and so $\pt\sq \GL(\C^{0|1})$ is not a Lie groupoid. However, it is a groupoid object in presheaves, which is enough to define a stack in the same fashion as Definition~\ref{defn:grpdstack}, and the result is indeed the stack classifying odd line bundles on supermanifolds. Furthermore, from the discussion in Example~\ref{ex:GL1}, ${\rm GL}(\C^{0|1})$ receives a homomorphism $\C^\times\to {\rm GL}(\C^{0|1})$. This defines maps of stacks
$$
[\pt\sq \R^\times]\to [\pt\sq \C^\times]\to [\pt\sq \GL(\C^{0|1})],
$$
which determine canonical odd line bundles over the stacks $[\pt\sq \R^\times]$ and $[\pt\sq \C^\times]$ by pullback. We call these line bundles the \emph{canonical odd line bundle} over the stacks $[\pt\sq \R^\times]$ and~$[\pt\sq \C^\times]$. 
\end{ex}

\begin{ex} 
Consider the groupoid ${\sf St}([N\sq G],{\sf Vect})$ of vector bundles on $[N\sq G]$. Using the adjunction~\eqref{eq:leftad}, we see that objects consist of a vector bundle on~$N$ with a $G$-action, i.e., a $G$-equivariant vector bundle on~$N$. In particular, a representation $\rho\colon G\to {\rm End}(V)$ defines a vector bundle $V_\rho$ on~$[N\sq G]$. When $V$ is 1-dimensional, sections of $V_\rho$ and $\Pi V_\rho$ are (e.g., see~\cite[Corollary 7.17]{HKST})
\beq&&\begin{array}{ccc}
\Gamma(N\sq G, V_\rho)&\cong& \{f\in C^\infty(N) \mid \mu^*(f)=p_1^*(\rho)\cdot p_2^*(f) \in C^\infty(G\times N)\},\\
\Gamma(N\sq G,\Pi V_\rho)&\cong& \{f\in C^\infty(N)^{\odd}\mid \mu^*(f)=p_1^*(\rho)\cdot p_2^*(f)\in C^\infty(G\times N)\}\end{array}
\label{compsections}
\eeq
where $p_1\colon G\times N\to G,$ $p_2\colon G\times N\to N$ are the projection maps, and $\mu\colon G\times N\to N$ is the action map. Note that for the trivial representation sections are precisely the $G$-invariant functions on~$N$.
\end{ex}

\subsection{Modular forms and Eisenstein series}\label{appen:modform} The following brief overview (and much more) can be found in~\cite[\S1-2]{ZagierMF}. For simplicity, in this section we lapse into complex coordinate notation, e.g., denoting a point in $\C^\times$ by $\mu$ rather than $(\mu,\bar\mu)$.

A 2-dimensional \emph{based, oriented lattice} is a monomorphism $\ell\colon \Z^2\to \R^2\cong \C$ such that the ratio of the image of the generators~$\ell_1\colon \{(1,0)\}\to \C$ and~$\ell_2\colon \{(0,1)\}\to \C$ defines a point $\frac{\ell_2}{\ell_1}\in \mathbb{H}\subset \C$ in the upper half plane. Let $\Lat$ denote the smooth manifold of these lattices; we have an evident diffeomorphism $\Lat\cong \mathbb{H}\times \C^\times$ that sends a pair of generators~$\ell_1,\ell_2$ to $(\ell_2/\ell_1,\ell_1)\in \mathbb{H}\times \C^\times$. There is an action of $\C^\times\times \SL_2(\Z)$ on~$\Lat$ by
\beq
\begin{array}{c}
\left(\mu,\left[\begin{array}{cc} a& b \\ c & d\end{array}\right],\ell_1,\ell_2\right)\mapsto (\mu^2(a\ell_1+b\ell_2),\mu^2(c\ell_1+d\ell_2)),\\ \mu\in \C^\times, \ \left[\begin{array}{cc} a& b \\ c & d\end{array}\right]\in \SL_2(\Z). \end{array}\label{eq:scaling}
\eeq

\begin{defn}\label{defn:wMF} \emph{Weak modular forms of weight $n/2$} are holomorphic functions $h$ on the space of based, oriented lattices $\Lat$ that are ${\rm SL}_2(\Z)$-invariant and have the property that $h(\mu\cdot \ell)=\mu^{-n}h(\ell)$ for $\mu\in \C^\times$. \end{defn}

There is an equivalent description of weak modular forms of weight~$k$ as holomorphic functions on the upper half plane $\mathbb{H}\subset \C$ satisfying
$$
f\left(\frac{a\tau+b}{c\tau+d}\right)=(c\tau+d)^kf(\tau)\quad \tau\in \mathbb{H}.
$$
A modular form in the sense of Definition~\ref{defn:wMF} gives a function on the upper half plane by restriction to lattices of the form $(\ell_1,\ell_2)=(1,\tau)$. Conversely, given a function on the upper half plane satisfying the modularity property above, we can extend it to a modular form as a function on based lattices using the formula $F(\ell_1,\ell_2)=(\ell_1)^{-k}f(\ell_2/\ell_1).$ As a function on the upper half plane, a modular form is invariant under the subgroup $\Z\hookrightarrow \SL_2(\Z)$ given by
$$
n\mapsto \left[\begin{array}{cc} 1 & n \\ 0 & 1\end{array}\right],
$$
and so defines a function on $\mathbb{H}/\Z$. The map $\tau\mapsto q=\exp(2\pi i\tau)$ for $\tau\in\mathbb{H}$ gives an analytic isomorphism between $\mathbb{H}/\Z$ and the punctured disk. Writing a weak modular form as a power series in the variables~$q$ and $q^{-1}$ is called \emph{$q$-expansion}. 

\begin{defn} \label{defn:MF}
A \emph{(weakly holomorphic) modular form} is a weak modular form whose $q$-expansion has only finitely many negative powers of~$q$. Taking products of holomorphic functions gives a graded ring, denoted ${\rm MF}$ whose degree $n$ piece, denoted ${\rm MF}_n$, are the weight $n/2$ modular forms. Define~$\MF^n:=\MF_{-n}$. 
\end{defn}


\bibliographystyle{amsalpha}
\bibliography{references}
\end{document}